\documentclass[letterpaper,  ngerman, openright, BCOR1cm, bibtotoc, 11pt, dvips]{scrbook}
\usepackage{diplomarbeit}

\usepackage[ngerman]{babel} %

\newcommand{\DatumInKlammern}{} %
\newcommand{\DatumErklaerung}{13. Juli 2007}
\newcommand{\DatumTitel}{Juli 2007}

\newcommand{\TitelrueckseiteOben}{}
\newcommand{\TitelrueckseiteUnten}{}

\begin{document}
\title{%
Torische Ideale von Flusspolytopen}
\author{Diplomarbeit von \\Matthias Lenz}
\publishers{Betreut von \\Dr.~Christian Haase\\ Fachbereich Mathematik und
Informatik\\Freie Universit"at Berlin
}
\date{\DatumTitel
\\[5cm]
	\scalebox{0.5}{
	\begin{picture}(0,0)%
\includegraphics{PullingWuerfel3sw2.pstex}%
\end{picture}%
\setlength{\unitlength}{4144sp}%
\begingroup\makeatletter\ifx\SetFigFont\undefined%
\gdef\SetFigFont#1#2#3#4#5{%
  \reset@font\fontsize{#1}{#2pt}%
  \fontfamily{#3}\fontseries{#4}\fontshape{#5}%
  \selectfont}%
\fi\endgroup%
\begin{picture}(4206,4206)(26473,-589)
\end{picture}%

	}
}

\uppertitleback{\TitelrueckseiteOben}%
\lowertitleback{\TitelrueckseiteUnten}%

\maketitle%

\setcounter{tocdepth}{2}
\tableofcontents

\setcounter{chapter}{-1}

\chapter{Einleitung \DatumInKlammern}

\section{Einführung}

Geschafft! Nach knapp fünf Jahren Studium bin ich nun an dessen Ende angelangt. 
Am meisten Spaß haben mir während dieser Zeit alle die Sachen gemacht, die mit Algebra und Kombinatorik zu tun hatten. Das Thema dieser Arbeit liegt im Grenzbereich dieser beiden Gebiete, wobei der Schwerpunkt auf der Kombinatorik liegt.

Zum Inhalt:
Im 1. Kapitel lernen wir die Begriffe und Objekte aus der diskreten Geometrie kennen, die wir in den folgenden Kapiteln benötigen. Zwischen den Klassen der wichtigsten dieser Objekte gilt folgende Beziehung:
\\[3pt]
$
\hspace*{0.7cm}\{\text{\emph{Polytope}}\} \supseteq \{\text{\emph{Gitterpolytope}}\} 
\supseteq
\{\text{\emph{Flusspolytope}}\} \supseteq\{\text{\emph{Transportpolytope}}\}
$\\[3pt]
Polytope sind konvexe, beschränkte Mengen im $\R^n$. Gitterpolytope sind solche mit ganzzahligen Ecken. Die letzten beiden Klassen erhält man als Lösungsmengen von Problemen aus der Graphentheorie (Fluss- und Transportprobleme).

In Kapitel 2 werden wir sehen, wie man Gitterpolytope in kleinere Gitterpolytope unterteilt, bis man schließlich eine Triangulierung, d.\,h. eine Zerlegung in Simplexe, erhält. Genauer untersuchen werden wir \emph{Pullingtriangulierungen} und \emph{reguläre Triangulierungen}.  

Die algebraischen Grundlagen werden in Kapitel 3 vermittelt. Dort wird definiert, wie man aus einer Punktmenge $\A\subseteq\Z^n$ (z.\,B. der Menge der Gitterpunkte eines Polytops) ein \emph{torisches Ideal} erhält. Außerdem beweisen wir eine Aussage über den Zusammenhang zwischen Gröbnerbasen von torischen Idealen auf der einen Seite  und regulären unimodularen Triangulierungen der zugehörigen Punktmenge auf der anderen Seite.  

In Kapitel 4 gehen wir folgenden Fragen nach: In welchem Grad sind die torischen 
Ideale von Transport- und Flusspolytopen erzeugt? Welche Gradschranken kann man für Gröbnerbasen angeben? Gibt es bessere Schranken, wenn man sich auf glatte Transportpolytope beschränkt?
Insbesondere werden wir beweisen, dass torische Ideale von Flusspolytopen alle im Grad drei erzeugt sind.

In Kapitel 1-3 werden (bis auf einen Teil von Abschnitt 
\ref{section:TransportKombinatorischGlatt}) bekannte Fakten wiedergegeben, teilweise für unseren Bedarf modifiziert. 
Bei Kapitel 4 handelt es sich bis auf Abschnitt \ref{section:Zellunterteilungsmethode}
um neue Erkenntnisse.
\smallskip

Die hier untersuchten Fragestellungen haben Verbindungen zu verschiedenen Gebieten der Mathematik. Transport- und Flussprobleme und die zugehörigen Polytope treten an verschiedenen Stellen in der kombinatorischen Optimierung und in vielen praktischen Anwendungen auf.
Transportpolytope kommen in der Statistik unter der Bezeichnung Kontingenztabellen vor.
Zu ihrer Untersuchung werden dort auch die zugehörigen torischen Ideale betrachtet.
In der algebraischen Geometrie interessiert man sich für glatte Polytope. 
Eine große Beispielmenge dafür sind die glatten  Flusspolytope.  

Gradschranken für Erzeugendensysteme und Gröbnerbasen bzw. Erzeugendensysteme und Gröbnerbasen in niedrigem Grad zu kennen ist u.\,a. deshalb interessant, weil man dann weiß, dass sich konkrete Berechnungen mit den torischen Idealen schnell durchführen lassen.

\smallskip
Zum Lesen (und Verstehen) dieser Arbeit werden nicht viele Vorkenntnisse benötigt. Der Leser sollte lediglich wissen, was ein Graph ist und die Grundlagen der (linearen) Algebra beherrschen.
 
  \medskip
  An dieser Stelle möchte ich ganz herzlich meinem Betreuer Christian Haase danken, der mich für dieses Thema begeistert hat und jederzeit für Fragen und Diskussionen zur Verfügung stand. Weiterhin danken möchte ich Martin Götze, der alle meine Fragen zu \LaTeX\ und Perl beantworten konnte, sowie René Birkner, der diese Arbeit Korrektur gelesen hat. 
\clearpage
\section{Notation}
In diesem Abschnitt wird kurz die verwendete Notation beschrieben.

Die \emph{Potenzmenge} einer  Menge $M$ bezeichnen wir mit $\pot(M)$.
Mit $[n]$ bezeichnen wir die Menge $\{1,\ldots,n\}$ und die \emph{natürlichen Zahlen} sind die Menge $\N:=\{0,1,2,3,\ldots\}$.

Für eine Menge $M$ bezeichnet $M^{m\times n}$ die Menge der $(m\times n)$-Matrizen mit Einträgen aus $M$.
Die $(n\times n)$-Einheitsmatrix bezeichnen wir mit $I_n$ oder $I$, falls die Dimension klar ist.
Matrizen werden in der Regel mit lateinischen Großbuchstaben bezeichnet.

Sei A eine $(m \times n)$-Matrix. Den Eintrag in der $i$-ten Zeile und $j$-ten Spalte bezeichnen wir mit $a_{ij}$.
Die transponierte Matrix bezeichnen wir mit $A^T$.

Sei $M$ eine Menge. Punkte bzw. Vektoren $\vek{m}\in M^n$ werden geschrieben als\linebreak[3] $\vek{m}=(m_1,\ldots, m_n)$. $\vek{m_i}$ bezeichnet also einen Vektor, wohingegen $m_i$ die $i$-te Komponente eines Vektors $\vek{m}$ bezeichnet.
 Wir schreiben $\vek{a}>\vek{b}$ genau dann, wenn $a_i>b_i$ für alle $i$ gilt. Analog definieren wir
$\vek{a}<\vek{b}$, $\vek{a}\ge\vek{b}$ und $\vek{a}\le\vek{b}$.
$\abs{\a}_1:=\sum_i a_i$ bezeichnet wie üblich die $1$-Norm eines Vektors.

Zu einem $K$-Vektorraum $V$ bezeichnet $V^*:=\{\varphi : V\to K\,|\, \text{$\varphi$ ist linear}\}$ den Dualraum. $(\Z^n)^*$ bezeichnet entsprechend die $\Z$-linearen Abbildungen von $\Z^n$ nach $\Z$. Alle vorkommenden Vektorräume werden Unterräume von $\R^n$ sein.

\smallskip

\emph{Ungerichtete Graphen} auf einer (endlichen) Knotenmenge $V$ mit Kantenmenge $E\subseteq  \binom{V}{2}$ werden wie üblich geschrieben als $G=(V,E)$. \emph{Gerichtete Graphen} auf der Knotenmenge $V$ mit Kantenmenge
$\vec{E}\subseteq V\times V$ werden geschrieben als $\vec{G}=(V,\vec{E})$. 
In ungerichteten Graphen wird die Menge der zu $v\in V$ inzidenten Kanten bezeichnet mit $\delta(v):=\{\{a,b\}\in E \,|\, v=a \text{ oder } v=b\}$. Im gerichteten Fall definieren wir
für einen Knoten $v$ die Mengen $\din(v):=\{(a,b)\in \vec{E}\,|\, v=b\}$ und $\dout(v):=\{(a,b)\in \vec{E}\,|\,$ $ v=a\}$ sowie für eine Kante $(u,v)$ die Mengen  
$\din((u,v))=v$ und $\dout((u,v))=u$.

Mit $K_n$ bezeichnen wir den ungerichteten vollständigen Graphen auf $n$ Ecken. %
 Mit $\vec{K}_{m,n}$ bezeichnen wir den gerichteten vollständig bipartiten Graphen, bei dem in der 1. Farbklasse $m$ Knoten und in der 2. Farbklasse $n$ Knoten sind und alle Kanten von der ersten Farbklasse in die zweite Farbklasse zeigen. 

Für einen gerichteten Graphen mit (geordneter) Knotenmenge $V=\{v_1,\ldots,v_n\}$ und
Kantenmenge $E=\{e_1,\ldots,e_m\}$ ist die \emph{Indzidenzmatrix} die ($n\times m$)-Matrix $A$ mit 
\begin{equation*}
 a_{ij}=\begin{cases}
         	-1 & \text{für } v_i=\dout(e_j)\\
		+1 & \text{für } v_i=\din(e_j)\\
		 \;\;\:0 & \text{sonst}
        \end{cases}\quad.
\end{equation*}

\chapter{Gitterpolytope \DatumInKlammern}
In diesem Kapitel besch\"aftigen wir uns mit Objekten aus der diskreten Geometrie. Wir definieren zun\"achst Kegel, Polytope und Gitterpolytope. Dann widmen wir uns den Transport- und Flusspolytopen.

\section{Polyeder, Polytope und Kegel} 

In diesem Abschnitt werden kurz die grundlegenden Begriffe und wichtige Sätze der Polytoptheorie erklärt. Weitergehende Informationen und Beweise findet man beispielsweise in 
\cite[Kapitel 0-2]{zieglerPolytopes}.

\begin{Definition}[Unterräume, affine und konvexe Hülle, Kegel]
$\,$\\
Sei $V=\{\vek{v_1},\ldots,\vek{v_l}\} \subseteq \R^n$.
Dann definieren wir
 \begin{itemize}
  \item den von $V$ erzeugten Unterraum $\lin(V):=\{\sum_{i=1}^l \l_i \vek{v_i} \,|\, \l_i\in\R\}$,
  \item die \emph{affine Hülle} von $V$ bzw. den von $V$ erzeugten \emph{affinen Raum}\\ $\aff(V):=\{\sum_{i=1}^l \l_i \vek{v_i} \,|\, \l_i\in\R, \sum_{i=1}^l \l_i=1\}$,
   \index{affine Hülle}
  \item den von $V$ aufgespannten \emph{Kegel}
   $\cone(V):=\{\sum_{i=1}^l \l_i \vek{v_i} \,|\, \l_i\in\R, \l_i\ge 0\}$ und
  \item die \emph{konvexe Hülle}  $\conv(V):=\{\sum_{i=1}^l \l_i \vek{v_i} \,|\, \l_i\in\R, \l_i\ge 0, \sum_{i=1}^l \l_i=1\}$ von $V$.
  \end{itemize}
\end{Definition}

Eine Menge $U\subseteq\R^n$ heißt Unterraum\,/\:affin\,/\:konvex, wenn eine Menge $V\subseteq\R^n$ existiert, sodass $U=\lin(V)$\,/\;$U=\aff(V)$\,/\;$U=\conv(V)$.

Ein Menge $\s\subseteq\R^n$ heißt Kegel, wenn eine endliche Menge $V\subseteq \Q^n$ existiert 
mit $\s=\cone(V)$.\footnote{
  Bei uns sind also alle Kegel endlich erzeugt und rational.}

Seien $A_1, A_2\subseteq \R^n$  zwei affine Räume. $\varphi : A_1 \to A_2$ heißt \emph{affine Abbildung}, wenn 
für $\vek{a_1},\ldots, \vek{a_k} \in A_1$ und $\l_1,\ldots,\l_k\in\R$ mit  $\sum_i\l_i=1$ gilt:
$\varphi(\sum_i \l_i \vek{a_i})= \sum_i \l_i\varphi(  \vek{a_i})$.

Eine Menge $\{\vek{v_1},\ldots, \vek{v_l}\}\subseteq \R^n$ heißt \emph{affin abhängig}, wenn 
\index{affin abhängig!geometrisch}
ein Vektor $\lambda\in\R^l\setminus \{\vek{0}\}$ existiert mit $\sum_{i=1}^l \l_i\vek{v_i}=0$ und $\sum_{i=1}^l \l_i=0$. 
Andernfalls heißt $V$ \emph{affin unabhängig}.
\index{affin unabhängig!geometrisch}

Das \emph{relativ Innere} einer konvexen Menge $P\subseteq \R^n$, bezeichnet mit $\relint(P)$, ist definiert als das Innere (im topologischen Sinn) von $P$, aufgefasst als (topologischer) Unterraum von $\aff(P)$.

Die \emph{Dimension} einer konvexen Menge $P$ ist definiert als die Dimension der affinen Hülle.

Sei $\varphi\in(\R^n)^*$ und $c\in\R$\hspace{0.2pt}.\footnote{
Wer will, kann sich ein Skalarprodukt wählen und dann $(\R^n)^*$
und $\R^n$ identifizieren. Wir werden dies gelegentlich tun, ohne es explizit zu erwähnen, z.\,B. um in den Abbildungen ein Polytop und seine Normalenkegel in das gleiche Bild zeichnen zu können und um in späteren Kapiteln die Notation zu vereinfachen.}
Dann erhalten wir den \emph{positiven Halbraum} $H^+_{\varphi,c}=H^+:=
\{\x\in\R^n\,|\,\varphi(\x)\ge c\}$, den \emph{negativen Halbraum} $H^-_{\varphi,c}=H^-:=
\{\x\in\R^n\,|\,$ $\varphi(\x)\le c\}$, sowie die \emph{affine Hyperebene}
$H_{\varphi,c}=H:=H^+\cap H^-$.

$P\subseteq \R^n$ heißt \emph{Polyeder}, wenn $P$ Schnitt von endlich vielen Halbräumen ist.
$P$ heißt \emph{Polytop}, wenn $P$ die konvexe Hülle einer endlichen Menge $V$ ist.
Dies ist äquivalent dazu, dass $P$ ein kompaktes Polyeder  ist. Alle Kegel sind Polyeder.

Für zwei Polyeder $P$ und $Q$ definieren wir die \emph{Minkowskisumme} $P+Q:=
\{ p +q \,|\,$ $ p\in P,\, q\in Q \}$. Für ein Polyeder $P$ und eine Zahl $c\in\R_{\ge 0}$ definieren wir $c\cdot P:=\{c\cdot\v\,|\,\v\in P\}$.

Ein Kegel heißt \emph{spitz}, wenn er keinen von $\{\vek{0}\}$ verschiedenen Unterraum enthält.

Sei $P\subseteq\R^n$ ein Polytop und $F\subseteq P$. $F$ heißt \emph{Seite} von $P$, geschrieben $F\seite P$, genau dann, wenn $\varphi \in(\R^n)^*$ und $c\in \R$ existieren, sodass  $F=H_{\varphi,c}\cap P$ und $P\subseteq H^+_{\varphi,c}$, d.\,h. 
$F=\{\vek{v}\in P \,|\, \varphi(\v)=\min_{\vek{x}\in P}\varphi(\vek{x})\}$.
 $\varphi$ heißt dann (innerer) \emph{Normalenvektor} an $F$.

Seiten eines Polytops sind wieder Polytope. Seiten von Kodimension eins heißen \emph{Facetten}, nulldimensionale Seiten heißen \emph{Ecken}. 
$F\seite P$ nennen wir \emph{echte Seite}, wenn $\emptyset\not=F\not=P$ gilt.
Die Menge der Ecken von $P$ bezeichnen wir mit $\eck(P)$.  Für jedes Polytop $P$ gilt $\conv(\eck(P))=P$.

Für Seiten der Form $S=\{ \x \in P \,|\, x_i \ge c\}$ für ein $i\in[n]$ verwenden wir die abkürzende Schreibweise $[ x_i \ge c]$.

Die Menge der Seiten eines Polytops ist mittels \glqq $\subseteq$\grqq\ halbgeordnet und bildet einen Verband. Aus $G\seite F\seite P$ folgt also $G\seite P$. 
Jede echte Seite $G$ eines Polytops lässt sich schreiben als Schnitt der Facetten, die $G$ enthalten.
Jede Seite von Kodimension zwei lässt sich schreiben als Schnitt von zwei Facetten.

Ein $d$-dimensionaler \emph{Simplex} ist ein Polytop, das sich als konvexe Hülle von $d+1$ affin unabhängigen Vektoren schreiben lässt.

Ein $d$-dimensionales Polytop  heißt \emph{einfach}, wenn
alle Ecken in genau $d$-Facetten enthalten sind (Beispiel: Würfel).
Ein $d$-dimensionales Polytop heißt \emph{simplizial}, wenn alle Facetten $(d-1)$-Simplexe sind (Beispiel: Oktaeder/Kreuzpolytop). Die Polytope, die sowohl einfach als auch simplizial sind, sind genau die Simplexe.

\section{Gitter und Gitterpolytope}
In diesem Abschnitt widmen wir uns einer speziellen Klasse von Polytopen: Den \emph{Gitterpolytopen}. Das sind Polytope, deren Ecken alle auf einem Gitter (bei uns $\Z^n$) liegen. 

Sei $B:=\{\vek{b_1},\ldots,\vek{b_d}\}\subseteq \R^n$ eine linear unabhängige Menge. Die Menge 
$\Lambda:= \{ \sum_i \lambda_i \vek{b_i} \,|\, \l_i\in\Z \}$ bezeichnen wir dann als $d$-dimensionales \emph{Gitter} und $B$ als \emph{Gitterbasis} von $\Lambda$. 
Gitter sind also nichts anderes als endlich erzeugte Untergruppen von $\R^n$.
Ein wichtiges Beispiel für Gitter ist $\Lambda=\Z^n$. 
Aus der Cramerschen Regel folgt, dass eine Menge $\{\vek{b_1},\ldots,\vek{b_n}\}\subseteq\Z^n$ genau dann Gitterbasis von $\Z^n$ ist, wenn $\abs{\det(\vek{b_1},\ldots,\vek{b_n})}=1$ gilt.

Ein \emph {Gitterhomomorphismus} zwischen zwei Gittern
$\Lambda\subseteq \R^n$ und  $\Lambda'\subseteq \R^{n'}$
ist eine Abbildung $\psi : \Lambda \to \Lambda'$, die sich zu einer linearen Abbildung
 $\phi: \lin(\Lambda)\to\lin(\Lambda')$ fortsetzen lässt.
Ein \emph{Gitterisomorphismus} ist ein bijektiver Gitterhomomorphismus.

Es sei $\SL_n(\Z):=\left\{M\in\Z^{n\times n} |\, \abs{\det M}=1\right\}$ die Menge der \emph{unimodularen Matrizen}. Ein  Gitter mit Basis
$\{\vek{b_1},\ldots,\vek{b_d}\}$ heißt \emph{unimodular}, wenn es eine Matrix $M\in\SL_n(\Z)$ gibt, mit $M \vek{b_i}=\vek{e_i}$ für $i=1,\ldots, d$. $\Z^n$ ist also unimodular.

Ein \emph{affines Gitter} ist ein um einen Vektor $\vek{v}\in\R^n$ verschobenes Gitter. 
Ein affines Gitter $\vek{v}+\Lambda$ (für ein Gitter $\Lambda$) heißt \emph {unimodular}, wenn $\Lambda$ unimodular ist.

Eine \emph{affine Gitterabbildung} zwischen den beiden affinen Gittern $\Lambda+\vek{v}\subseteq \R^n$ und  $\Lambda'+\vek{v'}\subseteq \R^{m}$ ist eine Abbildung $\psi : \Lambda+\vek{v} \to \Lambda'+\vek{v'}$, die sich zu einer affinen Abbildung $\phi : \aff(\Lambda+\vek{v}) \to \aff(\Lambda'+\vek{v'})$ fortsetzen lässt. Ein \emph{affiner Gitterisomorphismus}
ist eine bijektive affine Gitterabbildung.
Eine affine Abbildung $\varphi : \R^n\to \R^n$ heißt \emph{unimodular}, wenn 
$M\in\SL_n(\Z)$ und $\vek{b}\in\R^n$ existieren, sodass
$\varphi(\vek{v})=M\vek{v}+\vek{b}$ für alle $\v\in\R^n$ gilt.

\begin{Definition}[Gitterpolytop]
  Ein Polytop $P\subset \R^n$ heißt \emph{Gitterpolytop} bezüglich des Gitters $\Lambda\subseteq \R^n$, wenn $\eck(P)\subseteq \Lambda$ gilt.
\end{Definition}
Im weiteren Verlauf werden unsere Gitterpolytope stets Gitterpolytope bezüglich des Gitters $\Lambda=\Z^n$ sein. Das sind dann also genau die Polytope, bei denen alle Ecken ganzzahlige Koordinaten haben.

Die Menge der Gitterpunkte eines gegebenen Polytopes $P\subseteq \R^n$ wird im weiteren Verlauf dieser Arbeit eine wichtige Rolle spielen. Wir bezeichnen sie meist mit $\A_P=\A:=P\cap \Z^n$.

\begin{Definition}[Gitteräquivalenz von Polytopen und Kegeln]
$\,$\\
Seien $P\subseteq \R^{n}$ und $P'\subseteq \R^{m}$ Gitterpolytope. 
$P$ und $P'$ heißen \emph{gitteräquivalent}, wenn es eine affine Abbildung
$\varphi : \R^n\to \R^m$ gibt mit $\varphi\vert_{P}$ bildet $P$ bijektiv auf $P'$ ab
und $\varphi\vert_{\Z^n\cap \aff(P)}$ bildet ab nach $\Z^m\cap \aff(P')$ und ist affiner Gitterisomorphismus.

Analog heißen zwei Kegel $\s\subseteq \R^{n}$ und $\s'\subseteq \R^{m}$  \emph{gitteräquivalent}, wenn es eine lineare Abbildung
$\varphi : \R^n\to \R^m$ gibt mit $\varphi\vert_{\sigma}$ bildet $\sigma$ bijektiv auf $\sigma'$ ab
und $\varphi\vert_{\Z^n\cap \lin(\sigma)}$ bildet ab nach $\Z^m\cap \lin(\sigma')$ und ist  Gitterisomorphismus.
\end{Definition}

Ein \emph{simplizialer Kegel} ist ein Kegel, der von linear unabhängigen
Vektoren erzeugt ist.
$\v\in \Z^d$ heißt \emph{primitiv}, wenn $\{\l_i \v \,|\, 0 <\l_i< 1\}\cap \Z^d=\emptyset$ gilt. 
Für jeden Kegel ist es möglich primitive Erzeuger anzugeben: Man wähle ein beliebiges Erzeugendensystem und multipliziere dann jeden Vektor mit dem Kehrwert des größten gemeinsamen Teilers seiner Komponenten. Für einen $d$-dimensionalen simplizialen Kegel gibt es ein eindeutiges Erzeugendensystem mit $d$ primitiven Elementen.

Wir definieren nun auf kombinatorische Art das Volumen von Simplexen und simplizialen Kegeln: 
\begin{Definition}[Normalisiertes Volumen]
Sei $\s=\conv(\vek{v_0},\ldots, \vek{v_d})\subseteq \R^n$ ein $d$-dimesionaler Simplex und
$E:=\{ \sum_{i=0}^d \l_i (\vek{v_i}-\vek{v_0}) \,|\, 0\le \l_i < 1\}$ das von den Ecken von $\s$ aufgespannte halboffene Parallelepiped (nachdem $\vek{v_0}$ in den Ursprung verschoben wurde).

Wir definieren dann das \emph{normalisierte Volumen} von $\s$ als die Anzahl der Gitterpunkte in $E$:
\begin{equation*}
\vol(\s):= \abs{E\cap \Z^n}%
\end{equation*}

Sei $\s=\cone(\vek{v_1},\ldots, \vek{v_d})\subseteq \R^n$ ein simplizialer Kegel und seien 
$\vek{v_1},\ldots, \vek{v_d}$ primitiv. Dann definieren wir
das \emph{normalisierte Volumen} $\vol(\s)$ von $\s$ folgendermaßen:
\begin{equation*}
\vol(\s):= \vol(\conv(\vek{0},\vek{v_1},\ldots, \vek{v_d}))
\end{equation*}
\end{Definition}
Diese Definition des Volumens eines Simplex ist unabhängig davon, welcher Vektor den Index $0$ hat. Dies folgt beispielsweise aus dem Volumenlemma auf Seite 
\pageref{volumenlemma}.

\begin{Bemerkung}[Äquivalente Definitionen des normalisierten Volumens]
Für einen volldimensionalen simplizialen Kegel $\s\subseteq \R^n$   mit primitiven Erzeugern \linebreak[3]$\{\vek{v_1}, \ldots,\vek{v_n}\}$  gilt:
\begin{align*}
\vol(\s)&= \left|\det(\vek{v_1},\ldots, \vek{v_n})\right|  \\
        &= n! \cdot L(\conv(\vek{0},\vek{v_1},\ldots,\vek{v_n})) \text{ (wobei mit $L$ das Lebesguemaß bezeichnet wird)} \\
        &= \abs{\quot{$\Z^n$}{$U$}} \text{ (wobei $U$ das von $\vek{v_1},\ldots, \vek{v_n}$ erzeugte Untergitter bezeichnet)} 
\end{align*}

\end{Bemerkung}

\begin{Definition}[Unimodulare Kegel und Simplexe]
Ein Kegel $\sigma\subseteq\R^n$ mit $\dim(\s)=d$ heißt \emph{unimodular}, wenn ein 
$k\in\{0,\ldots, d\}$ existiert, sodass $\s$ gitteräquivalent ist zu dem Kegel
$\R_{\ge 0}^k \times \R^{d-k}$.

Ein $d$-dimensionaler Simplex $\sigma\subseteq\R^n$ heißt \emph{unimodular}, wenn  $\s$ gitteräquivalent ist zu dem Einheitssimplex
$\conv(\vek{0},\vek{e_1},\ldots,\vek{e_d})\subseteq \R^d$
\end{Definition}

\begin{Bemerkung}
Ein spitzer Kegel $\s\subseteq \R^n$ ist unimodular genau dann, wenn der Kegel simplizial ist und $\vol(\sigma)=1$ gilt.

Ein Simplex $\s\subseteq \R^n$ ist unimodular genau dann, wenn $\vol(\s)=1$ gilt.

\end{Bemerkung}

\begin{figure}[htb]
\begin{center}
 
\scalebox{0.7}{
	\input{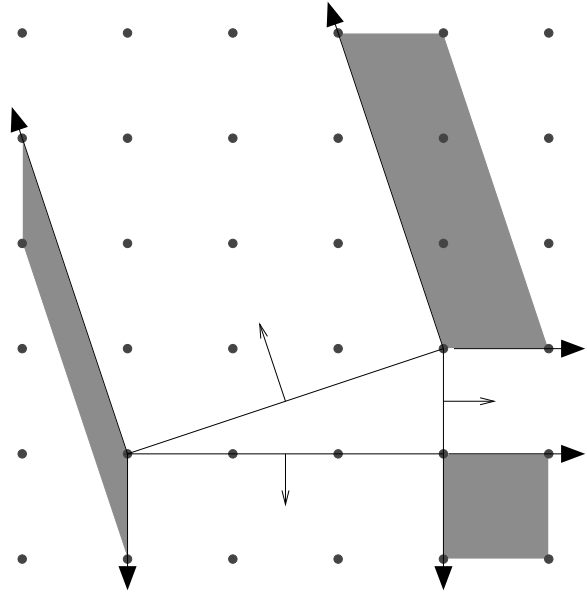}
} 

\caption{Ein Gitterpolytop im $\R^2$, das nicht glatt ist. Der Normalenkegel an der Ecke oben rechts ist nicht unimodular, denn das von seinen primitiven Erzeugern aufgespannte halboffene Parallelepiped enthält außer der Ecke noch zwei weitere Gitterpunkte. 
}
\label{figure:NichtGlatt}
\end{center}
\end{figure}

\begin{Definition}[Normalenkegel]
Für ein Polytop $P$ ist der (innere) \emph{Normalenkegel} an eine Seite  $S$  definiert als 
\[
\nc{S}{P}:=\left\{ \varphi \in (\R^n)^*\,\middle|\; \varphi(s)=\min_{x\in P} \varphi(x)\; \forall s\in S\right\}\,.
\]
\end{Definition}
In Abbildung  \ref{figure:NichtGlatt} ist ein Polytop zu sehen, bei dem an allen Ecken $\v$ der verschobene äußere Normalenkegel $-\nc{\v}{P}+\v$ eingezeichnet ist.

 \begin{Bemerkung}
 \begin{itemize}
	\item Der Normalenkegel ist stets volldimensional und enthält den Raum $U^\perp$, wobei mit $U$ der zu $\aff(P)$ gehörige Unterraum bezeichnet wird.
 	\item Sei $F$ eine Facette von $P$ und $n$ ein (innerer) Normalenvektor an $F$. Dann gilt $\nc{F}{P}=\cone(n)+U^{\perp}$.
 	\item Für Seiten $S\seite P$ lässt sich $\nc{S}{P}$  
 	schreiben als $\cone\{n_1,\ldots, n_k\} + U^{\perp}$, 
  		wobei $n_1,\ldots, n_k$  Normalenvektoren an die Facetten sind, die $S$ enthalten. %
 \end{itemize}
 \end{Bemerkung}

\begin{Definition}[glatt]
Ein Polytop heißt \emph{glatt}, wenn der Normalenkegel an allen Ecken unimodular ist. 
\end{Definition}

Die Frage, ob ein Gitterpolytop $P$ glatt ist, ist aus algebraischer Sicht deshalb interessant, weil diese Eigenschaft des Polytops äquivalent dazu ist, dass die zugehörige torische Varietät $X_P$ glatt ist (s. \cite[§2.1]{FultonToric}). In Abschnitt \ref{section:3Kreuz4Glatt} werden wir sehen, dass glatte Polytope auch aus kombinatorischer Sicht gute Eigenschaften haben.

\begin{Beispiel}
Betrachte das Polytop \[P=\conv(\vek{v_1},\vek{v_2},\vek{v_3}) \quad\text{mit }
\vek{v_1}=
\begin{bmatrix}0\\0\end{bmatrix},\;
\vek{v_2}=\begin{bmatrix}3\\0\end{bmatrix},\;
\vek{v_3}=\begin{bmatrix}3\\1\end{bmatrix}
 \subseteq \R^2\]
 (s. Abb. \ref{figure:NichtGlatt}). 
Die Normalenkegel an den drei Ecken haben die Determinanten
\begin{align*}
\det(\nc{\vek{v_1}}{P})&=\det
\begin{bmatrix}0 & -1  \\ -1 & 3\end{bmatrix}
=-1\,,
\\
\det(\nc{\vek{v_2}}{P})&=\det
\begin{bmatrix}1 & 0  \\ 0 & -1\end{bmatrix}
=-1 \text{ und}
\\
\det(\nc{\vek{v_3}}{P})&=\det
\begin{bmatrix}1 & -1  \\ 0 & 3\end{bmatrix}
=3\,.
\end{align*}
Der Normalenkegel an $\vek{v_3}$ ist also nicht unimodular und deshalb ist $P$ nicht glatt. 
\end{Beispiel}

\begin{Bemerkung}
\label{Bemerkung:GlattFolgtEinfach}
Jedes glatte Polytop ist auch einfach.
\end{Bemerkung}
\begin{proof}
Sei $P$ ein $d$-dimensionales glattes Polytop. Sei $\v$ eine Ecke, die in $k$ Facetten liegt.
Nach der Bemerkung von oben gilt $\nc{\v}{P}=\cone(n_1,\ldots, n_k) + U^{\perp}$.

$\nc{\v}{P}$ ist unimodular, also gitteräquivalent zur Minkowskisumme eines 
simplizialen Kegels und eines Unterraumes. Damit muss aber auch schon $\nc{\v}{P}$ die Minkowskisumme eines simplizialen Kegels und eines Unterraumes sein. 
Folglich ist $\cone(n_1,\ldots, n_k)$ simplizial und damit $k=d$.
$\,$
\end{proof}

\section{Flusspolytope}

In diesem Abschnitt definieren wir Flusspolytope und beweisen einige Sätze dazu, die wir später benötigen. 
Ein Flusspolytop ist die Menge aller Flüsse auf einem gegebenen Graphen, die bestimmte Randbedingungen erfüllen, wobei wir einen Fluss als Punkt im $\R^{\vec{E}}$ auffassen. Ein interessanter Spezialfall von Flusspolytopen sind Transportpolytope, die man erhält, wenn der dem Flusspolytop zugrunde liegende Graph vollständig bipartit ist.

\begin{Definition}[Flusspolytope]
\index{Flusspolytop}
Seien ein gerichteter Graph  $\vec{G} = (V,\vec{E})$  mit Inzidenzmatrix $I_{\vec{G}}$,
der \emph{Bedarfsvektor} $\vek{d} \in \Z^V\setminus\{\vek{0}\}$\footnote{Den Fall $\vek{d}=\vek{0}$ verbieten wir, weil $F$ dann i.\,A. nicht mehr homogen ist, d.\,h. es ist möglich, dass $f_1,f_2,f_3\in F$ existieren mit $f_1+f_2=f_3$. Im 3. und 4. Kapitel 
benötigen wir aber, dass genau so etwas nicht auftreten kann, d.\,h., dass unsere Polytope in einer affinen Hyperebene liegen, die nicht den Ursprung enthält. 

Alternativ könnte man auch $\vek{d}=\vek{0}$ zulassen und stattdessen fordern, dass $\vec{G}$ kreisfrei ist.
Falls man weder $\vek{d}=\vek{0}$ noch Kreise in $\vec{G}$ ausschließen möchte, 
so kann man statt $F$ einfach das zu $F$ isomorphe Polytop $F\times \{1\}\subseteq \R^{|\!\vec{E}|+1}$ betrachten.
},
sowie obere und untere Schranken $\vek{u},\vek{l}\in \Z^{\vec{E}}$ gegeben. 
Dann definieren wir das \emph{Flusspolytop} $F_{\vec{G},\vek{d},\vek{u},\vek{l}}$ folgendermaßen:

 \begin{align}
 F= F_{\vec{G},\vek{d},\vek{u},\vek{l}}:=&\left\{f : \vec{E} \to \R_{\ge 0} \,\Biggm|\, \sum_{e \in \din(v)}f(e) - \sum_{e \in \dout(v)}f(e) = d_v\right.\nonumber\\
&\hspace{3.2cm}  l_e \le f(e) \le u_e  \Biggr\}\\
   =& \left\{\vek{f} \in \R^{\vec{E}}_{\ge 0} \,\middle|\, I_G \cdot\vek{f} = \vek{d},\, \vek{l} \le \vek{f} \le \vek{u} \right\}
 \end{align}

\end{Definition}

Verzichtet man auf die obere Schranke bzw. setzt $\vek{u}=\vek\infty$, so erhält man ein Flusspolyeder, welches i.\,A. aber nicht beschränkt ist. Ist $\vec{G}$ kreisfrei, so erhält man weiterhin ein Flusspolytop, da dann der Fluss über alle Kanten z.\,B. durch 
$\sum_{v\in V} \abs{d_v}$ nach oben beschränkt ist.

Offensichtlich muss $\sum_{v\in V} d_v=0$ erfüllt sein, damit $F$ nicht leer ist.

\begin{Definition}[Transportpolytope]
\index{Transportpolytop}
\index{Polytop!Transportpolytop}
	Seien $m,n \in \N, \vek{r} \in \N^m, \vek{c} \in \N^n$. Die Menge der $(m \times n)$-Matrizen mit positiven Eintr\"agen, Zeilensummen $r_i$ und Spaltensummen $c_j$ bezeichnen wir als Transportpolytop $\trans{r}{c}$.
\[
\trans{r}{c} = \left\{ A \in \R_{\ge 0}^{m \times n} \,\middle|\, \sum_{i=1}^m a_{ij}=c_j, 
\sum_{j=1}^n a_{ij}=r_i \right\}
\]  
\end{Definition}

\begin{Beispiel}[Birkhoffpolytop]
\index{Birkhoffpolytop}
\index{Polytop!Birkhoffpolytop}
\index{Matrix!doppelt stochastische }
Das bekannteste Beispiel für Transportpolytope sind die \emph{Birkhoffpolytope} 
$B_n:=\trans{r}{c}$, mit $\vek{r}=\vek{c}=(1,\ldots,1)\in\R^n$.

$B_n$ ist gerade die Menge der doppelt-stochastischen ($n\times n$)-Matrizen.
Die Ecken von $B_n$ sind die ($n\times n$)-Permutationsmatrizen.
\end{Beispiel}

\begin{Bemerkung}
\label{Bemerkung:DimensionLieblingspolytope}
\begin{enumerate}[(i)]
  \item Für ein Flusspolytop  $F=F_{\vec{G},\vek{d},\vek{u},\vek{l}}$ gilt:\\
  $\dim(F)\le \abs{\vec{E}}-\abs{V}+ k$, wobei $k$ die Anzahl der Zusammenhangskomponenten von $\vec{G}$, aufgefasst als ungerichteter Graph, bezeichnet.
  \label{Bemerkung:DimensionFlusspolytop}
  Flusspolytope, die diese Schranke mit Gleichheit erfüllen, nennen wir \emph{maximaldimensional}.
 \item Transportpolytope sind Flusspolytope mit $\vec{G} = \vec{K}_{m,n}$,  $\vek{u}=\vek{\infty}$ und $\vek{l}=\vek{0}$ sowie $\vek{d}=(-r_1,\ldots,-r_m,c_1,\ldots, c_n)$.
 \item Als Spezialfall von (\ref{Bemerkung:DimensionFlusspolytop}) folgt
 für Transportpolytope $\dim (\trans{r}{c})\le (m-1)(n-1)$. Gleichheit gilt, falls $\vek{r},\vek{c}>0$.
  \label{Bemerkung:DimensionTransportpolytop}
\end{enumerate}

\end{Bemerkung}

Für einen Beweis von (\ref{Bemerkung:DimensionFlusspolytop}), siehe \cite[Seite 208]{schrijverCO}. Der zweite Satz von (\ref{Bemerkung:DimensionTransportpolytop}) ist Theorem 21.16 des gleichen Buches.

\begin{Definition}
Sei $A \in \R^{m\times n}$. $A$ hei\ss t \emph{vollst\"andig unimodular}, falls f\"ur alle quadratischen Untermatrizen $C$ gilt: $\det(C) \in \{0,1,-1\}$.
\end{Definition}

\begin{Lemma}[Erkennung vollst\"andig unimodularer Matrizen]
\label{Lemma:ErkennungTUM}
$\,$\\
Sei $A \in \R^{m\times n}$. $A$ ist vollst\"andig unimodular, falls $A$ die folgenden Bedingungen erf\"ullt.
	\begin{enumerate}[(i)]
		\item Alle Eintr\"age von $A$ sind aus $\{0,1,-1\}$.
		\item In jeder Spalte gibt es h\"ochstens zwei Eintr\"age ungleich Null.
		\item Die Menge der Zeilen von $A$ l\"asst sich partitionieren in zwei Mengen
		$I_1$ und $I_2$, sodass f\"ur jede Spalte mit zwei Eintr\"agen, die verschieden von Null sind, gilt:
		\begin{itemize}
			\item Haben die Eintr\"age ein verschiedenes Vorzeichen, so liegen die zugeh\"origen Zeilen in der gleichen Menge.
			\item Haben die Eintr\"age das gleiche Vorzeichen, so liegen die zugeh\"origen Zeilen in unterschiedlichen Mengen.
		\end{itemize}
	\end{enumerate}
\end{Lemma}

\begin{proof}
Der Beweis orientiert sich an \cite{moehringADM}.\par 
Induktion \"uber die Gr\"o\ss e $k$ der quadratischen Teilmatrix $C$. Der Fall $k=1$ ist klar.
Sei also $k\ge 2$:\par\smallskip
\begin{Fallunterscheidung}
\Fall{
Es gibt eine Spalte, die h\"ochstens einen Eintrag ungleich Null enth\"alt.
Entwickle $\det(C)$ nach dieser Spalte und wende die Induktionsvoraussetzung an.
}
\smallskip

\Fall{
Alle Spalten haben mindestens zwei Eintr\"age, die verschieden von Null sind. Betrachte die Aufteilung der Zeilen in $I_1$ und $I_2$. F\"ur jede Spalte $j$ gilt:
\[
\sum_{i \in I_1} c_{ij} = \sum_{i \in I_2} c_{ij}
\]
Damit erh\"alt man
\[
\sum_{i \in I_1} c_{i} = \sum_{i \in I_2} c_{i}\,,
\]
wobei $c_i$ die $i$-te Zeile von $C$ bezeichnet.
Die Zeilen von $C$ sind also linear abh\"angig, daraus folgt $\det(C)=0$. 
}
\end{Fallunterscheidung}
\vspace{-5pt}
$\,$
\end{proof}

\begin{Korollar}
	\label{Folgerung:UnimodGraphen}
		 Sei $A$ die Inzidenzmatrix eines gerichteten Graphen $\vec{G}= (V,\vec{E})$. Dann ist $A$ vollst\"andig unimodular.			
\end{Korollar}

\begin{proof}
Wir identifizieren die Menge der Zeilen von $A$ mit der Menge $V$ der Knoten von $\vec{G}$ und setzen  $I_1=V$, $I_2=\emptyset$.
$\,$
\end{proof}

\begin{Lemma}
\label{unimodLemma}
Sei $A$ vollst\"andig unimodular. Dann ist auch 
$\binom A I$
vollst\"andig unimodular.
\end{Lemma}

\begin{proof}
Betrachte eine quadratische Untermatrix $C$ von $\binom A I$. Enth\"alt $C$
nur Zeilen von $A$, so sind wir nach Voraussetzung fertig. Enth\"alt $C$ eine Zeile von $I$, 
so entwickeln wir nach dieser. Die Aussage folgt dann per Induktion.
\end{proof}

\begin{Satz}
\label{Satz:FlusspolytopeSindGP}
Flusspolytope sind Gitterpolytope.
\end{Satz}

\begin{proof}
Sei $F=F_{\vec{G},\vek{d},\vek{u},\vek{l}}$ ein Flusspolytop, $A$ die Inzidenzmatrix von $\vec{G}$ und $\v \in \R^{\vec{E}}$ eine Ecke von $F$, sowie $k:=\abs{\vec{E}}$. Wir wissen: 
\[
\v = \bigcap_{\substack{G \text{ Facette}\\\v \in G}}\hspace{-8pt} G
\]
Die Hyperebenen, in denen die Facetten liegen, sind alle durch Gleichungen der Form $v_e= l_e$ oder $v_e=u_e$ gegeben.
Wir wissen also, dass $\vek{v}$ durch die Gleichung $A\v=\vek{d}$ sowie einige der Gleichungen
der beiden linearen Gleichungssysteme $I_k\v=\vek{l}$ und $I_k\v=\vek{u}$ eindeutig bestimmt ist.
Folglich existiert eine $(k\times k)$-Untermatrix $B$ von $(A^T\: I_k\; I_k)^T$ mit vollem Rang und ein Vektor $\vek{c}\in \R^k$, dessen Eintr\"age die zu den $B$ Zeilen gehörenden Einträge aus $\vek{d}, \vek{l}$ und $\u$ sind, sodass $\v$ die eindeutig bestimmte L\"osung von $B\x=\vek{c}$ ist. Nach Korollar \ref{Folgerung:UnimodGraphen}, Lemma \ref{unimodLemma} und da Untermatrizen vollst\"andig unimodularer Matrizen wieder vollst\"andig unimodular sind, ist also $\abs{\det(B)}=1$. 
Aus der Cramerschen Regel folgt dann:
\[
v_i = \frac{\det(B^i)}{\det(B)} = \pm \det(B^i) \in \Z,
\]
wobei $B^i$ die Matrix bezeichnet, die aus $B$ entsteht, indem die $i$-te Spalte 
von $B$ durch den Vektor $\vek{c}$ ersetzt wird.
\end{proof}
Damit folgt natürlich automatisch, dass auch Transportpolytope Gitterpolytope sind.

\bigskip

Wir werden nun zeigen, dass sich für ein Flusspolytop $F$ und eine natürliche Zahl $k$ jeder Gitterpunkt aus 
$k\cdot F$ als Summe von $k$ Gitterpunkten aus $F$ schreiben lässt.
Diese Aussage ist eines unserer zentralen Hilfsmittel in Kapitel \ref{chapter:Gradschranken}. Zum Beweis benötigen wir den folgenden Satz:

\begin{Satz}[Existenzkriterium für Flüsse]
\label{Satz:ExistenzVonFluessen}
Seien ein gerichteter Graph  $\vec{G} = (V,\vec{E})$,
ein Bedarfsvektor $\vek{d} \in \Z^V\setminus \{\vek{0}\}$  mit $\sum_{v\in V} d_v=0$, sowie untere und obere Schranken $\vek{l},\vek{u} \in \Z^{\vec{E}}$ mit $\vek{l}\le\vek{u}$ gegeben.

Dann gilt $F_{\vec{G},\vek{d},\vek{u},\vek{l}}\not=\emptyset$
, d.\,h. es gibt einen Fluss zu diesen Parametern
 genau dann, wenn

\begin{equation}
\sum_{\substack{\din(e) \in U \\\dout(e)\not\in U}}u_e\,   -
\sum_{\substack{\dout(e) \in U \\\din(e)\not\in U}}l_e  
\ge \sum_{v\in U} d_v\quad\text{ für alle $U\subseteq V$ gilt.} 
\label{equation:FlussBedingung}
\end{equation}

\end{Satz}

\begin{figure}[htbp]
\centering
{

\scalebox{0.9}{
	\input{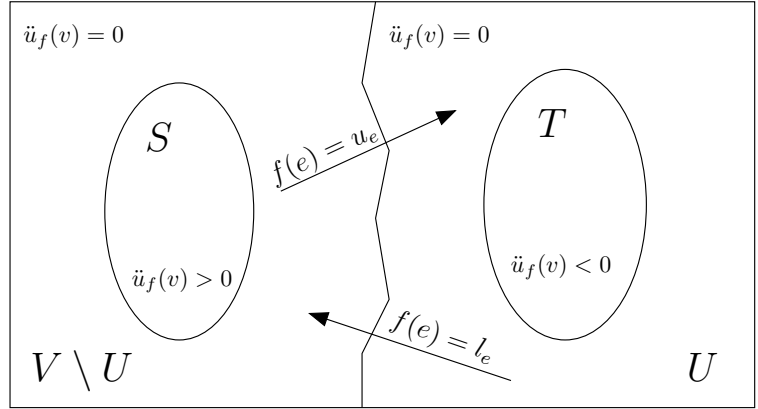}
}
}

\caption{Schematische Darstellung der Situation im Beweis von Satz \ref{Satz:ExistenzVonFluessen}}
\label{figure:Flussbeweis}
\end{figure}

\begin{proof}
Dieser Beweis ist eine an unsere Situation angepasste Version des Beweises von Satz 11.2 aus  
\cite{schrijverCO}.

\Hinrichtung Klar.

\Rueckrichtung
Für eine Funktion $f : \vec{E}\to \Z$ definieren wir die \emph{Überschussfunktion} $\ddot{u}_f : V\to \Z$ gemäß:
\begin{equation}
\ddot{u}_f(v):=\sum_{\din(e)=v}f(e)\,-\sum_{\dout(e)=v}f(e) - d_v
\end{equation}
 Gilt $l_e\le f(e)\le u_e\,\forall e\in \vec{E}$, so definieren wir den  \emph{Restgraphen} $G_f=(V,\vec{E}_f)$ mit $\vec{E}_f:=
\{(u,v)\,|\, e=(u,v)\in\vec{E}, f(e)<u_e\} \cup \{(v,u)\,|\, e=(u,v)\in\vec{E}, f(e)>l_e\}$.

Angenommen $F_{\vec{G},\vek{d},\vek{u},\vek{l}}=\emptyset$. 
Wegen $\vek{l}\le\vek{u}$ existiert eine Funktion $f : \vec{E}\to \Z$ mit $l_e\le f(e)\le u_e$ für alle $e\in\vec{E}$.
Wir wählen $f$ so, dass $\ddot{u}(f):=\sum_{v\in V}\abs{\ddot{u}_f(v)}$ minimal ist. Sei $S:=\{v\in V\,|\, \ddot{u}_f(v)>0\}$ und $T:=\{v\in V\,|\, \ddot{u}_f(v)<0\}$. 
Aus unserer Annahme folgt $\ddot{u}(f)>0$, denn wenn $\ddot u(f)=0$ gelten würde, so wäre $f$ in $F_{\vec{G},\vek{d},\vek{u},\vek{l}}$ enthalten. 
Da außerdem $\sum_{v\in V}\ddot{u}_f(v)=0$ gilt, sind $S$ und $T$ beide nicht leer.

Sei $U$ die Menge der Knoten von denen aus in $\vec{G}_f$ ein Weg zu einem Knoten in $T$ existiert.
Es gilt $U\cap S=\emptyset$. Gäbe es nämlich einen Weg von einem Knoten in $S$ zu einem Knoten in $T$, so könnten wir $f$ entlang dieses Weges erhöhen und dadurch  $\ddot{u}(f)$ reduzieren. 

Nach Wahl von $U$ gilt für Kanten, die $U$ mit $V\setminus U$ verbinden:
\begin{equation}
f(e)=\begin{cases}
l_e &  \text{für $\dout(e)\in U$ und $\din(e)\not\in U$ ($e$ ist aus $U$ ausgehende Kante)}   \\
u_e &  \text{für $\dout(e)\not\in U$ und $\din(e)\in U$ ($e$ ist nach $U$ eingehende Kante)}   
\label{equation:WerteKanten}
\end{cases}
\end{equation}
Für eine schematische Darstellung der Situation siehe Abbildung \ref{figure:Flussbeweis}.
Es folgt:
\begin{align*}
\sum_{\substack{\din(e) \in U \\\dout(e)\not\in U}}u_e   -
\sum_{\substack{\dout(e) \in U \\\din(e)\not\in U}}l_e  
&\gleich{(\ref{equation:WerteKanten})} 
\sum_{\substack{\din(e) \in U \\\dout(e)\not\in U}}f(e)   -
\sum_{\substack{\dout(e) \in U \\\din(e)\not\in U}}f(e)\\
&\gleich{$(*)$}\sum_{v\in U} \left(\sum_{\din(e)=v}f(e)   - \sum_{\dout(e)=v}f(e)\right)\\  
&\gleich{$\,$}\sum_{v\in U} (\ddot{u}_f(v) + d_v) \\
&\gleich{$\,$} \underbrace{\sum_{v\in T} \ddot{u}_f(v)}_{<0} + \underbrace{\sum_{v\in U\setminus T} \ddot{u}_f(v)}_{=0} + \sum_{v\in U}d_v \\
&\kleiner{$\,$} \sum_{v\in U}d_v \qquad\text{\wid\; zu (\ref{equation:FlussBedingung})}
\end{align*}
Gleicheit bei $(*)$ gilt, da nur Terme für Kanten $e$ mit $\din(e)\in U$ und $\dout(e)\in U$ hinzukommen. Für diese wird $f(e)$ einmal addiert und einmal subtrahiert.
$\,$
\end{proof}

Nun kommt der oben angekündigte Satz über Gitterpunkte im $k$-fachen eines Flusspolytops.
Der Satz verallgemeinert den Satz von Birkhoff und von Neumann, der besagt, dass sich jede quadratische Matrix, bei der alle Zeilen- und Spaltensummen den gleichen Wert $k\in \N$ haben, als Summe von $k$ Permutationsmatrizen schreiben lässt. 

\begin{Satz}[Verallgemeinerter Satz von Birkhoff und von Neumann]
$\,$\\
Sei $F=F_{\vec{G},\vek{d},\vek{u},\vek{l}}$ ein Flusspolytop wie oben. Sei $k\in \N$.
Sei $f\in (k\cdot F)\,\cap\,\Z^{\vec{E}}$. Dann existieren Flüsse $f_1,\ldots, f_k \in F\cap \Z^{\vec{E}}$ mit $f=f_1 + \ldots + f_k$.
\label{Satz:VerallgemeinterterBvN}
\end{Satz}

\begin{proof}
In \cite{DBLP:journals/jal/LewandowskiLL86}  wird ein algorithmischer Beweis für den Fall gegeben, dass $F$ ein Transportpolytop ist. Wir zeigen die Aussage unkonstruktiv für den allgemeinen Fall. 

Wir betrachten zunächst den Fall $\vek{l}=\vek{0}$.

Es genügt zu zeigen, dass ein ganzzahliger Fluss $g\in F$ existiert mit $0=l_e\le g(e)\le \min(u_e,f(e))$ für alle $e\in \vec{E}$. Dann folgt die Aussage per Induktion.

Also genügt es zu zeigen, dass für diese verschärfte obere Schranke Bedingung
(\ref{equation:FlussBedingung}) immer noch erfüllt ist.

Sei $U\subseteq V$. Dann gilt: 
\begin{align*}
 \sum_{v\in U} k\cdot d_v   
&=
\sum_{\substack{\din(e) \in U \\\dout(e)\not\in U}} f(e) -  
\sum_{\substack{\dout(e) \in U \\\din(e)\not\in U}} f(e) \\
&\le
\sum_{\substack{\din(e) \in U \\\dout(e)\not\in U}} f(e)\\
&\gleich{} \sum_{\substack{\din(e) \in U \\\dout(e)\not\in U}}  \min(k\cdot u_e, f(e)) 
\displaybreak[2]
\\    
\Longrightarrow\qquad
\sum_{v\in U}  d_v 
&\le  
 \sum_{\substack{\din(e) \in U \\\dout(e)\not\in U}}  \min(u_e, \frac 1 k f(e)) \\     
& \le 
\sum_{\substack{\din(e) \in U \\\dout(e)\not\in U}}  \min(u_e, f(e)) \\     
\end{align*}
Also existiert ein Fluss $g\in F$ mit der gewünschten Eigenschaft.
Damit ist der Fall $\vek{l}=\vek{0}$ erledigt. 

Den Fall $\vek{l}\not=\vek{0}$ kann man nun durch eine Transformation des Graphen darauf zurückführen.

Dazu konstruieren wir einen neuen Graphen $\vec{G}'=(V',\vec{E}')$, indem wir für jeden Knoten $v\in V$ zwei weitere Knoten $v'$ und $v''$, sowie Kanten $(v',v)$ und $(v,v'')$ hinzufügen (s. Abbildung \ref{figure:GraphtransformationUntereSchranken}).

\begin{figure}[htbp]
\centering
\subfloat[Vorher]
{

\scalebox{0.7}{
	\input{UntereSchrankenGraphA.pstex_t}

} 

}
\hspace*{2mm}
\subfloat[Nachher]
{

\scalebox{0.7}{
		\input{UntereSchrankenGraphB.pstex_t}
}

}
\caption{Eine Graphentransformation}
\label{figure:GraphtransformationUntereSchranken}

\end{figure}

Wir definieren nun  $\vek{d'}\in\Z^V$ und $\vek{u'}\in\Z^{\vec{E}}$:
\begin{align*}
 u_{(v',v)}&:=\sum_{\din(e)=v}l_e \; \text{und}
 &
 u_{(v,v'')}&:=\sum_{\dout(e)=v}l_e \text{ für } v\in V &
 u'_e &:=u_e-l_e \text{ für } e\in\vec{E} &  
\\
d'_v&:=d_v \text{ für } v\in V
&
d_{v'}&:=-u_{(v',v)}
\text{ für } v\in V
& 
d_{v''}&:=u_{(v,v'')}
\text{ für } v\in V
\end{align*}

Die Kapazität der Kanten wird also um den Mindestfluss reduziert und stattdessen kommt der Mindesteinfluss in $v$ von einem neuen Knoten $v'$ und der Mindestausfluss von $v$ wird zu einem neuen Knoten $v''$ geleitet.

Sei $\tilde{F}:=F_{\vec{G}',\vek{d'},\vek{u'},\vek{0}}$.
Wir definieren nun eine Funktion $\Phi : \bigcup_{k\in\N} k\cdot F\to \bigcup_{k\in\N} k\cdot F$, die einen Fluss $f\in k\cdot F$ auf einen Fluss
$\tilde f\in k\cdot\tilde{F}$ abbildet, der folgendermaßen definiert ist:

\begin{equation}
\tilde{f}(e):=\begin{cases}
		f(e)-k\cdot l_e & \text{für $e\in \vec{E}$} \\
		k\cdot u_e & \text{sonst}
              \end{cases}
\end{equation}
Wie man leicht sieht, ist dies eine bijektive Abbildung. Die Umkehrabbildung $\Phi^{-1}$ bildet $\tilde f\in k\cdot \tilde F$ auf $f\in k\cdot F$ mit $f(e):=\tilde{f}(e)+ k\cdot l_e$ ab. Außerdem sind $\Phi$ und $\Phi^{-1}$ additiv.

Für $f\in k\cdot F$ existieren also $\tilde{f}_1, \ldots, \tilde{f}_k \in \tilde{F}$ mit 
$\Phi(f)=\tilde{f}=\tilde{f}_1 +  \ldots +  \tilde{f}_k$. Durch Anwenden von $\Phi^{-1}$ erhalten wir  $f_1,\ldots,f_k\in F$ mit $f=f_1 +\ldots + f_k$.
$\,$
\end{proof}

\section{Ein kombinatorisches Kriterium für die Glattheit von Transportpolytopen}
\label{section:TransportKombinatorischGlatt}
Zu überprüfen, ob ein gegebenes Gitterpolytop $P$ glatt ist, ist i.\,A. eine recht komplizierte und rechenintensive Aufgabe. 
Wir zeigen in diesem Abschnitt, dass es für Transportpolytope ein leicht zu überprüfendes  kombinatorisches Kriterium dafür gibt.

Im folgenden Lemma zeigen wir zunächst, dass wir für die Vektoren $\vek{r}$ und $\vek{c}$ gewisse Eigenschaften \oBdA voraussetzen können, indem wir das Polytop ggf. solange verschieben, bis der gewünschte Fall eintritt.

\begin{Lemma}
\label{Lemma:BoeseGleichungenReduzieren}
Sei $\trans{r}{c}$ ein ($m\times n$)-Transportpolytop und sei
\begin{align}\nonumber
\label{equation:BoeseGleichungenReduzieren}
 \chi : \N^m \times \N^n &\rightarrow \pot\bigl((\pot([m])\times \pot([n])\bigr) \\\nonumber
(\vek{a},\vek{b}) & \mapsto 
\Bigl\{(I,J) \Bigm|\,
\emptyset \subsetneq I \subsetneq [m],\, \emptyset \subsetneq J \subsetneq [n],\\
&  \qquad \abs{I}\cdot\abs{\kompl{J}}=1 \textrm{ oder } \abs{\kompl{I}}\cdot\abs{J}=1,
\\\nonumber
&  \qquad {\textstyle\sum_{i\in I}a_i = \sum_{j \in J} b_j} 
\Bigr\}\;.\\\nonumber
\end{align}

Dann existiert ein $A\in\Z^{m\times n}$ mit $\trans{r}{c}+A=\trans{r'}{c'}$ für ein Transportpolytop $\trans{r'}{c'}$ mit 
$\chi(\vek{r'},\vek{c'}) = \emptyset $.

Insbesondere gilt, dass $\trans{r}{c}$ genau dann glatt ist, wenn $\trans{r'}{c'}$ glatt ist.
\end{Lemma}

\begin{proof}
Sei ein Transportpolytop $\trans{r}{c}$ gegeben. 
Sei $|\chi(\vek{r},\vek{c})|$ minimal für alle Paare $(\vek{r},\vek{c})$ aus
$\{(\vek{r'},\vek{c'}) \,|\, \exists A\in\Z^{m\times n} : 
\trans{r}{c}+A=\trans{r'}{c'} \}$.  

Angenommen $\chi(\vek{r},\vek{c})\not= \emptyset$.
Dann existiert $(I,J) \in \chi(\vek{r},\vek{c})$.
Sei \oBdA $J=\{1\}$ und $I=\{1,\ldots, m-1\}$, d.h es gilt:
\begin{equation}
r_1+\ldots + r_{m-1}=c_1
\label{boesegleichung}
\end{equation}

Definiere $\vek{r'}:= (r_1,\ldots, r_{m-1}, r_m + G)$, 
$\vek{c'}:= (c_1 + G , c_2, \ldots, c_n)$ für eine hinreichend große Konstante $G$, z.\,B. 
$G=m\cdot n\cdot \sum_i r_i$.
Dann gilt $\trans{r}{c}+A=\trans{r'}{c'}$ für

\begin{equation*}
A=
\begin{bmatrix}
0 &  0 & \ldots & 0 \\
\vdots &  \vdots & \ddots & \vdots \\
0 & 0   & \ldots & 0\\
G & 0 &  \ldots & 0 
\end{bmatrix}\:.
\end{equation*}
Für jedes $B \in \trans{r'}{c'}$ muss nämlich gelten $b_{m1} \ge G$. Die einzige Matrix für die Gleichheit gilt ist folgende: 
\begin{equation*}
B=
\begin{bmatrix}
r_1 & 0 & 0 & \ldots & 0 \\
r_2 & 0 & 0 & \ldots & 0\\
\vdots & \vdots & \vdots & \ddots & \vdots \\
r_{m-1} & 0 &0  & \ldots & 0\\
G & c_2 & c_3 & \ldots & c_n 
\end{bmatrix}
\end{equation*}

$\vek{r'}$ und $\vek{c'}$ erfüllen die Gleichung (\ref{boesegleichung}) nicht mehr und, da $G$ sehr groß gewählt wurde, auch keine zusätzliche Gleichung aus (\ref{equation:BoeseGleichungenReduzieren}). Also ist $\chi(\vek{r'},\vek{c'}) \subsetneq \chi(\vek{r},\vek{c})$. Dies ist ein Widerspruch zur 
Minimalität von $\abs{\chi(\vek{r},\vek{c})}$. 
$\,$
\end{proof}

\begin{Beispiel}
\label{Beispiel:IsoTPglattKomb}
Betrachte   
$\transv{(1,1,6)}{(2,2,2,2)}$ und $\transv{(1,1,10)}{(3,3,3,3)}$. 

Es gilt:
\begin{align*}
\chi((1,1,6)(2,2,2,2)) &= \bigl\{  (\{3\},\{1,2,3\}),(\{3\},\{1,2,4\}),(\{3\},\{1,3,4\}), \\
				&\qquad (\{3\},\{2,3,4\}), (\{1,2\},\{1\}), (\{1,2\},\{2\}), \\
				&\qquad  (\{1,2\},\{3\}), (\{1,2\},\{4\})  \bigr\} \\
\chi({(1,1,10)}{(3,3,3,3)})& = \emptyset
\displaybreak[1]
\\
 \transv{(1,1,10)}{(3,3,3,3)} 
&=
\begin{bmatrix}
 0 & 0 & 0 & 0 \\
 0 & 0 & 0 & 0 \\
 1 & 1 & 1 & 1 \\
\end{bmatrix} 
+
 \transv{(1,1,6)}{(2,2,2,2)}  
\end{align*}

In Beispiel \ref{Beispiel1110_3333} auf Seite \pageref{Beispiel1110_3333} werden alle ganzzahligen Punkte von $\transv{(1,1,6)}{(2,2,2,2)}$ aufgelistet.
\end{Beispiel}

\begin{Satz}
\label{Satz:KombinatorischGlatt}
Sei $\trans{r}{c}$ ein Transportpolytop. Dann sind die folgenden Bedingungen äquivalent:
\begin{enumerate}[(i)]
\item $\trans{r}{c}$ ist glatt \label{Satz:FallGlatt}
\item $\trans{r}{c}$ ist einfach \label{Satz:FallEinfach}
\item $\sum_{i \in I} r_i \not= \sum_{j\in J} c_j$ für alle Paare $(I,J)$ mit
$\emptyset \subsetneq I \subsetneq [m]$, $\emptyset \subsetneq J \subsetneq [n]$,
$|I|\cdot|\kompl{J}|>1$ und $|\kompl{I}|\cdot|J|>1$\label{Satz:FallKombinatorik}
\end{enumerate}

\end{Satz}

\begin{proof}
Dieser Satz ist eine korrigierte Fassung von Lemma 1.3. aus \cite{christian-andreas-GBTP}.
Teile des Beweises wurden von dort übernommen.

Im Beweis wird die Tatsache verwendet, dass alle Facetten von $\trans{r}{c}$ die Form 
$[a_{ij}\ge 0]$ haben.
\smallskip

\HinrichtungB{(\ref{Satz:FallGlatt})}{(\ref{Satz:FallEinfach})}
Gilt für alle Gitterpolytope, siehe Bemerkung \ref{Bemerkung:GlattFolgtEinfach}.
\smallskip

\HinrichtungB{(\ref{Satz:FallEinfach})}{(\ref{Satz:FallGlatt})} 
Gilt sogar für allgemeinere Polytope, wie das folgende Lemma zeigt:

\begin{Lemma}
Sei $A$ eine vollständig unimodulare $(r\times e)$-Matrix.
Sei $P=\{\x \in \R^e \,|\,$ $ \x\ge 0,\, A\x=\vek{b}\}$ ein einfaches Polytop der Dimension $d$. Sei $U$ der zu $\aff(P)$ gehörige Unterraum. Es gelte $U^\perp=[\text{von $A$ Zeilen aufgespannter Unterraum}]$.

Dann ist $P$ glatt.
\end{Lemma}
Nach Korollar \ref{Folgerung:UnimodGraphen} sind Inzidenzmatrizen von Graphen stets vollständig unimodular. Für Transportpolytope ist das erste Kriterium also stets erfüllt.
$U^{\perp}\supseteq[$von 
$A$ Zeilen aufgespannter Unterraum$]$ gilt immer. Für Transportpolytope mit $\vek{r},\vek{c}>0$ gilt aus Dimensionsgründen Gleichheit, da dann  $\dim(P)=\dim(U)=\dim(\ker(A))$ und damit $\rang(A)=n-d=\dim(U^\perp)$.

\begin{proof}

Sei $\v=(v_1,\ldots v_r)$ eine Ecke von $P$. Da $P$ einfach ist können wir \oBdA annehmen, dass $v_1=\ldots=v_d=0$ und $v_{d+1},\ldots,v_r > 0$ gilt.

Seien $a_1,\ldots,a_l$ die Zeilen von $A$. Fasst man die $a_i$ als Elemente von $(\R^n)^*$ auf, so sind sie auf ganz $P$ konstant. Es 
gilt $\nc{\v}{P}=\cone\{e_1^*,\ldots, e_d^*,\pm a_1,\ldots, \pm a_r\}$. 
Die Ecke $\v$ ist durch die Gleichung
\begin{equation}
\label{equation:MatrixMitId}
\setlength{\arraycolsep}{0pt}
\begin{bmatrix}
\multicolumn{2}{c}{
\text{\LARGE $A$}}\\
 \text{\large $I_d$ } &  \text{\large $0$ }
\end{bmatrix}
\v=
\begin{bmatrix}
\vek{b}\\
\vek{0}
\end{bmatrix}
\end{equation}
eindeutig bestimmt. Insbesondere kann man eine quadratische $(e\times e)$-Untermatrix $B$ mit vollem Rang auswählen, so dass $\v$ durch die entsprechenden Gleichungen immer noch eindeutig bestimmt ist. 

Nach Lemma \ref{unimodLemma} ist die Matrix in (\ref{equation:MatrixMitId}) vollständig unimodular. Daraus folgt $\abs{\det(B)}=1$. Die Zeilen von $B$  sind \oBdA:
$\{e_1^*,\ldots, e_d^*, a_{1},\ldots,   a_{n-d}\}$. 
Folglich erzeugen $a_{1},\ldots,   a_{n-d}$ den Raum $U^\perp$.

Der Normalenkegel wird also erzeugt von $\{e_1^*,\ldots, e_d^*, \pm a_{1},\ldots,  \pm a_{n-d}\}$ und ist unimodular. Wegen $\abs{\det(B)}=1$ ist $B^{-1}\in \SL_e(\Z)$ und liefert uns damit die Gitteräquivalenz zu $\R^d\times \R^{n-d}$. 
$\,$
\end{proof}

\HinrichtungB{(\ref{Satz:FallKombinatorik})}{(\ref{Satz:FallEinfach})} 

Angenommen $\trans{r}{c}$ ist nicht einfach, d.\,h.  es gibt eine Ecke $A$, die zu mindestens 
$(m-1)(n-1)+1$ Facetten gehört. Also hat $A$ mindestens soviele Nulleinträge,
also höchstens $mn-((m-1)(n-1)+1)=m+n-2$ Einträge, die verschieden von Null sind.

Wir betrachten nun den bipartiten Graphen $G$ mit $n+m$ Knoten und einer Kante zwischen zwei Knoten, wenn der zugehörige Eintrag von $A$ ungleich null ist.
$G$ ist unzusammenhängend, da $G$ weniger als $n+m-1$ Kanten hat.
Wähle für $I$ und $J$ die Farbklassen einer Zusammenhangskomponenten von $G$. Es gilt 
$\sum_{i\in I}r_i = \sum_{j \in J} c_j$.
Wegen Lemma \ref{Lemma:BoeseGleichungenReduzieren} können wir 
$|I|\cdot|\kompl{J}|>1$ und  $|I|\cdot|\kompl{J}|> 1$ annehmen. 
\wid

\medskip
\HinrichtungB{(\ref{Satz:FallEinfach})}{(\ref{Satz:FallKombinatorik})}
Angenommen, es gibt Mengen $\emptyset \subsetneq I \subsetneq [m]$ und  $\emptyset \subsetneq J \subsetneq [n]$ mit $\sum_{i\in I}r_i=\sum_{j\in J}c_j$,
 $|I|\cdot|\kompl{J}|>1$ und $|\kompl{I}|\cdot|J|>1$. Wir definieren $m':=|I|$, $m'':=|\kompl{I}|$, $n':=|J|$ und $n'':=|\kompl{J}|$.

Dann existieren Transportpolytope $\trans{r'}{c'}$ und  $\trans{r''}{c''}$ mit $\vek{r'}=(r_i)_{i\in I}$, $\vek{c'}=(c_j)_{j\in J}$,
$\vek{r''}=(r_i)_{i\in \kompl{I}}$ und $\vek{c''}=(c_j)_{j\in \kompl{J}}$. 
Seien $A'$ und $A''$ jeweils Ecken dieser beiden Transportpolytope.
Dann erhalten wir einen Punkt $A\in\trans{r}{c}$ auf die folgende Weise:

\begin{align}\nonumber
& \hspace*{36pt}J\hspace*{30pt} \kompl{J} \\\nopagebreak
\label{equation:MatrixAusZweiEcken}
A&:=
\left[
 \begin{tabular}{ccc}
    & \gc \\
 \hspace*{6pt}$A'$\hspace*{6pt}   & \gc\hspace*{6pt}$0$\hspace*{6pt} \\
    & \gc \\
 \gc&   \\
 \gc\hspace*{6pt}$0$\hspace*{6pt}&  \hspace*{6pt}$A''$\hspace*{6pt} \\
 \gc &   \\
 \end{tabular}
\right]
\begin{tabular}{c}
 \\ 
 $I$ \\ 
 \\ 
 \\
 $\kompl{I}$ \\
 \\
\end{tabular}
\end{align}

Wir wollen nun zeigen, dass $A$ in mehr als $\dim( \trans{r}{c})$ vielen Facetten liegt. Dazu benötigen wir die folgenden drei Lemmata.

\begin{Lemma}
\label{Lemma:TransportFacetten}
$[a_{ij}\ge 0]$ definiert eine Facette von $\trans{r}{c}$ genau dann, wenn
 eine Matrix  $B\in\trans{r}{c}$ existiert mit $b_{ij}=0$ und allen anderen Einträgen positiv
(aber nicht notwendig ganzzahlig).
\end{Lemma}

\begin{proof}
\Hinrichtung
Man wähle einen Punkt $B$ im relativ Inneren der Facette $[b_{ij}\ge 0]$. Dieser liegt also in genau einer Facette des Polytops. Da $[b_{ij}\ge 0]$ für alle $(i,j)$ eine Seite des Polytops definiert, hat $B$ die gewünschte Form.

\Rueckrichtung
\OBdA sei $i=j=1$. Sei $B$ eine Matrix, bei der alle Einträge bis auf $b_{11}$ positiv sind.
$[b_{11} \ge 0]$ definiert eine Seite $S$ des Polytops. 

Es genügt zu zeigen, dass $\dim(S)=\dim(\trans{r}{c})-1$ gilt.
Dazu zeigen wir, dass ein $\eps > 0$ existiert und eine Kugel $B_{\eps}(A) \subseteq S$ mit Radius $\eps$ um $B$ mit $\dim(B_{\eps}(A))= (m-1)(n-1)-1$.
Wählt man $\eps$ so, dass $\abs{\eps} \le \min_{(i,j)\not=(1,1)}b_{ij}$, so kann man die folgende Matrix (die Punkte stehen für Nullen) zu $B$ addieren und erhält einen Punkt aus $S$:

\[
\left[
\mbox{
 \begin{tabular}{ccccc}
 \multicolumn{1}{c}{0} & 
 \cellcolor{TabellenUnterlegfarbe}
  $\ldots$
  &  \gc $\ldots$ &  \gc $\ldots$  &  \\
   \gc  &  \gc $\vdots$ &  \gc $\vdots$ &  \gc $\vdots$ &  $\vdots$ \\ 
   \gc  $\vdots$ &  \gc $\ldots$ &  \gc $\pm\eps$ &  \gc $\ldots$ &  $\mp\eps$ \\ 
  \gc  &  \gc $\ldots$ &  \gc $ \vdots$  &  \gc $\ldots$ &  $\vdots$ \\
  & $\ldots$ & $\mp\eps$ &  & $\pm\eps$
 \end{tabular}
}
\right]
\]

Betrachtet man die Projektion auf den grau unterlegten Bereich, so sieht man, dass ein Würfel der Dimension $(m-1)(n-1)-1$ mit Kantenlänge $2\eps$ in $S$ enthalten ist. 
Also hat $S$ die gewünschte Dimension und ist damit eine Facette.
\end{proof}

\begin{Lemma}
\label{Lemma:UngleichungenA}
Die Ungleichungen $[a_{ij}\ge 0]$ für $(i,j) \in I \times \kompl{J} \cup \kompl{I} \times J$  definieren
Facetten von $\trans{r}{c}$. 
\end{Lemma}

\begin{proof}
Sei $(i,j) \in I \times J^c$. Wähle $A'$ und $A''$ so, dass alle Einträge positiv sind (z.\,B. $a'_{kl}:=\frac {r_k'\cdot c_l'}{s}$ mit $s:=\sum_k r_k'=\sum_l c_l'$) und konstruiere daraus die Matrix $A$, wie in (\ref{equation:MatrixAusZweiEcken}).
Addiere zu $A$ nun die folgende Matrix:

\begin{eqnarray*}
 J\hspace*{150pt} \kompl{J}\hspace*{110pt} \\\nopagebreak
\left[
 \begin{tabular}{ccc|c|ccc}
     $\ldots$ & $-\eps(n''-1)m''$\hspace*{5pt} & $\ldots$ & \gc $0$ & \gc$\ldots$ & \gc $+\eps m''n'$ & \gc$\ldots$ \\\hline
       &  $\ldots$   &  & \gc \vdots &  \gc & \ldots\gc & \gc \\
     $\ldots$ & $-\eps n''m''$ & 
     	$\ldots$ &  \gc $+\eps m''n'$ & \gc $\ldots$ & \gc $+\eps m''n'$ & \gc $\ldots$\\
     & $\ldots$ & &  \gc\vdots & \gc  & \gc $\ldots$ & \gc  \\\hline
 \gc &   \gc $\ldots$ & \gc    & \vdots &   & \ldots &   \\
  \gc $\ldots$ &\gc $+\eps(m'n''-1)$ & \gc $\ldots$ & $-\eps(m'-1)n'$  & $\ldots$  &$-\eps m'n'$ & \ldots \\ 
  \gc &    \gc $\ldots$ & \gc      & \vdots  & & \ldots &    \\
 \end{tabular}
\right]
\begin{tabular}{c}
 \\ 
 \\
 $I$ \\ 
 \\ \\  %
 \\ 
  $\kompl{I}$ \\
 \\
\end{tabular}
\end{eqnarray*}
\smallskip

Die Null sei an der Stelle $(i,j)$.
Da die Zeilen- und Spaltensummen dieser Matrix Null sind, erhalten wir für $\eps$ klein genug wieder eine Matrix $B\in\trans{r}{c}$. 

Nach Voraussetztung gilt $m'n''-1 \ge 1$. Demnach sind in $B$ alle Einträge außer $b_{ij}$ positiv. Wegen Lemma \ref{Lemma:TransportFacetten} folgt damit die Aussage . 
$\,$
\end{proof}

\begin{Lemma}
\label{Lemma:UngleichungenB}
Wenn die Ungleichung $[a_{ij} \ge 0]$ eine Facette von $\trans{r'}{c'}$ definiert, dann definiert 
sie auch eine Facette von $\trans{r}{c}$.
\end{Lemma}
\begin{proof}
Gilt $n'=1$ oder $m'=1$, so ist das Polytop $\trans{r'}{c'}$ nulldimensional und hat damit keine Facetten.
Wir können also $n',m'\ge 2$ voraussetzen.
 
Sei $[a_{ij} \ge 0]$ eine Facette von $\trans{r'}{c'}$.
Nach Lemma \ref{Lemma:TransportFacetten} gibt es eine Matrix $A'$, deren einziger Nulleintrag an der Stelle $(i,j)$ ist. W\"ahle $A''$ so, dass alle Eintr\"age positiv sind und konstruiere $A$ wie in
(\ref{equation:MatrixAusZweiEcken}). Wir k\"onnen nun wieder die Werte der Nulleintr\"age in $A$ erh\"ohen und die Werte der \"ubrigen Eintr\"age reduzieren, sodass $a_{ij}$
der einzige Nulleintrag ist.

Dazu addieren wir zu $A$ die folgende Matrix, deren Zeilen- und Spaltensummen Null sind:

\begin{eqnarray*}
 J\hspace*{145pt} \kompl{J}\hspace*{85pt} \\\nopagebreak
    \left[
    \begin{tabular}{c|ccc|ccc}
       0 & $\ldots$ & $-\eps m''n''$ & $\ldots$ & \gc $\ldots$ & \gc $+\eps m''(n'-1)$ & \gc $\ldots$\\\hline      
      \vdots & & $\ldots$ & & \gc & \gc $\ldots$& \gc \\      	
          $-\eps m''n''$ & \ldots &  $-\eps m''n''$ & \ldots & \gc $\ldots$ & \gc $+\eps m''n'$& \gc $\ldots$ \\      	    
      \vdots && \ldots && \gc & \gc$\ldots$ & \gc \\\hline
       \gc\vdots & \gc& \gc $\ldots$& \gc & & $\ldots$ &\\      	          
      \gc $+\eps(m'-1)n''$ & \gc $\ldots$ & \gc $+\eps m'n''$ & \gc $\ldots$& $\ldots$ &
	  $ -\eps ( m'n'-1)$ & $\ldots$\\
      \gc \vdots  & \gc  & \gc $\ldots$& \gc & & $\ldots$ &\\      	        
    \end{tabular}
\right]
\begin{tabular}{c}
 \\ 
 \\
 $I$ \\ 
 \\ \\
 \\
 $\kompl{I}$ \\
 \\
\end{tabular}
\end{eqnarray*}
Die Null sei an der Stelle $(i,j)$. Wir erhalten also für kleines $\eps$ eine Matrix $B\in\trans{r}{c}$ mit $b_{ij}=0$ und positiven Einträgen an allen anderen Stellen.

$\,$
\end{proof}

\noindent Nun können wir den Beweis beenden. 
\OBdA sei $\trans{r}{c}$ volldimensional. Dann liegt 
der Punkt $A$ insgesamt in
\begin{eqnarray*}
&&\underbrace{(m'-1)(n'-1)+(m''-1)(n''-1)}_{\text{Lemma \ref{Lemma:UngleichungenB}}} + \underbrace{(m'n''+m''n')}_{\text{Lemma \ref{Lemma:UngleichungenA}}} \\
&=& m'n' - m' - n' + 1 + m''n'' - m'' - n'' + 1 + m'n'' + m''n' \\
&=& (m'+m''-1)(n'+n''-1)+1 \\
&=& (m-1)(n-1)+1     
\end{eqnarray*}
vielen Facetten. Damit ist $\trans{r}{c}$ nicht einfach. \wid 
$\,$
\end{proof}

\chapter{Unterteilungen und Triangulierungen \DatumInKlammern}

In diesem Kapitel besch\"aftigen wir uns mit Unterteilungen und Triangulierungen. Wir interessieren uns dabei besonders f\"ur Pullingtriangulierungen und regul\"are Triangulierungen.

\section{Punktkonfigurationen, Unterteilungen und Triangulierungen}

In diesem Abschnitt definieren wir Unterteilungen und Triangulierungen einer Punktkonfiguration auf kombinatorische Art. 
Die Definitionen orientieren sich an \cite{TriangulationsBook}.

Zur Motivation zunächst eine geometrische Definition:
\begin{Definition}[Unterteilungen und Triangulierungen (geometrisch)]
$\,$\\
Sei $P\subseteq \R^n$ ein Polytop. Eine Menge $\Delta$ von Polytopen im $\R^n$ heißt Unterteilung von $P$, wenn gilt:

\begin{itemize}
\item $Q\in\Delta, Q'\seite Q\Rightarrow Q'\in\Delta$
\item $Q,Q'\in \Delta\Rightarrow Q\cap Q'\seite Q$
\item $\bigcup_{Q\in \Delta}Q = P$
\end{itemize}

Eine Unterteilung $\Delta$ heißt Triangulierung, wenn für jedes $Q\in\Delta$ gilt:
$Q$ ist Simplex.

\end{Definition}

Wir werden eine kombinatorische Definition von Unterteilungen verwenden,
 die bei den von uns betrachteten Polytopen zur geometrischen im wesentlichen äquivalent ist, aber technische Vorteile hat. Dazu betrachten wir anstatt des Polytops selbst die Menge seiner Gitterpunkte bzw. die Indexmenge seiner Gitterpunkte.

\begin{Definition}[Punktkonfigurationen]
Eine \emph{Punktkonfiguration} im $\Z^n$ mit Indexmenge $I$ ($\left|I\right|<\infty$) ist 
eine Familie $\A=\{\vek{a_i} \in \Z^n\,|\,i\in I\}$.\footnote{Wer von der Indexmenge verwirrt ist, kann sich in den meisten Fällen statt einer Menge $J\subseteq I$ einfach die Menge $\{\vek{a_j}\,|\,j\in J\}$ oder sogar die Menge $\conv\{\vek{a_j}\,|\,j\in J\}$ vorstellen.
}

Eine Punktkonfiguration $\A\not=\{\vek{0}\}$ heißt \emph{homogen}, wenn ein $\varphi\in(\Z^n)^*$ und ein $c\in \Z\setminus\{ 0\}$ existieren, mit $\varphi(\a)=c$ für alle $\a\in\A$. Geometrisch bedeutet dies, dass alle Punkte von $\A$ auf einer affinen Hyperebene liegen, die nicht den Ursprung enthält.
\index{Punktkonfiguration}
\index{homogen!homogene Punktkonfiguration}
\end{Definition}

In Kapitel \ref{chapter:Gradschranken} werden unsere Punktkonfigurationen die Gitterpunkte eines Flusspolytopes $F=F_{\vec{G},\vek{d},\vek{u},\vek{l}}$ sein. Diese sind homogen.
Sei $v$ ein Knoten des zugrundeliegenden Graphen und $d_v$ eine Komponente von $\vek{d}$, die verschieden von Null ist.
Die Funktion $\varphi_v \in(\Z^{\vec E})^*$ mit  $\varphi_v(f):=\sum_{\din(e)=v}f(e)\,-\sum_{\dout(e)=v}f(e)$ hat auf ganz $F$ den Wert $d_v\not=0$.

Für unsere Definition von Unterteilungen benötigen wir kombinatorische Analoga von konvexgeometrischen Begriffen. Alle Begriffe sind so definiert, wie man sie erwartet, d.\,h. wenn man von einer Indexmenge von Gitterpunkten zu der konvexen Hülle der zugehörigen Gitterpunkte übergeht, erhält man die geometrische Definition. 
\begin{Definition}
\label{Definition:KonvexgeometrieKominatorisch}
Sei $\A=\{\vek{a_i}\,|\,i\in I\}$ eine Punktkonfiguration im $\Z^n$. Sei $J\subseteq I$. Dann
definieren wir:
\begin{itemize}
\item
Die \emph{konvexe Hülle} von $J$ in $\A$: 
\[
\conv\nolimits_\A(J):= \conv\{\vek{a_j} \,|\, j\in J \}
\]
\item
Den \emph{Kegel} von $J$ in $\A$: 
\[
\cone\nolimits_\A(J):= \cone\{\vek{a_j} \,|\, j\in J \}
\]

\item
Das \emph{relativ Innere} von $J$ in $\A$:
\[
\relint\nolimits_\A(J):= \relint(\conv_\A(J))
\]
\item
$J$ ist \emph{affin abhängig} bzw. \emph{affin unabhängig},
 wenn die entsprechende Eigenschaft für die Familie
 $\{\vek{a_j} \,| \, j\in J\}$ von Punkten im $\Z^n$ gilt. 
\index{affin abhängig!kombinatorisch}
\index{affin unabhängig!kombinatorisch}
\item Die \emph{Dimension} von $J$ ist definiert als die Dimension  von $\conv\nolimits_\A(J)$.

\item 
$F\subseteq J$ heißt \emph{Seite} von $J$, wenn eine Seite $K$ von $\conv_\A(J)$ existiert, so dass für alle $j\in J$ gilt: $j\in F \Leftrightarrow \vek{a_j}\in K$.
Wir schreiben dann $F\seite J$. 
Wir nennen $\varphi\in(\Z^n)^*$ einen \emph{Normalenvekor} an $J$, wenn es einen Normalenvektor $\psi\in(\R^n)^*$ an $K$ gibt mit $\varphi(\v)=\psi(\v)$ für alle $\v\in\Z^n$.

\item
Eine Seite $F\seite J$ heißt \emph{echte Seite} von $J$, falls
 $\emptyset\not=F\not=J$ gilt. Eine Seite von Kodimension eins heißt $\emph{Facette}$. 
Nulldimensionale Seiten heißen \emph{Ecken} und die Menge der Ecken von $J$ bezeichnen wir mit $\eck(J)$. %

\item $J$ heißt $d$-dimensionaler \emph{Simplex}, wenn $\abs{J}=d+1$ und $J$ affin unabhängig ist. Insbesondere folgt dann natürlich $\dim(J)=d$. 

\item
Für einen Simplex $J$ definieren wir das \emph{Volumen} $\vol(J)$ als das normalisierte Volumen von $\conv_\A(J)$. 
Wir nennen $J$ \emph{unimodular}, wenn $\vol(J)=1$ gilt.

\item Zu einer Menge $J\subseteq I$ definieren wir eine Art kombinatorischen Abschluss. Dieser enthält die Indizes aller Gitterpunkte, die in der konvexen Hülle der zu $J$ gehörenden Gitterpunkte liegen, also $\overline{J}:=\{i \in I \,|\, \vek{a_i}\in \conv_\A(J)\cap\Z^n\}$.

\end{itemize}

\end{Definition}

\begin{Definition}[Unterteilungen]
Sei $\A$ eine Punktkonfiguration im $\Z^n$ mit Indexmenge $I$.
Eine Menge $\Delta\subseteq \pot(I)$ heißt \emph{Unterteilung} von $\A$,
wenn folgendes gilt:
\begin{itemize}
\item $J\in \Delta, K\seite J \Rightarrow K\in\Delta$
\item $J,K\in \Delta \Rightarrow J\cap K \in \Delta $
\item $\bigcup_{J\in \Delta}\overline{J}=I$
\end{itemize}
Die Elemente von $\Delta$ nennen wir \emph{Zellen}. Zellen, die die gleiche Dimension wie $\A$ haben nennen wir maximal, nulldimensionale Zellen heißen Ecken. Wenn $J\subseteq I$ nicht in $\Delta$ enthalten ist, so nennen wir $J$ eine Nichtseite von $\Delta$.
\end{Definition}

\begin{Beispiel}
Sei $\A$ ein Punktkonfiguration mit Indexmenge $I$.
Das einfachste Beispiel einer Unterteilung von $\A$ ist die Menge
$\Delta_0:=\{ F\subseteq I\,|\, F\seite I\}$.
\end{Beispiel}

\begin{Definition}[Verfeinerungen]
Seien $\Delta$ und $\Delta'$ Unterteilungen einer Punktkonfiguration $\A$. Dann heißt $\Delta'$ \emph{Verfeinerung} von $\Delta$, wenn für jedes $J'\in\Delta'$ ein
$J\in\Delta$ existiert mit $J'\subseteq J$.
\end{Definition}

\begin{Definition}[Triangulierungen]
Sei $\Delta$ eine Unterteilung von $\A$. $\Delta$ heißt \emph{Triangulierung} von $\A$, 
wenn für jedes $\s\in\Delta$ gilt: $\s$ ist Simplex.
\end{Definition}

\begin{Definition}[Unimodulare Triangulierungen]
Wir nennen eine Triangulierung $\Delta$  \emph{unimodular}, wenn alle Simplexe in $\s$ unimodular sind, d.\,h. $\vol(\s)=1$ gilt.
\end{Definition}
Für eine unimodulare Triangulierung genügt es zu fordern, dass alle maximalen Simplexe unimodular sind, weil Seiten unimodularer Simplexe wieder unimodular sind.
Dies folgt beispielsweise aus dem Volumenlemma im nächsten Abschnitt.

\begin{Bemerkung}
Sei $\A$ eine Punktkonfiguration mit Indexmenge $I$ und $\Delta$ eine Triangulierung von $\A$.

Dann gilt:
\begin{equation*}
\conv_\A(I)=\BigDisjUnion_{\sigma\in\Delta} \relint(\conv_\A(\sigma))
\end{equation*}
\label{Bemerkung:PolytopInneresTriangulierung}
\end{Bemerkung}

\section{Pullingunterteilungen}
Eine für uns sehr wichtige Klasse von Unterteilungen sind die Pullingunterteilungen.
Geometrisch erhält man diese, indem man an einer Ecke des Polytops \glqq zieht\grqq. Das Ziehen geht folgendermaßen: Man bettet das Polytop $P\subseteq \R^n\cong \R^n\times\{0\}\subseteq\R^{n+1}$ (bzw. dessen Unterteilung $\Delta$) mittels $(\mathrm{id},0)$ in den $\R^{n+1}$ ein. Dann zieht man eine Ecke $\v$ ein kleines Stück nach unten, d.\,h. man verringert $v_{n+1}$, und betrachtet die konvexe Hülle. Schaut man nun von unten darauf, sieht man, dass das Polytop weiter unterteilt wurde.
Diese \glqq untere konvexe Hülle\grqq\ projiziert man wieder in den $\R^n$ und erhält so eine Pullingunterteilung von $P$ bzw. eine Pullingverfeinerung von $\Delta$.
Die folgende kombinatorische Definition ist äquivalent dazu
(s. z.\,B. \cite{LeeRegular}).
 
Des Weiteren werden wir zeigen, dass man  eine Triangulierung erhält, wenn man mit der trivialen Unterteilung $\Delta_0$ beginnt und diese immer weiter verfeinert, indem man nacheinander an allen Ecken zieht. 

\begin{Definition}[Pullingverfeinerungen]
Sei $\A$ %
eine Punktkonfiguration mit Indexmenge $I$.
Sei $i\in I$   und $\Delta$ eine Unterteilung von $\A$.

Wir definieren die \emph{Pullingverfeinerung} $\pull(\Delta,i)$ von $\Delta$ durch Ziehen an $i$ (bzw. $\vek{a_i}$) folgendermaßen:

\smallskip
\noindent $\sigma\in\pull(\Delta,i)  \text{ genau dann, wenn eine der beiden folgenden Bedingungen erfüllt ist}$:
\begin{compactenum}[(i)]
 \item 
 $i\not\in\sigma  \text{ und }  \sigma\in\Delta$ 
\item $i\in\sigma \text{ und es existieren $\s'\in\Delta$ und $F\seite\s'$ mit $i\in\s'$ und } \sigma=\{i\} \cup F$ 
\end{compactenum}
\end{Definition}

Wie man leicht überprüfen kann ist $\pull(\Delta,i)$ tatsächlich eine Unterteilung von $\A$ und eine Verfeinerung von $\Delta$.
\smallskip

Für ein Tupel $J=(j_1,\ldots,j_k) \subseteq I$ definieren wir  
\begin{align*}
\pull(\Delta,J)&:=
\pull(\ldots(\pull(\pull(\Delta,j_1),j_2)\ldots),j_k)\\
\text{ und }\quad\pull(\A,J)&:=\pull(\Delta_0,J).
\end{align*}

\begin{Satz}[Pullingtriangulierungen]
\label{Satz:PullingIstTriangulierung}
Sei $\A$ eine Punktkonfiguration mit Indexmenge $I$.

Dann ist $\Delta_{\mathrm{pull}}:=\pull(\A,\eck(I))$ eine Triangulierung.

Wir nennen $\Delta_{\mathrm{pull}}$ eine \emph{Pullingtriangulierung} von $\A$. Diese ist i.\,A. nicht eindeutig, sondern abhängig davon, in welcher Reihenfolge an den Ecken gezogen wird.
\end{Satz}

Wir verschieben den Beweis nach hinten (s. Seite \pageref{Beweis:PullingIstTriangulierung}),
da wir dafür das Pullinglemma (Lemma \ref{pullinglemma}) benötigen. 

\begin{figure}[p]
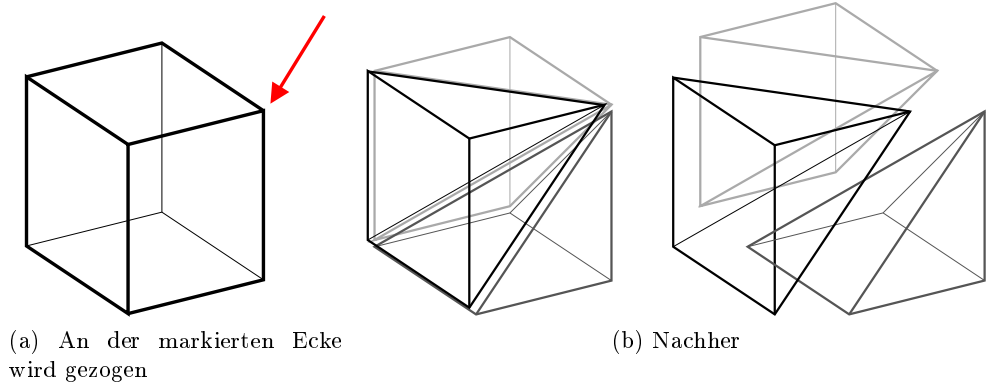

\centering
\subfloat[An der markierten Ecke wird gezogen]
{

\scalebox{0.45}{
	\input{PullingWuerfel1.pstex_t}
} 

}
\subfloat[Nachher]
{

\scalebox{0.45}{
	\input{PullingWuerfel2sw.pstex_t}
}
\hspace*{2mm}
\scalebox{0.45}{
	\input{PullingWuerfel3sw2.pstex_t}
	} 
}
\caption{Pullingunterteilung eines Würfels}
\label{figure:PullingWuerfel}
\end{figure}

Für Beispiele der in diesem Abschnitt vorgestellten Konzepte betrachte man die Abbildungen \ref{figure:Triangulierungen} und \ref{figure:PullingWuerfel}.

\begin{figure}[p]
\centering
\subfloat[Im $\R^2$ haben alle Polytope eine unimodulare Triangulierung.]
{

\scalebox{0.45}{
	\input{TriangulierungPolygon.pstex_t}

} 

}
\hspace*{2mm}
\subfloat[Pulling-Tri\-an\-gu\-lie\-rung eines Polygons. Die großen grauen Punkte markieren die Ecken der Unterteilung, die kleinen die Gitterpunkte.]
{

\scalebox{0.45}{
		\input{PullingTriangulierungPolygon.pstex_t}
}

}
\hspace*{1mm}
\subfloat[Diese Triangulierung eines Würfels mit Kantenlänge $1$ ist nicht unimodular.]
{

\scalebox{0.56}{
	\input{NichtUnimodularWuerfel2.pstex_t}
	} 
}

\caption{Drei Triangulierungen}
\label{figure:Triangulierungen}

\end{figure}

\subsection*{Pacos Lemma}
In diesem Unterabschnitt werden wir sehen, dass
Pullingtriangulierungen für uns eine sehr günstige Eigenschaft haben. Unter gewissen Voraussetzungen sind nämlich alle Zellen einer Pullingtriangulierung unimodular.

\begin{Definition}[Weite]
Sei $\A=\{\vek{a_i}\,|\,i\in I\}$ eine Punktkonfiguration im $\Z^n$ und $J\subseteq I$. %

\begin{itemize}
\item Für $u\in (\Z^n)^*$ definieren wir die \emph{Weite} von $J$ bezüglich $u$ als:
\[
\weite{u}{J}:=\max_{ j\in J}u(\vek{a_j}) - \min_{j\in J}u(\vek{a_j})
\]
\item Die \emph{Weite} von $J$ bezüglich einer Facette $F\seite J$ definieren wir als: 
\[
\weite{F}{J}:=\min \{ \weite{u}{J} \,|\, \text{$u$ Normalenvektor an $F$}\} 
\]
\item $J$ hat \emph{Facettenweite 1} $\Leftrightarrow$ $\weite{F}{J}=1$ für jede Facette $F\seite J$. %
\end{itemize}

\end{Definition}

Wir kennen aus der Analysis folgende Formel für die rekursive Berechnung des Volumens bzw. des Lebesguemaßes eines $d$-dimensionalen Simplex: 
\begin{equation*}
\vol\nolimits_L(\text{\it Simplex})=\frac{1}{d!} \vol\nolimits_L(\text{\it Grundfläche}) \cdot \text{\it Höhe}
\end{equation*}

Eine ähnliche Formel werden wir nun für das normalisierte Volumen eines Gittersimplex beweisen. 

\begin{Lemma}[Volumenlemma]
\label{volumenlemma}
Sei $\A=\{\vek{a_i}\,|\,i\in I\}$ eine Punktkonfiguration im $\Z^n$.
Sei $\sigma=\{i_0,\ldots, i_d\}\subseteq I$ ein Simplex  und $F\seite \sigma$ eine Facette. %

Dann gilt für das normalisierte Volumen von $\s$:
\begin{equation*}
\vol(\sigma)=\weite{F}{\s}\cdot \vol(F)
\end{equation*}

\end{Lemma}

\begin{figure}[htbp]
\begin{center}

\scalebox{0.7}{
\input{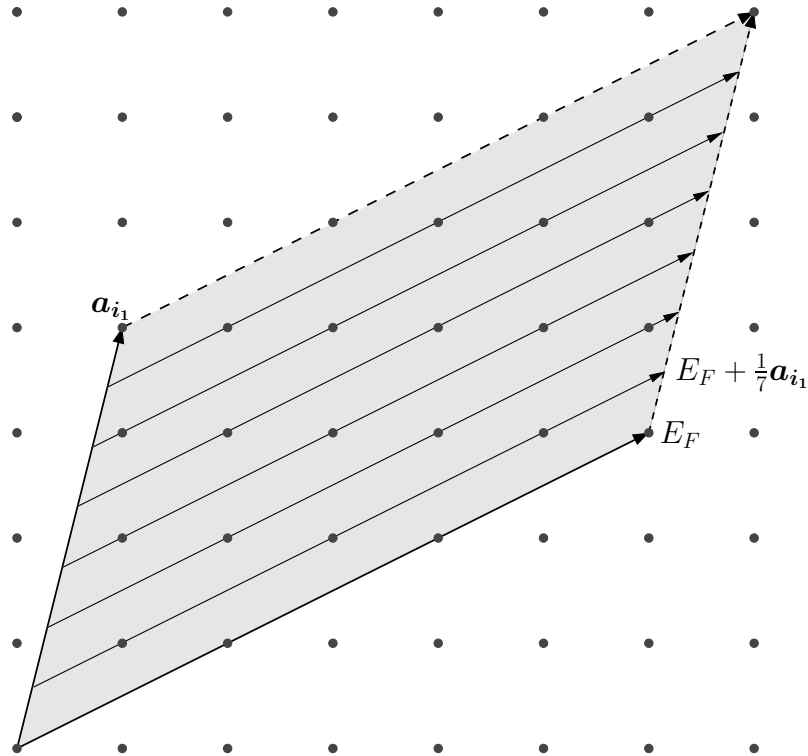}
}

\caption{
Alle Gitterpunkte des halboffenen Parallelepipeds liegen in verschobenen Kopien von $E_F$.
}
\label{figure:VolumenLemma}
\end{center}
\end{figure}

\begin{proof}
\OBdA sei $\vek{a_{i_0}}=\vek{0}$, $F=\{i_0,i_2,\ldots,i_d\}$.
  $E:=\{ \sum_{k=1}^d \l_k \vek{a_{i_k}} \,|\,$ $ 0\le \l_k < 1\}$ sei das von $\{\vek{a_{i_1}},\ldots,\vek{a_{i_d}}\}$ aufgespannte halboffene Parallelepiped und $E_F$ das von 
$\{\vek{a_{i_2}},\ldots,\vek{a_{i_d}}\}$ aufgespannte halboffene Parallelepiped.
Sei $u \in (\Z^n)^*$ ein  Normalenvektor an $F$ mit $\weite{u}{\s}=\weite{F}{\s}$.

Wir werden zeigen $\abs{E\cap \Z^n}=\weite{u}{\s}\cdot\abs{E_F\cap \Z^n}$.
Dazu beweisen wir zunächst folgendes Lemma:

\begin{Lemma}
Sei $\vek{v}\in\R^n$, sodass $(\aff(E)+\v)$ einen Gitterpunkt  enthält.

Dann gilt $\abs{E\cap \Z^n}=\abs{(E+\vek{v})\cap \Z^n}$.
\end{Lemma}
\begin{proof}
Wegen $\vek{0}\in E$ gilt $\lin(E)=\aff(E)$.
Betrachte zunächst den Fall $\vek{v}\in\aff(E)$. 
Da die Menge $\{\vek{a_{i_1}},\ldots, \vek{a_{i_d}}\}$ linear unabhängig ist, bildet sie eine Basis von $\lin(E)$. Zu $\vek{w}\in\lin(E)$ existieren also eindeutige Koeffizienten $\mu_k\in\R$ mit $\vek{w}= \sum_{k=1}^d \mu_k\vek{a_{i_k}}$.
Wir definieren nun eine Abbildung $\phi : (E+\vek{v})\cap\Z^n \to E\cap\Z^n$
gemäß 
\begin{equation*}
\phi(\vek{w}):=\sum_{k=1}^d \underbrace{(\mu_k - \lfloor \mu_k \rfloor)}_{=:\l_k}\vek{a_{i_k}}\:.
\end{equation*}
Wir werden zeigen, dass $\phi$ bijektiv ist, woraus folgt, dass $E$ und $E+\v$ gleichviele Gitterpunkte enthalten.

Seien $\a,\vek{b}\in (E+\v)\cap\Z^n$ mit $\a=\sum_{k=1}^d \mu_k^a\vek{a_{i_k}}$ und $\vek{b}=\sum_{k=1}^d \mu_k^b\vek{a_{i_k}}$. Sei $\phi(\a)=\phi(\vek{b})$, d.\,h. 
$\sum_{k=1}^d \l_k^a\vek{a_{i_k}}=\sum_{k=1}^d \l_k^b\vek{a_{i_k}}$.
Da die Menge $\{\vek{a_{i_1}},\ldots, \vek{a_{i_d}}\}$ linear unabhängig ist, folgt $\l_k^a=\l_k^b$ und damit $\mu_k^a-\mu_k^b\in\Z$ für alle $k$. Gleichzeitig gilt $\mu_k^a-\mu_k^b\in(-1,1)$,  wegen $\a,\vek{b}\in E+\v$. Daraus folgt $\mu_k^a=\mu_k^b$ für alle $k$ und damit $\a=\vek{b}$. Also ist $\phi$ injektiv.

Sei $\v=\sum_{k=1}^d \eta_k\vek{a_{i_k}}$ und $\a=\sum_{k=1}^d \l_k\vek{a_{i_k}}  \in E\cap\Z^n$. Definiere einen Punkt $\vek{b}=\sum_{k=1}^d \mu_k\vek{a_{i_k}}$, wobei

\begin{equation*}
 \mu_k:=\begin{cases}
         \lceil\eta_k\rceil + \l_k & \text{für } \lceil\eta_k\rceil-\eta_k+\l_k<1 \\
	 \lceil\eta_k\rceil  + \l_k  -1 & \text{für } \lceil\eta_k\rceil-\eta_k+\l_k\ge 1\quad.
        \end{cases}
\end{equation*}
Man kann sich leicht überzeugen, dass $\vek{b}$ ein Gitterpunkt in $E+\vek{v}$ ist, der von $\phi$ auf $\a$ abgebildet wird. Damit ist gezeigt, dass $\phi$ surjektiv ist.

Betrachte nun den allgemeinen Fall. 
Sei $\vek{w}$ ein Gitterpunkt in  $(\aff(E)+\v)$. Verschiebe  $E+\v$ um $-\vek{w}$. 
Dabei werden Gitterpunkte bijektiv auf Gitterpunkte abgebildet und es gilt $E+\v-\vek{w}\subseteq \aff(E)$. Damit haben wir den allgemeinen Fall auf den oben behandelten Spezialfall zurückgeführt.
$\,$
\end{proof}
Und nun zurück zum Beweis des Volumenlemmas.
Es gilt (s. auch Abb. \ref{figure:VolumenLemma}):
\begin{equation}
E \cap \Z^n =  \BigDisjUnion_{0\le i < \weite{u}{\s}} \left(\left(E_F + \frac{i}{\weite{u}{\s}} \vek{a_{i_1}}\right)   \cap \Z^n\right) 
\end{equation}
Die Inklusion \glqq$\supseteq$\grqq\ ist klar. Die Inklusion \glqq$\subseteq$\grqq\ folgt, da $u(\a)$ für jedes $\a\in E\cap \Z^n$ ganzzahlig sein muss. 

Nach obigem Lemma ist $\abs{\left(\left(E_F + \frac{i}{\weite{u}{\s}} \vek{a_{i_1}}\right)   \cap \Z^n\right)}=\abs{F \cap \Z^n}$ für alle $i$ und damit folgt:

\begin{equation*}
\abs{E \cap \Z^n}   
= \weite{u}{\s} \cdot  \abs{E_F \cap \Z^n}
=\weite{u}{\s} \cdot  \vol(F)
\end{equation*}
$\,$
\end{proof}

\begin{Lemma}[Facettenlemma]
\label{facettenlemma}
Sei $\A\subseteq \Z^n$ eine Punktkonfiguration mit Indexmenge $I$ und $J\subseteq I$.
$J$ habe Facettenweite 1 und sei $F \seite J$ eine Facette. Dann hat $F$ ebenfalls Facettenweite 1.
\end{Lemma}

\begin{proof}
Sei $G$ eine Facette von $F$. Dann ist $G=F \cap F'$ für eine Facette $F'\seite J$.
Sei $u\in(\Z^n)^*$ ein Normalenvektor an $F' \seite J$ mit $\weite{u}{J}=\weite{F'}{J}=1$. Dann ist $u$
auch Normalenvektor an $G$, aufgefasst als Facette von $F$. Also gilt:
\[
0<\weite{G}{F} \le \weite{u}{F} \le \weite{F'}{J} =  1
\]
Daraus folgt $\weite{G}{F}=1$.
$\,$
\end{proof}

\begin{Lemma}[Pullinglemma]
\label{pullinglemma}
Sei $\A$ eine Punktkonfiguration mit Indexmenge $I$ und $J=\{j_1,\ldots, j_m\}\subseteq \eck(I)$ eine (geordnete) Teilmenge der Ecken.

Dann haben die maximalen Zellen von $\pull(\A,J)$ die Form 
$\{j_1\} \cup \sigma$, wobei $\sigma$ in $\pull(F,J \cap F)$ enthalten ist für eine Facette  $F$  von $I$  mit $j_1 \not\in F$. 
\end{Lemma}

\begin{figure}[htbp]
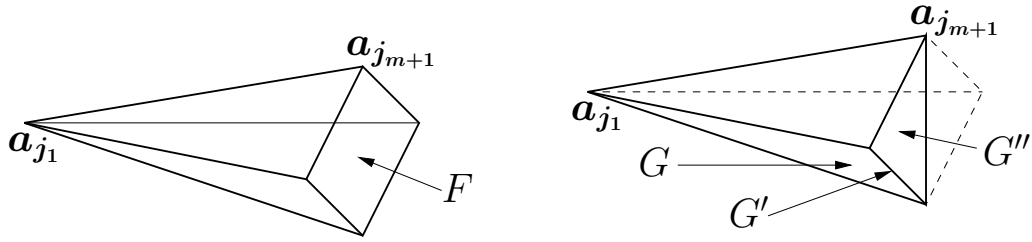


\scalebox{0.75}{
	\input{pullingLemma2komb.pstex_t}
}
\hfill
\scalebox{0.75}{
	\input{pullingLemma3komb.pstex_t}
}

\caption{Beispiel für den Induktionsschritt beim Pullinglemma. Links die Zelle $Z$, rechts $Z$'.
}
\label{figure:PullingLemmaBild}
\end{figure}

\begin{proof}

\IA{$m=1$}
 $\pull(\A,j_1)$ hat maximale Zellen $\{j_1\} \cup G$, wobei
 $G\in\pull(F,\emptyset)=\{\text{Seiten von $F$}\}$ für eine Facette $F$ von $\A$ mit $j_1\not\in F$.

\IS{$m\rightarrow m+1$} 
Sei $Z'\in\pull(\A,J \cup \{j_{m+1}\})$  eine maximale Zelle.
Diese ist in einer maximalen Zelle $Z$ von $\pull(\A,J)$ enthalten.
\begin{Fallunterscheidung}
\Fall{  $j_{m+1}\not\in Z$, d.h. $Z=Z'$. Fertig.
}
\Fall{$j_{m+1}\in Z$. 
Eine Skizze, die die Konstruktion in diesem Fall darstellt, findet man in Abbildung \ref{figure:PullingLemmaBild}.

Wendet man das Lemma
für $J=\{j_{m+1}\}$ auf $Z$ an, so erhält man
$Z'=\{j_{m+1}\} \cup G$
 für eine Facette $G$ von $Z$ mit $j_{m+1} \not\in G$. %

Es gilt nach Induktionsvoraussetzung $Z=\{j_1\} \cup F$ für eine Facette $F$ von $Z$ mit 
$j_{1}\not\in F$ und damit $j_{m+1}\in F$. Da $j_{m+1} \not\in G$ und $j_1\in Z'$ (wegen der Maximalität von $Z'$) folgt 
$j_1 \in G$ und damit
$G = \{j_1\} \cup G'$ für die Seite $G'=G\cap F$ von $Z$.

Damit erhalten wir $Z'=\{{j_{m+1}}\} \cup G = 
{\{{j_1},{j_{m+1}}\} \cup G'}
= {\{{j_1}\}\cup G''}$ für $G'':= {\{{j_{m+1}}\} \cup G'}$.

Wir müssen nun nur noch zeigen, dass $G''\in \pull(\A,(J  \cup \{j_{m+1}\})\cap H)$ gilt, wobei $H$ die Facette von $\A$ bezeichnet, die $F$ enthält. 
 Dies ist erfüllt, denn aus  $G'\seite F$ folgt $G'\in \pull(\A,J \cap H)$ und damit $G'' \in \pull(\A,(J \cup \{j_{m+1}\}) \cap H)$.

}
\end{Fallunterscheidung}

$\,$
\end{proof}

Nun können wir endlich Satz \ref{Satz:PullingIstTriangulierung} beweisen, d.\,h. 
 zeigen, dass Pullingtriangulierungen tatsächlich Triangulierungen sind:
\begin{proof}
\label{Beweis:PullingIstTriangulierung}
Induktion über $\dim(\A)$: Für $\dim(\A)\le 1$ ist die Aussage klar.

Sei $\dim(\A)\ge 2$:
Sei $J=\eck(I)=\{j_1,\ldots, j_m\}$ die (geordnete) Menge der Ecken von $\A$. Wir wissen aufgrund des Pul\-ling\-lem\-mas, dass die 
maximalen Zellen von $\pull(\A,J)$ die Form  ${\{j_1\} \cup \sigma}$
 haben, mit $\sigma \in \pull(F,J \cap F)$, wobei $F$ eine Facette von $I$ ist mit $j_1 \not\in F$. $\s$ ist nach Induktionsvoraussetzung ein Simplex, also auch ${\{{j_1}\}\cup\sigma}$.
$\,$
\end{proof}

\begin{Satz}[Pacos Lemma]
\label{Satz:PacosLemma}
Sei $\A=\{\vek{a_i}\,|\,i\in I\}$ eine Punktkonfiguration. 

Dann gilt:
\[
\text{$\A$ hat Facettenweite 1}
\Leftrightarrow
\text{Alle Pullingtriangulierungen von $\A$ sind unimodular.}
\]
\end{Satz}

\begin{proof}

\Hinrichtung 
Sei $\eck(I)=\{i_1,\ldots,i_m\}$ die (geordnete) Menge der Ecken von $I$. Betrachte die Pullingtriangulierung bezüglich dieser Eckenreihenfolge. 
Wir zeigen, dass alle Zellen unimodular sind per Induktion über die Dimension der Zellen. Nulldimensionale Zellen sind unimodular.

 Sei $Z$ eine $k$-dimensionale Zelle von
$\pull(\A,\eck(I))$. 
Aufgrund des Pullinglemmas %
wissen wir, dass $Z={\{i_1\} \cup\sigma}$ für $\sigma 
\in \pull(\A,\eck(I) \cap F)$ und eine Facette $F$ von $I$
 gilt. $F$ hat nach dem Facettenlemma 
 Facettenweite 1 und deshalb ist $\sigma$ nach Induktionsvoraussetzung ein unimodularer Simplex. Damit gilt:
\[
\vol(Z) \gleich{Lemma \ref{volumenlemma}}   \weite{F}{Z} \cdot\vol(\s)= 1 \cdot 1 = 1
\]

\Rueckrichtung
Angenommen $\A$ habe nicht Facettenweite 1. Dann gibt es eine Facette $F$ und eine Ecke $i$, sodass $\abs{ u(F) - u(\vek{a_{i}}) }  \ge 2$ für alle Normalenvektoren $u\in(\Z^n)^*$ von $F$ gilt. 
Ordne die Ecken von $\A$ so, dass zunächst $i$ kommt, dann die Ecken von $F$ und dann die restlichen Ecken von $\A$. Aus dem Pullinglemma folgt, %
dass es in der Pullingtriangulierung von $\A$ bezüglich dieser Ordnung der Ecken eine Zelle 
$Z = {\{i\}\cup\sigma}$ gibt, wobei $\sigma$ ein Zelle von $\pull(F,\eck(I) \cap F)$ ist.
Mit Hilfe des  Volumenlemmas %
folgt nun:
\begin{equation*}
\vol(Z) = \underbrace{\vol(\sigma)}_{=1}\cdot \underbrace{\weite{Z}{F}}_{\ge 2} \ge 2 \quad \text{\wid\, zu $Z$ unimodular.}
\end{equation*}
$\,$
\end{proof}

In Kapitel \ref{chapter:Gradschranken} werden wir sehen, 
dass alle Flusspolytope eine unimodulare Triangulierung haben. Dazu werden wir sie in Zellen mit Facettenweite 1 zerschneiden und dann Pacos Lemma anwenden.

\section{Reguläre Unterteilungen}
In diesem Abschnitt definieren wir eine andere wichtige Klasse von Unterteilungen, die regulären Unterteilungen, und zeigen, dass Pullingunterteilungen ein Spezialfall davon sind.

Die folgende Definition stammt aus \cite{sturmfelsGBCP}.
\begin{Definition}[Reguläre Unterteilungen]
Sei $\A=\{\vek{a_i}\,|\,i\in I\}$ eine Punktkonfiguration im $\Z^n$.  Ein Vektor $\omega\in\R^I$ induziert dann eine Unterteilung $\Delta_\omega$ von $\A$, gemäß $F={\{i_1,\ldots, i_r\}}
\in\Delta_\omega$ %
\begin{eqnarray}
\label{equation:RegTriA}
\Longleftrightarrow \exists\, \varphi_F : \R^n\rightarrow \R \text{ affin linear,  mit: }\quad
\varphi_F(\vek{a}_{i_j})&=&\omega_{i_j} \text{ für } i_j\in F \\
\varphi_F(\vek{a}_{i_j})  &<&\omega_{i_j} \text{ für } i_j\not\in F
\label{equation:RegTriB}
\end{eqnarray}

Eine Unterteilung $\Delta$ heißt \emph{regulär}, wenn ein $\omega\in\R^I$ existiert mit $\Delta=\Delta_\omega$.
\end{Definition}

\begin{Bemerkung}
 Ist $\omega$ hinreichend generisch, so ist $\Delta_\omega$ eine Triangulierung.
\end{Bemerkung}

\begin{Beispiel}[Reguläre Unterteilung]
$\,$\\
\parbox{9.5cm}{
Rechts sehen wir eine reguläre Unterteilung eines Rechtecks mit zwei maximalen Zellen. Die Ecken der Unterteilung sind die Ecken des Rechtecks. Der Gitterpunkt oben in der Mitte ist Teil der oberen maximalen Zelle, wohingegen der Gitterpunkt unten in der Mitte in keiner Zelle enthalten ist.
}
\hspace{0.5cm}
\parbox{5cm}{
\scalebox{0.7}{
	\input{RegulaereUnterteilung.pstex_t}
\\
$\,$
}
}

\end{Beispiel}

Wie bei den Pullingunterteilungen gibt es auch für die regulären Unterteilungen eine äquivalente geometrische Definition: Man bettet die Punktkonfiguration wieder in den $\R^{n+1}$ ein. Dazu bildet man den Punkt $\vek{a_i}$ ab auf den Punkt $(\vek{a_i},\o_i)\in \R^{n+1}$  und projiziert wieder die \glqq untere konvexe Hülle\grqq\ der Punkte zurück in den $\R^n$ und erhält so eine Unterteilung. Anhand dieser anschaulichen Definition ist relativ klar zu sehen, dass Pullingunterteilungen ein Spezialfall von regulären Unterteilungen sind.

Dies werden wir nun auch formal zeigen. Genau genommen wollen wir zeigen, dass  alle Pullingtriangulierungen regulär sind.

\begin{Satz}[Pullingverfeinerung regulär]
\label{Satz:PullingverfeinerungRegulaer}
Sei $\A=\{\vek{a_i}\,|\,i\in I\}$ eine Punktkonfiguration im $\Z^n$. Sei $i\in \eck(I)$ und $\Delta_\omega$ eine reguläre Unterteilung von 
$\A$. 

Dann ist $\pull(\Delta_\o, i)$ auch eine reguläre Unterteilung von $\A$.
\end{Satz}

\begin{proof}
\OBdA sei $I=[m]$ und $i=m$. Definiere $\o'\in\R^m$ gemäß $\o_j':=\o_j$ für $1\le j\le m-1$ und $\o_m':=\o_m- \eps$ für ein sehr kleines $\eps>0$.

\emph{Behauptung:} Dann gilt $\pull(\Delta_\o, m)=\Delta_{\o'}$, d.\,h. $\pull(\Delta_\o, m)$ ist regulär.

Für den Beweis führen wir folgende neue Bezeichnung ein:
Für eine reguläre Unterteilung $\Delta_\o$, eine Seite $\s\in\Delta_\o$ und eine affin lineare Funktion $\varphi$ sagen wir $\varphi$ \emph{rechtfertigt} $\sigma$ für $\Delta_\o$, genau dann, wenn
 $\varphi$ (\ref{equation:RegTriA}) und (\ref{equation:RegTriB}) für $\sigma$ und $\Delta_\o$
 erfüllt. 

\smallskip
\RechtsLinks Sei $\sigma\in\Delta_{\o'}$. 

\begin{Fallunterscheidung}
 \Fall{ $m\not\in\sigma$. Sei $\varphi_\sigma$ eine affin lineare Funktion, die 
 $\sigma$ für $\Delta_{\o'}$ rechtfertigt. Dann rechtfertigt $\varphi_\s$ auch $\sigma$ für $\Delta_\o$. Da $m \not\in\s$ folgt $\s\in\pull(\Delta_\o, m)$.
 }
 \Fall{$m\in\sigma$. Sei $\varphi_\sigma$ eine affin lineare Funktion, die 
 $\sigma$ für $\Delta_{\o'}$ rechtfertigt. $\varphi_\s$ rechtfertigt dann $\s\setminus \{m\}$ für $\Delta_\o$. %
Da wir $\eps$ sehr klein gewählt haben, liegt $m$ in einer Zelle $Z$ von $\Delta_\o$ mit  $\sigma\setminus \{m\}\seite Z$. 
 Nach Definition der Pullingtriangulierung folgt damit $\sigma={\{m\} \cup (\sigma\setminus \{m\}})\in
 \pull(\Delta_\o, m)$.
}
\end{Fallunterscheidung}

\medskip
\LinksRechts Sei $\sigma\in\pull(\Delta_\o, m)$.
\begin{Fallunterscheidung}
\Fall{ $m\not\in\sigma$. Also wurde $\sigma$ beim Ziehen an $m$ nicht verändert und es gilt $\sigma\in \Delta_\o$.
Es gibt also eine affin lineare Funktion $\varphi_\s$, die $\s$ für $\Delta_\o$ rechtfertigt. Insbesondere erfüllt diese Funktion $\varphi_\s(\vek{a_m})< \o_m$.
Da wir $\eps$ hinreichend klein gewählt haben gilt auch $\varphi_\s(\vek{a_m})<\o_m-\eps$
(für jedes $\tau\in\Delta$ existiert ein $\eps_\tau>0$ mit $\varphi_\tau(\vek{a_m})<\o_m-\eps_\tau$, wähle $\eps\le\min_{\tau\in\Delta_\o} \eps_\tau$).
}
\smallskip
\Fall{ $m\in\sigma$. 
Da $\pull(\Delta_\o, m)$ eine Verfeinerung von $\Delta_\o$ ist, existiert ein $\tau\in\Delta_\o$ mit $\s\subseteq \tau$. Sei $\varphi_\tau$ die affin lineare Funktion, die 
$\tau$ für $\Delta_\o$ rechtfertigt. Insbesondere gilt $\varphi_\tau(\vek{a_m})=\o_m$.

Nach Konstruktion der Pullingtriangulierung gilt $\sigma={\{m\} \cup F}$, wobei $F$ eine  Seite von $\tau$ ist. Folglich existiert eine affin lineare Funktion $\varphi_F : \R^n\rightarrow \R$, mit $\varphi_F(F)=0$, $\varphi_F(\tau\setminus F)<0$ und $\varphi_F(\vek{a_m})=-\eps$. 

Definiere nun $\varphi_\sigma=\varphi_\tau+\varphi_F$. Nach Konstruktion erfüllt $\varphi_\sigma$  für $\sigma$ und $\Delta_\o'$ (\ref{equation:RegTriA}). Es ist also noch zu zeigen, dass für $j\not\in\sigma$  gilt: $\varphi_\s(\vek{a_j})<\o_j$.

Für $j\in \tau\setminus (F \cup \{m\})$ gilt 
\begin{align*}
\varphi_\sigma(\vek{a_j})&=\underbrace{\varphi_\tau(\vek{a_j})}_{=\o_j}
+\underbrace{\varphi_F(\vek{a_j})}_{<0}<\o_j \\
\intertext{
und für $j\not\in \tau$ gilt
}
\varphi_\sigma(\vek{a_j})&=
\underbrace{\varphi_\tau(\vek{a_j})}_{<\o_j}+
\underbrace{\varphi_F(\vek{a_j})}_{\text{klein}}
\kleiner{$(*)$}\o_j,
\end{align*} 
wobei $(*)$ gilt, da wir $\eps$ sehr klein gewählt haben und damit auch $\varphi_F(\vek{a_j})$ sehr klein ist.
}
\end{Fallunterscheidung}
$\,$
\end{proof}

\begin{Korollar}
Sei $\A$ eine Punktkonfiguration im $\Z^n$ mit Indexmenge $I$. Dann ist jede Pullingtriangulierung $\Delta_{\mathrm{pull}}=\pull(\A,\eck(I))$ von $\A$ regulär.
\end{Korollar}

\begin{proof}
$\Delta_0$ ist regulär mit $\omega=\vek{0}$. $\Delta_{\mathrm{pull}}$ wird konstruiert, indem man nacheinander an allen Ecken zieht. In jedem Schritt erhalten wir nach obigem Satz wieder eine reguläre Unterteilung, sodass die Aussage per Induktion folgt.
$\,$
\end{proof}

\begin{Definition}[Hyperebenenverfeinerung]
\label{Definition:Hyperebenenverfeinerung}
Sei $\A=\{\vek{a_i}\,|\,i\in I\}$ eine Punktkonfiguration im $\Z^n$
und $\Delta$ eine Unterteilung von $\A$.
Sei $H$ eine affine Hyperebene, sodass für jedes $\sigma\in\Delta$ gilt:
$\conv_\A(\sigma)\cap H=\conv_\A(J)$ für ein $J\subseteq I$.\footnote{Das
bedeutet also, dass sich die beiden Hälften der Polytope aus $\Delta$, die zerschnitten werden, auch als konvexe Hülle der Punkte aus $\A$ darstellen lassen müssen. 

Verzichten wir auf diese Voraussetzung, so erhalten wir zwar eine Unterteilung im geometrischen Sinne, aber keine, die mit unserer kombinatorischen Definition verträglich ist.}

$\Delta'$ heißt \emph{Hyperebenenverfeinerung}
von $\Delta$ für die affine Hyperebene $H$, wenn gilt:
\begin{align*}
\sigma' \in \Delta' \Longleftrightarrow &\:\sigma'\subseteq \sigma \text{ für ein $\sigma\in\Delta$}\\
&\,\text{und }
 \left(\conv_\A(\s')=H^+\cap\conv_\A(\sigma) \text{ oder }
\conv_\A(\sigma')=H^- \cap\conv_\A(\sigma)\right)\\
\end{align*}
\end{Definition}

\begin{Satz}
\label{Satz:HyperebenenunterteilungRegulaer}
Sei $\A=\{\vek{a_i}\,|\,i\in I\}$ eine Punktkonfiguration im $\Z^n$, $\Delta_\o$ eine reguläre Unterteilung von $\A$ und $\Delta'$ eine Hyperebenenverfeinerung von $\Delta_\o$
bezüglich der affinen Hyperebene $H$. 

Dann ist $\Delta'$ ebenfalls regulär.
\end{Satz}

\begin{figure}[htbp]
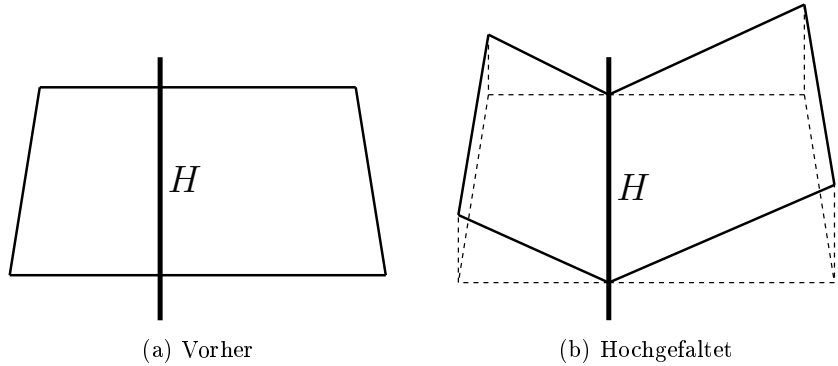

\centering
\subfloat[Vorher]
{

\scalebox{1}{
		\input{HyperplaneSubdivRegularA.pstex_t}
}

}
\hspace*{1mm}
\subfloat[Hochgefaltet]
{

\scalebox{1}{
		\input{HyperplaneSubdivRegularB.pstex_t}
	} 
}

\caption{Unsere Konstruktion entspricht in der geometrischen Definition dem Hochfalten des Polytops auf beiden Seiten von $H$, wobei in der Zeichnung vorher alle Punkte auf Höhe $0$ eingebettet waren.}
\label{figure:HyperplaneSubdivRegular}

\end{figure}

\begin{proof}
Sei  $H=\{\x\in\R^n\,|\,\psi(\x)= c\}$ für $\psi\in(\R^n)^*$ und $c\in\R$.

Wähle $\delta>0$ sehr klein und definiere $\omega'$ gemäß $\omega_j':=\omega_j + \delta|\psi(\vek{a_j})-c|$ für $j\in J$ (s. Abb. \ref{figure:HyperplaneSubdivRegular}).

Wie man leicht nachprüfen kann, gilt dann: $\Delta'=\Delta_{\omega'}$. 
$\,$
\end{proof}

Ein konkretes Beispiel befindet sich auf Seite \pageref{figure:RegulaereTriangulierungKonkret} in Abbildung \ref{figure:RegulaereTriangulierungKonkret}. Dort wird eine Punktkonfiguration zunächst durch eine Hyperebene unterteilt und anschließend wird an zwei Ecken gezogen, sodass wir eine (reguläre) Triangulierung erhalten.

\chapter{Algebra \DatumInKlammern}
In diesem Kapitel lernen wir einige algebraische Konzepte kennen. Wir lernen, was Gr\"obnerbasen und torische Ideale sind und stellen einen Zusammenhang zwischen regul\"aren unimodularen Triangulierungen einer Punktkonfiguration und den Gr\"obnerbasen des zugeh\"origen torischen Ideals her.
 
Die Hauptquellen sind das erste, vierte und achte Kapitel von \cite{sturmfelsGBCP} sowie das erste Kapitel von \cite{cloUsing}. Dort befinden sich auch Beweise von S\"atzen, die hier nur zitiert werden.

\section{Ideale und Gröbnerbasen}
In diesem Abschnitt definieren wir u.\,a. Termordnungen und Gröbnerbasen.

Sei $k$ ein Körper und $k[\vek{x}]:=k[x_1,\ldots,x_n]$ der Polynomring über $k$ in $n$ Variablen. Monome in $k[\vek{x}]$ werden mit $\x^\a:=x_1^{a_1}\cdots x_n^{a_n}$ bezeichnet.\footnote{In späteren Abschnitten werden wir unsere Variablen meistens mit einer anderen (geordneten) Menge als der Menge $[n]$ indizieren, z.\,B. der Indexmenge $I$ einer Punktkonfiguration $\A$. Um die Notation zu erleichtern bleiben wir in diesem Abschnitt aber bei der Indexmenge $[n]$. Da es eine Bijektion zwischen $I$ und $[\,\abs{I}\,]$ gibt ist dies auch völlig in Ordnung.}

Es existiert eine kanonische Bijektion zwischen $\N^n$ und den Monomen in $k[\x]$.
Dies erlaubt uns Gitterpunkte und Monome miteinander zu identifizieren. 
Beispielsweise entspricht $(3,0,2)\in\N^3$ dem Monom $x_1^3x_3^2\in k[x_1,x_2,x_3]$. 
Wir werden im weiteren Verlauf regelmäßig zwischen beiden Darstellungen hin- und herspringen.

Für einen Vektor $\v\in\Z^n$ bzw. ein Monom $\x^\v\in k[\x]$ definieren wir den \emph{Support} von $\v$ bzw. von $\x^\v$ als $\supp(\v)=\supp(\x^\v):=\{i\in [n]\,|\, v_i\not=0\}$.

Der \emph{Grad} eines Monomes $\x^\a$ ist definiert als $\deg(\x^\a):=\sum_{i=1}^n a_i$.
Für ein Erzeugendensystem $A$ eines Ideals $I\subseteq k[\x]$ definieren wir den Grad als $\deg(A):= \max_{f\in A}\deg(f)$. Ein Erzeugendensystem mit minimalem Grad nennen wir kurz \emph{minimales Erzeugendensystem}. Wir sagen, dass ein  Ideal $I$ im Grad $k$ erzeugt ist, genau dann, wenn der Grad eines minimalen Erzeugendensystems von $I$ höchstens $k$ ist.

\begin{Definition}[Termordnungen]
Eine Ordnung $\TO$ auf der Menge $\N^n$ heißt \emph{Term\-ord\-nung}, wenn je zwei Elemente vergleichbar sind, $\vek{0}$ das eindeutige Minimum ist und aus $\a\TO \vek{b}$ folgt $\a+\vek{c} \TO \vek{b}+\vek{c}$ für beliebiges $\vek{c}\in \N^n$.
\index{Termordnung}
\end{Definition}

\begin{Beispiel}[Termordnungen]
\begin{compactitem}
	\item Lexikographische Ordnung: 
	\index{Termordnung!Lexikographische Ordnung}
	\begin{align*}
	\vek{a}\lex\vek{b} \Leftrightarrow \exists k \in [n]:
	a_i=b_i \text{ für } i\le k-1 \text{ und } a_k < b_k \hspace*{130pt}
	\end{align*}
	\item Gradiert lexikographische Ordnung: 
	\begin{align*}
	\vek{a}\grlex\vek{b} \Leftrightarrow \deg(\vek{a})<\deg(\vek{b})\text{ oder } (\deg(\vek{a})=\deg(\vek{b}) \text{ und }\vek{a}\lex\vek{b}) \hspace*{70pt}
	\index{Termordnung!Gradiert lexikographische Ordnung}	
	\end{align*}
	\item Gradiert umgekehrt lexikographische Ordnung:
	\index{Termordnung!Gradiert umgekehrt lexikographische Ordnung}	
	\begin{align*}
	\vek{a}\grevlex\vek{b} \Leftrightarrow &\deg(\vek{a})<\deg(\vek{b}) \text{ oder }\\
	&\deg(\vek{a})=\deg(\vek{b}) \text{ und } \exists k\in [n]: a_i=b_i
	\text{ für } i\ge k+1 \text{ und } a_k>b_k	
	\end{align*}
\end{compactitem}
\end{Beispiel}

Entsprechend der Bijektion zwischen Gitterpunkten und Monomen in $k[\x]$ kann man den Begriff Termordnung natürlich auch auf die Monome in $k[\x]$ übertragen.
\smallskip

Sei $f\in k[\vek{x}]$. Es existiert eine eindeutige Punktkonfiguration $\A_f=\{\vek{a_i} \,|\,i\in I   \}\subseteq \N^n$, bei der alle Elemente paarweise verschieden sind und eine eindeutige Familie $\{\l_i\,|\,i\in I\} \subseteq \R$, sodass sich $f$ schreiben lässt als $f=\sum_{i\in I}\l_i \x^\vek{a_i}$.

Wir nennen das Polynom $f$ \emph{homogen} genau dann, wenn die Punktkonfiguration $\A_f$ homogen ist bezüglich der Funktion $\varphi\in(\Z^n)^*$, die alle Koordinaten aufsummiert.
\index{homogen!homogenes Polynom} 
 Homogene Polynome sind also solche, bei denen alle Monome den gleichen Grad haben.
Ein Ideal $I\subseteq k[\x]$ nennen wir \emph{homogen}, wenn es ein Erzeugendensystem von $I$ gibt, dass aus homogenen Polynomen besteht.
\index{homogen!homogenes Ideal}

Sei eine Termordnung $\TO$ vorgegeben und sei wieder $f\in k[\vek{x}]$ und $\A_f$ die Punktkonfiguration mit Indexmenge $I$, die die Exponenten der Monome von $f$ enthält.
Es gibt ein eindeutiges $k\in I$, sodass ${\vek{a_k}} $ maximal bezüglich der Termordnung $\TO$ ist.

Wir definieren nun den \emph{Leitterm} von $f$ als $\lt{f}{\TO}=\ltt{f}:= \lambda_k \x^{\vek{a_k}}$.
Zu einem Ideal $I\subseteq k[\vek{x}]$ definieren wir das \emph{Leitideal} von $I$ als
$\lt{I}{\TO}=\ltt{I}:=\{ \lt{f}{\TO} \,|\, f\in I\}$.

Die Monome von $I$, die nicht in $\lt{I}{\TO}$ enthalten sind, nennen wir \emph{Standardmonome}. 

Ein Erzeugendensystem $A$ für ein Ideal $I$ nennen wir \emph{quadratfrei}, wenn für jedes 
$f\in A$ gilt:  Der Exponent von $\ltt{f}$ liegt in $\{0,1\}^n$, d.\,h. im Leitterm tritt jede Variable höchstens einmal auf.    
\smallskip

Aus dem Hilbertschen Basissatz folgt, dass $k[\vek{x}]$ noethersch ist. Folglich sind alle Ideale $I\subseteq k[\vek{x}]$ endlich erzeugt. 
Sei $I$ ein Ideal und $A=\{g_1,\ldots,g_k\}$ ein Erzeugendensystem. Es gilt stets 
$\lt{I}{\TO}\supseteq \langle\lt{g_1}{\TO},\ldots,\lt{g_k}{\TO}\rangle$. Im Allgemeinen gilt aber keine Gleichheit. Dies führt uns zu der folgenden Definition:

\begin{Definition}[Gröbnerbasen]
Eine endliche Menge $\G=\{g_1,\ldots,g_k\}\subseteq I$ heißt \emph{Gröbnerbasis} von $I$ bezüglich $\TO$, genau dann, wenn die Leitterme der $g_i$ das Leitideal von $I$ erzeugen, d.\,h. wenn $\lt{I}{\TO}=\langle\lt{g_1}{\TO},\ldots,\lt{g_k}{\TO}\rangle$ gilt.

Dies ist äquivalent dazu, dass für jedes $f\in I$ ein $g\in \G$ existiert mit $\ltt{g}\,|\ltt{f}$.
\end{Definition}

Wenn für ein Ideal $I$ und $\G=\{g_1,\ldots,g_k\}\subseteq I$ gilt  
$\langle\ltt{g_1},\ldots,\ltt{g_k}\rangle=\ltt{I}$, so wird $I$ von $\G$ erzeugt, d.\,h. $\G$ ist
eine Gröbnerbasis von $I$. Dies folgt aus dem Divisionsalgorithmus für Polynome.

Eine Gröbnerbasis $\G=\{g_1,\ldots,g_k\}$ heißt \emph{minimal}, wenn für $i=1,\ldots, k$ gilt: $\ltt{g_i}\not\in\langle\ltt{g_1},\ldots,\widehat{\ltt{g_i}},\ldots,\ltt{g_k}\rangle$, d.\,h. kein Element von $\G$ ist redundant.
Eine Gröbnerbasis heißt \emph{reduziert}, wenn für zwei verschiedene Elemente $g,g'\in \G$ kein Term aus $g'$ durch $\lt{g}{\TO}$ teilbar ist.
Zu einem Ideal $I$ und einer Termordnung $\TO$ gibt es eine (bis auf skalare Vielfache) eindeutige reduzierte Gröbnerbasis. 
\smallskip

Sei $\omega=(\omega_1,\ldots,\omega_n)\in\R^n$. 
Mit $\,\cdot\,$ bezeichnen wir das kanonische Skalarprodukt auf dem $\R^n$.
Für ein Polynom $f=\sum_{i=1}^n c_i \vek{x}^{\vek{a_i}}$ definieren wir die Leit\-form 
\begin{equation}
\lt{f}{\omega}:=\sum_{i \in M} c_i \x^\vek{a_i},\; \text{ wobei } 
M:=\{i \in [n] \,|\, \omega\cdot \vek{a_i}=\max_{j\in[n]} (\omega\cdot  \vek{a_j}) \}\,.
\end{equation}

Ein Beispiel dafür ist
$\lt{x_1^3x_2x_3 + x_1x_2x_3^2+ x_1x_2^8}{(1,0,2)}=x_1^3x_2x_3 + x_1x_2x_3^2$.

Für ein Ideal $I$ definieren wir das Leitideal 
$\lt{I}{\omega}:=\langle \lt{f}{\omega} \,| \, f\in I \rangle$. 
Im Allgemeinen ist $\lt{I}{\omega}$ kein Monomideal (d.\,h. ein von Monomen erzeugtes Ideal).
Wenn man $\omega$ hinreichend generisch wählt, aber schon.

Sei $\omega\ge 0$ und $\TO$ eine beliebige Termordnung. Dann definieren wir folgendermaßen eine neue Termordnung $\TOO$:
\[
\vek{a}\TOO \vek{b} \Leftrightarrow \omega\cdot \vek{a} < \omega \cdot \vek{b} \text{ oder } (\omega\cdot\vek{a}=
\omega\cdot\vek{b} \text{ und } \vek{a}\TO\vek{b} )
\]

\begin{Satz}
Sei $\TO$ eine Termordnung und $I$ ein Ideal.

Dann existiert ein $\omega\in\N^n$ mit 
$\lt{I}{\omega}=\lt{I}{\TO}$. 
\end{Satz}
Wir sagen dann,  $\omega$ \emph{repräsentiert}  $\TO$ für $I$. 
Ein Beweis befindet sich in \cite{sturmfelsGBCP} auf Seite 4.

\section{Torische Ideale}
In diesem Abschnitt beschäftigen wir uns mit einer speziellen Klasse von Idealen im Ring $k[\vek{x}]$, den \emph{torischen Idealen}.

Sei $\A=\{\vek{a_i}\,|\,i\in I\}$ eine Punktkonfiguration im  $\Z^n$. Die Punkte $\vek{a_i} \in \Z^n$ identifizieren wir mit
Monomen $\vek{t}^{\vek{a}_i}$ im Laurentpolynomring $k[\vek{t}^{\pm1}]=k[t_1,\ldots, t_n,t_1^{-1},\ldots, t_n^{-1}]$.
Betrachte den Halbgruppenhomomorphismus 
\[
\pi : \N^I\rightarrow \Z^n,\qquad \vek{u}\mapsto \sum_{i\in I}u_i \vek{a_i}\quad.
\]
Das Bild von $\pi$ ist die Halbgruppe $\N\A=\{\sum_{i\in I}\lambda_i\vek{a_i} \,| \, 
\lambda_i\in \N\}$.

Die Abbildung $\pi$ lässt sich zu einem Homomorphismus von Halbgruppenalgebren hochheben:

\[
\hat{\pi} : k[\vek{x}]\rightarrow k[\vek{t}^{\pm 1}], \quad x_i\mapsto \vek{t}^{\vek{a_i}}
\]

\begin{Definition}[Torische Ideale]
 $I_{\A}:=\ker \hat\pi\subseteq k[\x]$ nennen wir das \emph{torische Ideal} von $\A$.

Für ein Gitterpolytop $P\subseteq\R^n$ und $\A:=P\cap\Z^n$ nennen $I_{\A}$ das torische Ideal von $P$.
\end{Definition}

Wie man leicht nachprüfen kann, ist $I_{\A}$ tatsächlich ein Ideal und sogar ein Primideal, da 
$k[\vek{t}^{\pm 1}]$ ein Integritätsbereich ist.

Zu $I_\A$ kann man ein geometrisches Objekt definieren, die \emph{affine torische Varietät} $V(I_\A):=\{\x\in k^n\,|\, \mbox{$f(\x)=0$}\; \mbox{$\forall f\in  I_\A$} \}$. Darum werden wir uns außer im nächsten Beispiel aber nicht weiter kümmern.
\begin{Beispiel}
Sei $\A=\{1,2\}$. Dann gilt $I_\A=\langle x^2-y \rangle$ und die zugehörige affine torische Varietät ist die Parabel $V(I_\A)=\{(x,y)\in k^2 \,|\, y=x^2 \}$.
\end{Beispiel}

Wir werden nun zeigen, dass jedes torische Ideal ein Erzeugendensystem hat, das aus Binomen besteht:
\begin{Satz}
\label{Satz:ESausBinomenFuerTI}
Das torische Ideal $I_\A$ wird als $k$-\!Vektorraum erzeugt von der Menge
\[
A:=
\{
\x^\u-\x^\v \,|\, \u,\v\in\N^n \text{ mit } \pi(\u)=\pi(\v)
\}\,.
\]
\end{Satz}

\begin{proof}
\OBdA sei $[n]$ die Indexmenge von $\A$.
Die Inklusion $A\subseteq I_\A$ ist klar. Wir müssen also nur noch zeigen, dass sich jedes Monom in  $I_\A$ als $k$-Linearkombination von Elementen der Menge $A$ darstellen lässt.

Angenommen nicht. Sei $f$ das Polynom aus $I_\A$, das sich nicht darstellen lässt, mit der Eigenschaft, dass  $\ltt{f}=\x^\u$ über alle diese Polynome minimal ist. Da $f\in I_\A$ ist, gilt $f(\vek{t}^{\vek{a_1}},
\ldots, \vek{t}^{\vek{a_n}})=0\in k[\vek{t}^{\pm 1}]$. Da der Term $\vek{t}^{\pi(\u)}$ als Summand in $f(\vek{t}^{\vek{a_1}},
\ldots, \vek{t}^{\vek{a_n}})$ auftritt, muss der gleiche Term auch ein weiteres mal auftreten (mit negativem Vorzeichen), d.\,h. in $f$ muss es ein Monom $\x^\v$ geben mit $\pi(\u)=\pi(\v)$. Da $\x^\u$ der Leitterm von $f$ ist gilt $\x^\v\TO\x^\u$. 
$f':=f-\x^\u+\x^\v$ ist dann  in $I_\A$ enthalten und lässt sich ebenfalls nicht als $k$-Linearkombination von Elementen aus $A$ darstellen, aber $\ltt{f'}\TO\ltt{f}$. \;\wid\: zu $\ltt{f}$ minimal.

$\,$
\end{proof}

 Aus obigem Satz folgt, dass 
$A:=$
$\{\mbox{$\x^\u-\x^\v$} \,|\,\u,\v\in\N^n \text{ mit } \pi(\u)=\pi(\v)\}$
das torische Ideal $I_\A$  
  erzeugt. Insbesondere hat $I_\A$ ein minimales Erzeugendensystem, das in $A$ enthalten ist und jede reduzierte Gröbnerbasis von $I_\A$ ist in $A$ enthalten.
\medskip

Torische Ideale haben wir oben algebraisch definiert als Kern des Al\-ge\-bra\-ho\-mo\-mor\-phis\-mus $\hat\pi$. Man kann torische Ideale aber auch rein kombinatorisch (bzw. geometrisch) sehen, und zwar als  Menge von Relationen zwischen den Punkten einer Punktkonfiguration $\A$.

Sei wie üblich $\A$ eine Punktkonfiguration mit Indexmenge $I$ und $k[\x]=k[x_i]_{i\in I}$ ein Polynomring. Wir identifizieren ein Monom $\x=x_i$ mit dem Punkt $\vek{a_i}$ und ein Monom $\x^\u$ mit der Summe $\sum_{i\in I} u_i \vek{a_i}$. 

Es gibt eine Eins-zu-Eins Beziehung zwischen Binomen in $I_\A$ und
Relationen zwischen den Gitterpunkten von $\A$ (eine Relation sind zwei Familien von Punkten aus $\A$, deren Summe gleich ist): Nach Definition ist
 $\x^\u-\x^\v$ in $I_\A$ enthalten genau dann, wenn $\sum_{i\in I} u_i \vek{a_i}=\sum_{i\in I} v_i \vek{a_i}$ gilt. Wie wir oben gesehen haben, haben torische Ideale ein Erzeugendensystem aus Binomen. Daher kann man $I_\A$ und die Menge der Relationen von $\A$ in gewisser Weise als (fast) gleich ansehen.

Die Begriffe Erzeugendensystem und Gröbnerbasis lassen sich auf diese Weise kombinatorisch 
definieren: Eine Menge $M$ von Relationen  von $\A$  bildet ein Erzeugendensystem, wenn sich alle Relationen von $\A$ durch Relationen in $M$ darstellen lassen.\footnote{
Ein Beispiel zum \glqq sich darstellen lassen\grqq:
Die Relation $a+b+c=d+e+f$ lässt sich darstellen durch
$a+b=g+e$ ($*$) und  $g+c=d+f$ ($**$), denn:
$a+b+c\gleich{$(*)$}
g+e+c \gleich{$(**)$} 
d+e+f
$.
}

Die Termordnung $\TO$ auf $k[\x]$ lässt sich auf kanonische Weise auf $\A$ übertragen 
($\vek{a_i}\TO\vek{a_j} \Leftrightarrow x_i\TO x_j$)
und so wird eine der beiden Familien jeder Relation zum Leitterm.
Eine Menge $\G$ von Relationen von $\A$ heißt Gröbnerbasis, wenn sie ein Erzeugendensystem  ist und alle Leitterme von Relationen von $\A$ Obermenge eines Leitterms einer Relation aus $\G$ sind.

Im \ref{chapter:Gradschranken}. Kapitel werden wir sowohl die kombinatorische, als auch die algebraische Darstellung verwenden, je nachdem welche gerade für uns günstiger ist.

\begin{Satz}
\label{Satz:SeiteTorischesIdeal}
Sei $\A=\{\vek{a_i} \,|\,i\in I\}$ eine homogene Punktkonfiguration im $\Z^n$, $J\seite I$ eine Seite und ${\cal B}:=\{\vek{a_j}\,|\, j\in J\}$ die Menge der Punkte in der Seite.  

 Seien $A$ und $B$ minimale Erzeugendensysteme von $I_\A$ und $I_{\cal B}$, die aus Binomen bestehen. Sei $a$ der Grad von $A$ und $b$ der Grad von $B$.

Dann gilt $I_{\cal B}\subseteq I_\A$ und $b\le a$.  
\end{Satz}

\begin{proof}
$I_{\cal B}\subseteq I_\A$ ist klar. Angenommen $b>a$.

 Falls $\cal B=\emptyset$, so ist die Aussage klar. Sei also $\cal B\not=\emptyset$.
\OBdA haben wir $B$ so gewählt, dass sich maximal viele Binome aus $I_\A$ durch Binome in $B$ vom Grad kleiner gleich $a$ darstellen lassen.

Nach unseren Voraussetzungen existiert ein Binom $\x^\u-\x^\v \in B$ vom Grad $b$, das sich nicht durch Binome kleineren Grades in $\cal B$ darstellen lässt, aber durch Binome aus $A$, die alle höchstens Grad $a$ haben:\footnote{Geometrisch ist klar, dass dies nicht geht. Wenn sich eine Relation, die in einer Seite von $\A$ liegt durch Relationen von $\A$ ausdrücken lässt, so müssen diese Relationen auch alle in dieser Seite von $\A$ liegen. }

 \begin{equation}
 \x^\u-\x^\v=\sum_{i=1}^k \x^{\vek{w^i}} \underbrace{(\x^{\vek{u^i}}-{\x^\vek{v^i}})}_{\in A}   
\label{equation:SeitenTeleskop}
 \end{equation}
\OBdA sei $\x^\u=\x^{\vek{w^1}}\x^{\vek{u^1}}$, $\x^\v=\x^{\vek{w^k}}\x^{\vek{v^k}}$
und $\x^{\vek{w^i}}\x^{\vek{v^i}}=\x^{\vek{w^{i+1}}}\x^{\vek{u^{i+1}}}$ für $i\in[k-1]$.
(\ref{equation:SeitenTeleskop}) ist also eine Teleskopsumme.

Nach Annahme muss ein $l$ existieren, sodass  nicht alle drei Vektoren  $\vek{w^l}, \vek{u^l}$ und $\vek{v^l}$ in $\cal B$ enthalten sind.  

Sei $l$ minimal mit dieser Eigenschaft. Betrachte nur den Fall $l\ge 2$. Der Fall $l=1$ geht analog.
Da $J\seite I$, existiert ein $\varphi \in (\R^n)^*$ und ein $c\in \R$ mit $\varphi(\vek{b})=c$ für alle $\vek{b}\in \cal B$ und $\varphi(\a) \ge c $ für alle $\a\in\A$. %

Wegen $\vek{w^{l-1}}+\vek{v^{l-1}}=\vek{w^{l}}+\vek{u^{l}}$
 folgt daraus:
\begin{align*}
c\cdot\left(\abs{\vek{w^{l-1}}}_1+\abs{\vek{v^{l-1}}}_1\right)&\gleich{\quad}\varphi\left(\sum_{i\in I} w_i^{l-1}\vek{a_i}\right) +
\varphi\left(\sum_{i\in I} v_i^{l-1}\vek{a_i}\right) 
\\
&\kleinergleich{($*$)}
\varphi\left(\sum_{i\in I} w^l_i\vek{a_i}\right) + \varphi\left(\sum_{i\in I} u^l_i 
\vek{a_i}\right)\\
&\kleinergleich{\quad}c \cdot \left(\abs{\vek{w^{l}}}_1+\abs{\vek{u^{l}}}_1\right)
\end{align*}

Da der erste und der letzte Term gleich, sind muss bei $(*)$ sogar Gleichheit gelten und wir erhalten 
$\vek{u^l},\vek{w^l}\in \cal B$. Analog folgt $\vek{v^l}\in\cal B$.
Also sind alle drei Vektoren in $\cal B$ enthalten.\:
\wid
$\,$
\end{proof}

Analog lässt sich beweisen, dass der Grad der reduzierten Gröbnerbasis von $I_{\cal B}$ kleiner gleich dem Grad der reduzierten Gröbnerbasis von $I_\A$ ist.

\begin{Satz}
\label{Satz:AdditivFolgtTorischeIdealeGleich}
Seien $\A\subseteq \Z^n$ und $\A'\subseteq \Z^m$ zwei homogene Punktkonfigurationen auf der gleichen Indexmenge $I$. 
Seien $H$ und $H'$ die von $\A$ und $\A'$ erzeugten additiven Halbgruppen.
Wenn eine bijektive Abbildung  $\phi : \A \to \A'$ existiert, die sich
 zu einem  Halb\-grup\-pen\-iso\-mor\-phis\-mus ${\psi} : H \to H'$ fortsetzen lässt,
 dann folgt $I_\A=I_{\A'}$.
\end{Satz}
\begin{proof}

$\sum_{i} \vek{a_i}=\sum_{i} \vek{b_i}$ ist eine Relation von Punkten aus $\A$, genau dann, wenn  
$\sum_{i} \psi(\vek{a_i})=\sum_{i} \psi(\vek{b_i})$ eine Relation von Punkten aus $\A'$ ist.
Damit folgt $I_\A \subseteq I_{\A'}$.

Andererseits ist
$\sum_{i} \vek{a_i'}=\sum_{i} \vek{b_i'}$ eine Relation von Punkten aus $\A'$, genau dann, wenn  
$\sum_{i} \psi^{-1}(\vek{a_i})=\sum_{i} \psi^{-1}(\vek{b_i})$ eine Relation von Punkten aus $\A$ ist.
Damit folgt $I_{\A'} \subseteq I_{\A}$.

$\,$
\end{proof}

\begin{Korollar}
\label{Korollar:HomogenePkTranslationsinvariant}
Sei $\A=\{\vek{a_i}\,|\,i\in I\}$ eine homogene Punktkonfiguration im $\Z^n$. 
Sei $\v\in \Z^n\setminus \aff(-\A)$ und $\A':=\{\vek{a_i}+\v\,|\,i\in I\}$ die um $\v$ verschobene Punktkonfiguration. 
Dann gilt $I_\A=I_{\A'}$.
\end{Korollar}

\begin{proof}
Wegen $\v\not\in\aff(-\A)$ ist $\A'$ weiterhin homogen, denn die Punkte von $\A'$ liegen weiterhin in einer Hyperebene, die nicht den Ursprung enthält. 

Sei $\phi : \A\to \A'$ die bijektive Abbildung, die $\vek{a_i}$ auf $\vek{a_i}+\v$ abbildet. 
Seien $H$ und $H'$ die von $\A$ bzw. $\A'$ erzeugten Halbgruppen.
Definiere $\psi : H\to H' $ gemäß $\psi(\sum_{l=1}^k\vek{a_{i_l}}):=\sum_{l=1}^k\vek{a_{i_l}}+k\cdot\v$. Wegen der Homogenität von $\A$ ist $\psi$ wohldefiniert. $\psi$ ist ein Halbgruppenisomorphismus und eine Fortsetzung von $\phi$ und damit folgt die Aussage aus dem eben bewiesenen Satz.
$\,$
\end{proof}

\begin{Korollar}
\label{Korollar:MinusAIdealGleich}
Sei $\A=\{\vek{a_i}\,|\,i\in I\}$ eine homogene Punktkonfiguration im $\Z^n$ und $\A':=-\A$.

Dann gilt $I_\A=I_{\A'}$.
\hfill $\Box$
\end{Korollar}

Wir werden nun noch einen Satz zitieren, der uns eine obere Schranke für den Grad von reduzierten Gröbnerbasen torischer Ideale liefert (s. \cite[Proposition 13.15]{sturmfelsGBCP}):

\begin{Satz}
Sei $\A$ eine homogene Punktkonfiguration im $\Z^d$ mit Indexmenge $I$ und sei $\TO$ eine Termordnung auf $k[\vek{x}]$,
für die
das Leitideal $\lt{I_\A}{\TO}$ quadratfrei ist.

Dann ist der Grad der reduzierten Gröbnerbasis von $I_\A$ bezüglich $\TO$ höchstens $d$.
\end{Satz}

\begin{Korollar}
\label{Korollar:GradschrankeTI}
Sei $\A$ die Menge der Gitterpunkte in einem (maximaldimensionalen) $(m\times n)$-Transportpolytop $\trans{r}{c}$. Dann gibt es eine Termordnung $\TO$, sodass der Grad der reduzierten Gröbnerbasis von $I_\A$ bezüglich $\TO$ 
kleiner gleich $\dim(A)+1=(m-1)(n-1)+1$ ist.
\end{Korollar}

\begin{proof}
$\trans{r}{c}$ ist gitteräquivalent zu einem homogenen Polytop $P\subseteq \R^{(m-1)(n-1)+1}$.

Um eine Gitteräquivalenz zu erhalten, projizieren wir $\trans{r}{c}$ zunächst in den Raum $\R^{(m-1)(n-1)}$, indem wir bei jeder Matrix die letzte Zeile und die letzte Spalte vergessen. Das so enstandene Polytop  ist gitteräquivalent zu
$\trans{r}{c}$, aber i.\,A. nicht homogen. Die Homogenität erreichen wir dadurch, dass wir das Polytop in Höhe Eins in den Raum $\R^{(m-1)(n-1)+1}$ einbetten, d.\,h. jeder Punkt im Polytop erhält eine zusätzliche letzte Koordinate mit dem Wert Eins.

Wie wir an die Termordnung $\TO$ kommen, sodass $\lt{I_\A}{\TO}$ quadratfrei ist, werden wir im nächsten Abschnitt sehen. 
$\,$
\end{proof}
Diese Schranke werden wir in Abschnitt \ref{section:SchrankeGB} ungefähr um den Faktor zwei verbessern.

\section{Gröbnerbasen und reguläre Triangulierungen}
In diesem Abschnitt stellen wir einen Zusammenhang zwischen 
diskreter Geometrie und Algebra her, genau genommen zwischen 
Triangulierungen einer Punktkonfiguration und dem torischen Ideal von dieser Punktkonfiguration.
Wir werden zeigen, dass die minimalen Nichtseiten einer regulären unimodularen Triangulierung den Leittermen einer Gröbnerbasis  entsprechen.

Dazu benötigen wir das folgende Lemma:
\begin{Lemma}
Sei $\A$ eine homogene Punktkonfiguration im $\Z^n$ mit Indexmenge $I$ und $\Delta$ eine unimodulare Triangulierung von $\A$.

Sei $\vek{b}\in\cone_\A(I)\cap \Z^n$. Dann existiert ein eindeutiges $\sigma\in\Delta$, sodass $\vek{b}$ im relativ Inneren des Kegels $\cone_\A(\s)$ liegt und dazu ein eindeutiger Vektor $\l \in \Z^I_{\ge 0}$ mit 
$\supp(\l)\subseteq \sigma$ und 
\[
\vek{b}=\sum_{i \in I}\l_i \vek{a_i}\,.
\]
\end{Lemma}

\begin{proof}
Da $\A=\BigDisjUnion_{\sigma\in\Delta} \relint(\conv_\A(\sigma))$
(s. Bemerkung \ref{Bemerkung:PolytopInneresTriangulierung})
 existiert ein eindeutiges $\sigma\in\Delta$, sodass  $\vek{b}$ im relativ Inneren des Kegels über $\cone_\A(\s)$ liegt. Daraus folgt, dass
ein Vektor
$\l \in \R^I_{\ge 0}$ mit 
$\supp(\l)\subseteq \sigma$ und $\vek{b}=\sum_{i\in I}\l_i \vek{a_i}$ existiert.
Wir müssen also lediglich noch die Ganzzahligkeit und die Eindeutigkeit von $\l$ zeigen.

Die Menge $\{\vek{a_i}\,|\,i\in\s\}$ ist affin unabhängig, da $\s$ ein Simplex ist.
Die Menge ist sogar linear unabhängig. Sei nämlich $\sum_{i\in\s}\mu_i\vek{a_i}=0$.
Wegen der Homogenität von $\A$ folgt $\sum_{i\in\s}\mu_i=0$ und wegen der affinen Unabhängigkeit folgt $\mu_i=0$ für alle $i\in \s$.
Aus der linearen Unabhängigkeit folgt die Eindeutigkeit des Vektors $\l$.

Betrachte den Punkt $\a=\sum_{i\in I} (\l_i-\lfloor \l_i \rfloor) \vek{a_i}$. Dies ist ein Gitterpunkt im von $\{\vek{a_i}\,|\,i\in\s\}$ aufgespannten Parallelepiped. Wegen $\vol(\s)=1$ folgt $\a=\vek{0}$. Damit folgt $\lfloor\l_i\rfloor=\l_i$ für alle $i\in I$, also $\l\in\Z^n_{\ge 0}$.
$\,$
\end{proof}

Sei $J\subseteq I$. Dann definieren wir
 $\vek{x}^J:=\prod_{j\in J}x_j$. Zu einer Nichtseite $F$ bezeichnen wir mit $s(F)$ den eindeutigen Vektor $\l$ aus obigem Lemma. 
 Für jede Nichtseite $F$ gilt also: $\x^F-\x^{s(F)}\in I_\A$.
 
 Wir werden nun die erstaunliche Tatsache beweisen, dass Binome dieser Form sogar eine Gröbnerbasis von $I_\A$ bilden. Bei dem folgenden Satz handelt es sich um einen Spezialfall der Korollare 8.4 und 8.9 aus \cite{sturmfelsGBCP}.
 
\begin{Satz}
\label{Satz:GroebnerbasisNS}
Sei $\A$ eine homogene Punktkonfiguration im $\Z^n$ und
$\Delta_\o$ eine reguläre unimodulare Triangulierung von $\A$.
Dann ist $\G_\Delta:=\{\x^F-\x^{s(F)}\,|\,\text{$F$ ist minimale }\linebreak[2]\text{Nichtseite von $\Delta$}\}$ eine reduzierte Gröbnerbasis von 
$I_\A$ bezüglich~$\TO_\o$.
\label{Satz:NSliefernInitialterme}
\end{Satz}

\begin{proof}
Wir wissen, dass $I_\A$ ein Erzeugendensystem aus Binomen hat (Satz \ref{Satz:ESausBinomenFuerTI}).
Sei $\x^\u-\x^\v$ ein Erzeuger von $I_\A$ 
 und $\x^\u$ der Leitterm. 
Es gilt also 
\begin{equation}
\sum_{i\in I} \omega_i u_i =\o\cdot\u \ge \o\cdot\v=\sum_{i\in I} \omega_i v_i\,.
\label{equation:uLT}
\end{equation}

Um zu sehen, dass $\G$ eine Gröbnerbasis von $I_\A$ ist, zeigen wir, dass eine minimale Nichtseite $F$ von $\Delta$ existiert mit $\x^F|\,\x^\u$.
\OBdA können wir voraussetzen, dass $\supp(\u) \cap \supp(\v)=\emptyset$ gilt.
\begin{Fallunterscheidung}
\Fall{$\supp(\u)$ ist Nichtseite:
Dann existiert eine minimale Nichtseite $F$ mit $F\subseteq \supp(\u)$ und folglich 
gilt $x^F|\,x^\u$.
}
\smallskip
\Fall{
$\supp(\u)$ ist Seite: Wir werden zeigen, dass dieser Fall nicht eintreten kann.

Da $\Delta$ eine reguläre Triangulierung ist, gibt es einen Vektor $\vek{c}$ mit
$\vek{a_j}\cdot \vek{c}=\omega_{j}$ für $j\in F\ (*)$ und $\vek{a_j}\cdot \vek{c}<\omega_{j}$ für  $j\not\in F\ (**)$.
Es gilt also:
\begin{align*}
 \vek{c}\cdot\left(\sum_{i\in I} u_i \vek{a_i}\right)
  \gleich{\hspace*{1cm}} \sum_{i\in I}  (\vek{c}\cdot \vek{a_i})u_i
 &\gleich{$(*)$}  \sum_{i\in I} \omega_i u_i \\ 
  &\groessergleich{(\ref{equation:uLT})}\sum_{i\in I} \omega_i v_i \\  
 &\groesser{$(**)$} \sum_{i\in I} (\vek{c}\cdot\vek{a_i})v_i \\
 &\gleich{\hspace*{1cm}} \vek{c}\cdot\left(\sum_{i\in I} v_i \vek{a_i}\right) 
\end{align*}

Nach Voraussetzung gilt aber $\sum_{i\in I} u_i \a_i = \sum_{i\in I} v_i \a_i$.\;\wid 

}
\end{Fallunterscheidung}
\smallskip

$\G$ ist sogar reduziert, denn $\supp(s(F))$ ist für beliebiges $F$ ein Simplex, enthält also keine Nichtseite. Also wird kein Term in $\G$ vom Initialterm eines anderen geteilt.

$\,$
\end{proof}

\begin{Korollar}
Torische Ideale von Flusspolytopen haben eine quadratfreie Gröbnerbasis.
\end{Korollar}

Wir wissen nun also, dass jede  homogene Punktkonfiguration  $\A$, die eine reguläre unimodulare Triangulierung besitzt, eine quadratfreie reduzierte Gröbnerbasis $\G$ hat. 

Diese tatsächlich auszurechnen ist nicht schwierig, wie das folgende Beispiel zeigt:

\begin{Beispiel}
Sei $\A$ die Punktkonfiguration aus Abbildung \ref{figure:RegulaereTriangulierungKonkret}.
Für das Ideal $I_\A$ können wir die reduzierte Gröbnerbasis $\G$ bezüglich der Termordnung, die vom Vektor $\o=(9,0,10,9,0,10)$ induziert wird, direkt aus Abbildung \ref{figure:RegulaereTriangulierungKonkretC}
ablesen.

Es gilt
$\G= \{ \underline{bf} - ae, \underline{ce}-bd, \underline{cf}-be, \underline{ad} - be,
\underline{ac} - b^2, \underline{df}-e^2 
\}$.
\end{Beispiel}

\begin{figure}[htbp]
\centering
\subfloat[Punktkonfiguration mit trivialer Unterteilung]
{
	\scalebox{0.76}{
		\input{RegulaereTriangulierung02.pstex_t}
	\label{figure:RegulaereTriangulierungKonkretA}
	}
}
\hfill
\subfloat[Hyperebenenunterteilung an der Hyperebene durch die beiden mittleren Punkte]
{
	\scalebox{0.76}{
		\input{RegulaereTriangulierungA2.pstex_t}

	}
	\label{figure:RegulaereTriangulierungKonkretB}
}
\hfill
\subfloat[Durch Ziehen an den Ecken unten links und oben rechts erhalten wir eine reguläre Triangulierung.]
{

\scalebox{0.76}{
		\input{RegulaereTriangulierungB.pstex_t}
		\label{figure:RegulaereTriangulierungKonkretC}
	} 
}
\caption{Eine reguläre Triangulierung einer Punktkonfiguration}
\label{figure:RegulaereTriangulierungKonkret}
\end{figure}

\chapter{Gradschranken f\"ur torische Ideale von Flusspolytopen \DatumInKlammern}
\label{chapter:Gradschranken}

In diesem Kapitel beweisen wir einige Gradschranken f\"ur Gr\"obnerbasen und Erzeugendensysteme von torischen Idealen von Transport- und Flusspolytopen. 

Unsere Ergebnisse sind im Einzelnen: Torische Ideale von Flusspolytopen sind im Grad drei erzeugt (Abschnitt \ref{section:GradDrei}).
Die reduzierte Gr\"obnerbasis von ($m\times n$)-Trans\-port\-poly\-to\-pen bez\"uglich einer beliebigen umgekehrt lexikographischen Termordnung hat h\"ochstens Grad $\left\lfloor\frac{m\cdot n}{2}\right\rfloor$ (Abschnitt \ref{section:SchrankeGB}) und es gibt Termordnungen und Transportpolytope, f\"ur die diese Schranke ann\"ahernd scharf ist (Abschnitt \ref{section:KonstruktionSchlechteGBs}). Glatte ($3\times 4$)-Transportpolytope sind sogar im Grad zwei erzeugt (Abschnitt \ref{section:3Kreuz4Glatt}).

\section{Die Zellunterteilungsmethode}
\label{section:Zellunterteilungsmethode}

In diesem Abschnitt beschreiben wir die Hauptmethode, die wir verwenden, um Gradschranken zu beweisen. Diese stammt aus \cite{christian-andreas-GBTP}.
Wir werden unsere Punktkonfigurationen entlang von affinen Hyperebenen der Form $\{\a \,|\,a_{i}=k\}$ (für $k\in\Z$) in Zellen schneiden und zeigen, dass es genügt die Grade der Erzeugendensysteme und  Gröbnerbasen der Zellen zu betrachten. 

Dies erleichert uns die Arbeit enorm. 
Zu einem fest gewählten Graphen $\vec{G}$ gibt es unendlich viele Flusspolytope. 
In diesen Flusspolytopen treten aber nur endlich viele verschiedene Zelltypen auf. 
Die Zellen haben eine einfachere Struktur als die Flusspolytope, da nach Translation alle Gitterpunkte aus den Zellen nur noch Einträge aus $\{0,1\}$ haben.

Für feste Graphen kann man also alle Zelltypen einzeln durchgehen. Wir werden dies in Abschnitt \ref{section:3Kreuz4Glatt} für  ($3\times 4$)-Transportpolytope tun.
\smallskip

Sei $F= F_{\vec{G},\vek{d},\vek{u},\vek{l}}$ ein Flusspolytop. %
Für $\vek{k}\in\Z^{\vec{E}}$
definieren wir eine \emph{Zelle von $F$} als:
\begin{equation}
Z_{F}(\vek{k}):=\{f\in F\,|\, k_e \le f(e) \le k_e+1   \}
\end{equation}

\begin{Bemerkung}[Zellen sind Gitterpolytope]
\label{Bemerkung:ZellenSindGitterpolytope}
 Die Zellen eines Flusspolytopes sind Flusspolytope mit anderen oberen und unteren Schranken und damit nach Satz \ref{Satz:FlusspolytopeSindGP} wieder Gitterpolytope.   
\end{Bemerkung}

Die verschobene Zelle $Z_{F}(\vek{k})-\vek{k}$ ist im Einheitswürfel $[0;1]^{\vec{E}}$ enthalten. Sie ist wieder ein Flusspolytop zum gleichen Graphen $\vec{G}$ mit neuem Bedarfsvektor $\vek{d'}$ und neuen Schranken $\vek{u'},\vek{l'}\in\{0;1\}^{\vec{E}}$. 
Die verschobene Zelle bezeichnen wir mit $Z_{\vek{d'}}=Z_{\vek{d'},\vek{u'},\vek{l'}}$. Wir sagen dann, dass die Zelle 
$Z_{F}(\vek{k})$ vom \emph{Typ} $Z_{\vek{d'},\vek{u'},\vek{l'}}$ ist.

 Nach Korollar \ref{Korollar:HomogenePkTranslationsinvariant} sind das torische Ideal von der Zelle $Z_F(\vek{k})$ und das von der verschobenen Zelle $Z_{\vek{d'},\vek{u'},\vek{l'}}$ gleich.

 \begin{Satz}
 \label{Satz:WannIstZelleVolldimensional}
 Sei $Z_\vek{d}\subseteq [0;1]^{\vec{E}}$ eine (verschobene) Zelle eines Flusspolytops zu einem Graphen $\vec{G}=(V,\vec{E})$.
 Dann gilt:
 \begin{enumerate}[(i)]
 \item Ist $Z_\vek{d}$ nicht leer, so ist $-\dout(v)\le d_v \le \din(v)$ für alle $v\in V$.
 \item Ist $Z_\vek{d}$ eine volldimensionale Zelle eines maximaldimensionalem Flusspolytops, so ist
$-\dout(v)+1\le d_v \le \din(v)-1$ für alle Ecken, die in einem ungerichteten Kreis von $\vec{G}$ enthalten sind.  
 \end{enumerate}
 
 \end{Satz}
 \begin{proof}
 Der erste Teil ist klar. 
Ist  $Z_\vek{d}$ volldimensional, so darf  einen Knoten $v\in V$ keine der beiden Schranken mit Gleichheit erfüllt sein. Andernfalls wäre nämlich für alle zu $v$ inzidenten Kanten $e$ der Wert von $f_e$ konstant für alle $\vek{f}\in Z_\vek{d}$. Damit wäre  $Z_\vek{d}$ nicht volldimensional.  Es muss also gelten $-\dout(v) < d_v < \din(v)$.
Wegen $\vek{d}\in\Z^V$ folgt daraus die Aussage.
$\,$
\end{proof}

Im Spezialfall, dass $F$ ein ($m\times n$)-Transportpolytop $\trans{r}{c}$ ist,  verwenden wir für die Zellen folgende Notation: Für eine Matrix $K\in \Z^{m\times n}$ schreiben wir 
$Z_{\vek{rc}}(K):=$ $\{M \in\trans{r}{c}\,|\,$ $ k_{ij} \le m_{ij} \le k_{ij}+1   \}$.

Verschiebt man die Zelle, so erhält man für  geeignete Vektoren $\vek{r'}$ und $\vek{c'}$: 

\begin{equation*}
Z_{\vek{rc}}(K)-K  = Z_{\vek{r-r'}\!,\vek{c-c'}}(0)=:\cell{\vek{r-r'}}{\vek{c-c'}}
\end{equation*}
Wir sagen dann, dass die Zelle $Z_{\vek{rc}}(K)$ vom Typ $\cell{\vek{r-r'}}{\vek{c-c'}}$ ist.

\medskip
Sei $\Delta_H$ die Unterteilung von $\A$, die wir erhalten, wenn wir mit der Unterteilung $\Delta_0=\{\text{Seiten von $\A$}\}$ beginnen und unsere Unterteilung  
entlang der Hyperebenen der Form $H_{e_i^*,k}$ für $k\in\Z$ (d.\,h. der affinen Hyperebenen, die senkrecht auf einem Einheitsvektor stehen und Gitterpunkte enthalten) verfeinern. Dass es sich dabei tatsächlich um  Hyperebenenverfeinerungen wie in Definition \ref{Definition:Hyperebenenverfeinerung} handelt (d.\,h. die Hyperebene schneidet Zellen nur in Seiten der Zellen), folgt aus 
Bemerkung \ref{Bemerkung:ZellenSindGitterpolytope}.

Die maximalen Zellen von $\Delta_H$ sind dann gerade die volldimensionalen Zellen von $F$.
Nach Satz \ref{Satz:HyperebenenunterteilungRegulaer} ist $\Delta_H$ regulär.

Verfeinert man die Unterteilung $\Delta_H$ nun zu einer regulären Triangulierung, so hat diese eine sehr angenehme Eigenschaft, wie der folgende Satz zeigt:

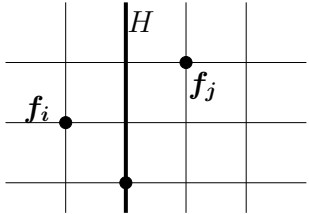
\vspace{0.1cm}

\begin{minipage}{4.5cm}

\scalebox{0.8}{
\input{MinimaleNS.pstex_t}
} 
\\
{Nicht in einer Zelle enthaltene, minimale Nichtseiten haben Kardinalität zwei.}
\label{figure:MinimaleNS}
\end{minipage}
\hfill
\begin{minipage}{8cm}
\begin{Satz}
Sei $F\subseteq \R^n$ ein Flusspolytop und $\A:=F\cap \Z^n$ die Menge der Gitterpunkte von $F$ sowie $\Delta$ eine reguläre Triangulierung von $\A$, die eine Verfeinerung der oben beschriebenen Hyperebenenunterteilung $\Delta_H$ ist.

Sei $\sigma$ eine minimale Nichtseite von $\Delta$, die in keiner Zelle von $F$ enthalten ist. 
Dann gilt $|\sigma|=2$.
\end{Satz}
\end{minipage}
\vspace{0.2cm}

\begin{proof}

Nach Voraussetzung muss es eine Hyperebene $H$  und zwei Gitterpunkte $\vek{f_i},\vek{f_j}\in\A$ mit $i,j\in \sigma$ geben, sodass 
$\vek{f_i}\in H^+\setminus H$ und $\vek{f_j}\in H^-\setminus H$, d.\,h. $\vek{f_i}$ und $\vek{f_j}$ liegen  auf verschiedenen Seiten von $H$ (s. Abb. oben). 

$\{i,j\}$ ist dann eine Nichtseite und wegen der Minimalität von $\sigma$ folgt $\sigma=
\{i,j\}$.

$\,$
\end{proof}

\begin{Korollar}
\label{Korollar:ZellenReichen}
Sei $F\subseteq \R^n$ ein Flusspolytop, $k\ge 2$ und $\Delta$ eine reguläre unimodulare Triangulierung von $F\cap \Z^n$, die eine Verfeinerung der oben beschriebenen Hyperebenenunterteilung $\Delta_H$ ist.
 
Wenn für alle Zellen $Z$ von $F$ gilt, dass das torische Ideal $I_Z$ im Grad $k$ oder kleiner erzeugt ist, so ist $I_F$ im Grad  $k$ oder kleiner erzeugt.

Zu einer Zelle $Z$ von $F$ bezeichnen wir die Einschränkung von $\Delta$ auf $Z\cap\Z^n$ mit $\Delta_Z$.%
Wenn für alle volldimensionalen Zellen $Z$ von $F$ gilt, dass die Gröbnerbasis $\G_{\Delta_Z}$ von $I_Z$  aus 
Satz \ref{Satz:GroebnerbasisNS} höchstens Grad $k$ hat, so hat die Gröbnerbasis 
$\G_\Delta$ von $I_F$  ebenfalls höchstens Grad $k$.
\end{Korollar}
 
\begin{proof}
Wir betrachten zunächst den Fall, dass für alle volldimensionalen Zellen $Z$ von $F$ die
Gröbnerbasen $\G_{\Delta_Z}$ höchstens Grad $k$ haben.

Nach Voraussetzung und Satz \ref{Satz:GroebnerbasisNS} haben dann alle minimalen Nichtseiten von $\Delta$ (sowohl die, die innerhalb einer Zelle von $F$ liegen, als auch die anderen) höchstens Kardinalität $k$. Deshalb wissen wir, wiederum wegen Satz \ref{Satz:GroebnerbasisNS}, dass 
$\G_\Delta=\mbox{$\{\x^F-\x^{s(F)}\,|\,$}\linebreak[2]\text{$F$ ist minimale}$ $\text{Nichtseite von $\Delta$}\}$ eine quadratfreie reduzierte Gröbnerbasis von $I_\A$ ist, deren Grad höchstens $k$ ist.

Im Fall, dass die torischen Ideale $I_Z$ von
 allen Zellen $Z$ lediglich ein Erzeugendensystem im Grad $k$ oder kleiner haben,  betrachten wir die Vereinigung dieser Erzeugendensysteme mit der Menge $\{\x^F-\x^{s(F)}\,|\,\text{$F$ ist minimale}$ Nichtseite von $\Delta$ und $F$ 
$\text{ist}$ $\text{nicht}$ $\text{in einer Zelle enthalten}\}$. Diese Menge erzeugt $\G_\Delta$ und damit auch das Ideal $I_F$.   
$\,$
\end{proof}

 Wenn wir Gradschranken für torische Ideale von Flusspolytopen zeigen wollen genügt es also, Gradschranken für alle auftretenden Zelltypen zu beweisen. Da alle Zellen Seiten von volldimensionalen Zellen sind, reicht es wegen Satz \ref{Satz:SeiteTorischesIdeal} aus, sich dabei auf volldimensionale Zellen zu beschränken.
 
\subsection*{Beispiele}

\begin{Beispiel}
\label{Beispiel1110_3333}
Betrachte das Transportpolytop $\transv{(1,1,10)}{(3,3,3,3)}$. Dieses ist isomorph zu 
$\transv{(1,1,6)}{(2,2,2,2)}$ (vgl. Beispiel \ref{Beispiel:IsoTPglattKomb}
auf Seite \pageref{Beispiel:IsoTPglattKomb}) und hat fünf volldimensionale Zellen.
Das Polytop enthält die folgenden 16 Gitterpunkte:

\begin{align*}
M_{1}&=\begin{dreiviermatrix}
1 & 0 & 0 & 0  \\
1 & 0 & 0 & 0  \\
1 & 3 & 3 & 3  \\
\end{dreiviermatrix}
&
M_{2}&=\begin{dreiviermatrix}
1 & 0 & 0 & 0  \\
0 & 1 & 0 & 0  \\
2 & 2 & 3 & 3  \\
\end{dreiviermatrix}
&
M_{3}&=\begin{dreiviermatrix}
1 & 0 & 0 & 0  \\
0 & 0 & 1 & 0  \\
2 & 3 & 2 & 3  \\
\end{dreiviermatrix}
&
M_{4}&=\begin{dreiviermatrix}
1 & 0 & 0 & 0  \\
0 & 0 & 0 & 1  \\
2 & 3 & 3 & 2  \\
\end{dreiviermatrix}
\displaybreak[1]
\\
M_{5}&=\begin{dreiviermatrix}
0 & 1 & 0 & 0  \\
1 & 0 & 0 & 0  \\
2 & 2 & 3 & 3  \\
\end{dreiviermatrix}
&
M_{6}&=\begin{dreiviermatrix}
0 & 1 & 0 & 0  \\
0 & 1 & 0 & 0  \\
3 & 1 & 3 & 3  \\
\end{dreiviermatrix}
&M_{7}&=\begin{dreiviermatrix}
0 & 1 & 0 & 0  \\
0 & 0 & 1 & 0  \\
3 & 2 & 2 & 3  \\
\end{dreiviermatrix}
&
M_{8}&=\begin{dreiviermatrix}
0 & 1 & 0 & 0  \\
0 & 0 & 0 & 1  \\
3 & 2 & 3 & 2  \\
\end{dreiviermatrix}
\displaybreak[3]
\\
M_{9}&=\begin{dreiviermatrix}
0 & 0 & 1 & 0  \\
1 & 0 & 0 & 0  \\
2 & 3 & 2 & 3  \\
\end{dreiviermatrix}
&
M_{10}&=\begin{dreiviermatrix}
0 & 0 & 1 & 0  \\
0 & 1 & 0 & 0  \\
3 & 2 & 2 & 3  \\
\end{dreiviermatrix}
&
M_{11}&=\begin{dreiviermatrix}
0 & 0 & 1 & 0  \\
0 & 0 & 1 & 0  \\
3 & 3 & 1 & 3  \\
\end{dreiviermatrix}
&
M_{12}&=\begin{dreiviermatrix}
0 & 0 & 1 & 0  \\
0 & 0 & 0 & 1  \\
3 & 3 & 2 & 2  \\
\end{dreiviermatrix}
\displaybreak[3]
\\
M_{13}&=\begin{dreiviermatrix}
0 & 0 & 0 & 1  \\
1 & 0 & 0 & 0  \\
2 & 3 & 3 & 2  \\
\end{dreiviermatrix}
&
M_{14}&=\begin{dreiviermatrix}
0 & 0 & 0 & 1  \\
0 & 1 & 0 & 0  \\
3 & 2 & 3 & 2  \\
\end{dreiviermatrix}
&
M_{15}&=\begin{dreiviermatrix}
0 & 0 & 0 & 1  \\
0 & 0 & 1 & 0  \\
3 & 3 & 2 & 2  \\
\end{dreiviermatrix}
&
M_{16}&=\begin{dreiviermatrix}
0 & 0 & 0 & 1  \\
0 & 0 & 0 & 1  \\
3 & 3 & 3 & 1  \\
\end{dreiviermatrix}
\end{align*}

\begin{align*}
K_\alpha&= \begin{dreiviermatrix}
	0 0 0 0 \\
	0 0 0 0 \\
	1 2 2 2 \\
     \end{dreiviermatrix}
\,
\leadsto
Z(K_\alpha)\cap \Z^{3\times 4}=\{M_1,M_2,M_3,M_4,M_5,M_9,M_{13}\}&\text{Typ: $\cell{1,1,3}{2,1,1,1}$}
\displaybreak[1]
\\
K_\b&= \begin{dreiviermatrix}
	0 0 0 0 \\
	0 0 0 0 \\
	2 1 2 2 \\
     \end{dreiviermatrix}
\,
\leadsto
Z(K_\b)\cap \Z^{3\times 4} = \{M_2,M_5,M_6,M_7,M_8,M_{10}, M_{14}  \}&\text{Typ: $\cell{1,1,3}{1,2,1,1}$}
\displaybreak[2]
\\
K_\gamma&= \begin{dreiviermatrix}
	0 0 0 0 \\
	0 0 0 0 \\
	2 2 1 2 \\
     \end{dreiviermatrix}
\,
\leadsto
Z(K_\gamma)\cap \Z^{3\times 4} = \{M_3,M_7,M_9,M_{10},M_{11},M_{12},M_{15}  \}&\text{Typ: $\cell{1,1,3}{1,1,2,1}$}
\displaybreak[1]
\\
K_\delta&= \begin{dreiviermatrix}
	0 0 0 0 \\
	0 0 0 0 \\
	2 2 2 1 \\
     \end{dreiviermatrix}
\,
\leadsto
Z(K_\delta)\cap \Z^{3\times 4} = \{M_4,M_8,M_{12},M_{13},M_{14},M_{15},M_{16}  \}&\text{Typ: $\cell{1,1,3}{1,1,1,2}$}
\\
K_\eps&= \begin{dreiviermatrix}
	0 0 0 0 \\
	0 0 0 0 \\
	2 2 2 2 \\
     \end{dreiviermatrix}
\,
\leadsto
Z(K_\eps)\cap \Z^{3\times 4} = \{M_2,M_3,M_4,M_5,M_7,M_8,M_{9},\\&\hspace{6.7cm}M_{10},M_{12},M_{13},M_{14},M_{15}  \}&\text{Typ: $\cell{1,1,2}{1,1,1,1}$}
\end{align*}

In den Beispielen am Ende von Abschnitt \ref{section:3Kreuz4Glatt} geben wir ein minimales Erzeugendensystem für das torische Ideal von diesem Transportpolytop an. 

\end{Beispiel}

\subsubsection*{Flusspolytope vom $\bm{K_4}$}
\begin{figure}[htbp]
\begin{center}
\begin{minipage}{11.4cm}
\begin{center}

	\input{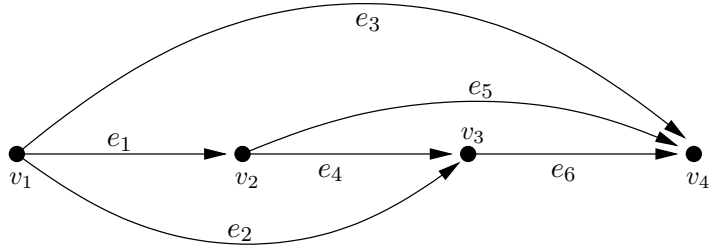}

\caption{Der $K_4$ mit der azyklischen Orientierung 
}
\label{figure:K4orientiert}
\end{center}

\end{minipage}
\end{center}
\end{figure}

Die Kanten des $K_4$ lassen sich (bis auf Isomorphie) auf vier verschiedene Arten orientieren. Genau eine dieser Orientierungen ist azyklisch (s. Abb. \ref{figure:K4orientiert}). Sei nun $\vec{G}$ der vollständige Graph auf vier Ecken mit der azyklischen Orientierung.

 Seien  $\vek{d}\in \Z^4$ und $\vek{u},\vek{l}\in\Z^6$ so gewählt, dass $ F= F_{\vec{G},\vek{d},\vek{u},\vek{l}}$ maximaldimensional, d.\,h. dreidimensional ist. Wir werden nun die Zellen des  Flusspolytops
 $F$ untersuchen und zeigen, dass $F$ eine Gröbnerbasis im Grad zwei hat.

Wir wollen zunächst herausfinden, welche Zelltypen in $F$ auftreten.
Sei $Z$ eine volldimensionale Zelle von $F$. Aus
Satz \ref{Satz:WannIstZelleVolldimensional} folgt:

\begin{displaymath}
\begin{array}{rclll}
-2 & \le & d_1 & \le & -1   \\ 
-1 & \le & d_2 & \le & 0 \\ 
0 & \le & d_3 & \le & 1 \\ 
 1 & \le & d_4 & \le & 2
\end{array}
\end{displaymath}
Es gilt also $d_1 \le d_2 \le d_3 \le d_4$. Da außerdem $d_1+d_2+d_3+d_4=0$ gelten muss, können  höchstens die folgenden sechs Zellen volldimensional sein: 
\begin{equation*}
Z_{(-2,-1,1,2)},\: 
Z_{(-2,0,1,1)},\:
Z_{(-2,0,0,2)}, \:
Z_{(-1,-1,1,1)}, \:
Z_{(-1,-1,0,2)},\:%
Z_{(-1,0,0,1)} 
\end{equation*}

Diese werden wir nun näher untersuchen. Wir bestimmen dazu die Gitterpunkte der Zellen:

\begin{align*}
Z_{(-2,-1,1,2)}\cap \Z^{\vec{E}}&=\{(1,1,0,1,1,1), (1,0,1,1,1,0), 
(0,1,1,0,1,0), (0,1,1,0,1,0)\} \\
Z_{(-2,0,1,1)}\cap \Z^{\vec{E}}&=\{(1,1,0,0,1,0), (1,1,0,1,0,1), 
(1,0,1,1,0,0),(0,1,1,0,0,0)\} \\
Z_{(-2,0,0,-2)}\cap \Z^{\vec{E}}&=\{(1,1,0,0,1,1), (1,0,1,1,0,1), 
(1,0,1,0,1,0), (0,1,1,0,0,1)\} \\
Z_{(-1,-1,1,1)}\cap \Z^{\vec{E}}&=\{(1,0,0,1,1,0), (0,1,0,0,1,0), 
(0,1,0,1,0,1), (0,0,1,1,0,0)\} \\
Z_{(-1,-1,0,2)}\cap \Z^{\vec{E}}&=\{(1,0,0,1,1,1), (0,1,0,0,1,1), 
(0,0,1,1,0,1), (0,0,1,0,1,0)\} \\
Z_{(-1,0,0,1)}\cap \Z^{\vec{E}}&=\{(1,0,0,1,0,1), (1,0,0,0,1,0), 
(0,1,0,0,0,1), (0,0,1,0,0,0)\} \\
\end{align*}

Wie man sieht, enthalten alle Zellen genau vier affin unabhängige Gitterpunkte. Es handelt sich also um dreidimensionale Simplexe und es gibt keine Relationen zwischen den Gitterpunkten der Zellen.

Mit Korollar \ref{Korollar:ZellenReichen} folgt daraus, dass Flusspolytope, die vom $K_4$ mit der azyklischen Orientierung stammen, eine Gröbnerbasis im Grad zwei haben.

\section{Eine scharfe obere Schranke für Erzeugendensysteme}
\label{section:GradDrei}
\subsection*{Schranke für Transportpolytope}
In diesem Abschnitt  beweisen wir die erstaunliche Tatsache, 
dass die torischen Ideale von allen Flusspolytopen im Grad drei  erzeugt sind.

Wir beweisen die Aussage zunächst in dem folgenden Satz für Transportpolytope:

\begin{Satz}
\label{Satz:TPsVonGradDreiErzeugt}
Torische Ideale von Transportpolytopen sind im Grad drei  erzeugt.
\end{Satz}

Wir zeigen dafür, dass die torischen Ideale von den  Zellen von ($m\times n$)-Transport\-poly\-to\-pen (in Zukunft ($m\times n$)-Zellen genannt)
im Grad drei erzeugt sind. Dies genügt nach Korollar \ref{Korollar:ZellenReichen}.
Wegen Satz \ref{Satz:SeiteTorischesIdeal} können wir uns dabei auf volldimensionale Zellen beschränken. 

Sei also $Z$ eine volldimensionale Zelle eines maximaldimensionalen ($m\times n$)-Trans\-port\-poly\-tops und sei $\A:=Z\cap\Z^{m\times n}$. Wir definieren nun eine 
Abstandsfunktion \mbox{$d: \A\times \A \to \N$}  gemäß $d(M,N):=\left|\{(i,j)\,|\,m_{ij}\not=n_{ij}\}\right|$. $d$~zählt also, an wievielen Stellen sich zwei Matrizen unterscheiden.
Für $\u,\v\in\R^\A$ mit $\x^\u-\x^\v\in I_\A$ definieren wir  
den Abstand $\bar d(\u,\v):=\min\{d(M,N)\,|\linebreak[2]\,M\in\supp(\u),N\in\supp(\v)\}$.\footnote{Hier und an einigen anderen Stellen in diesem Kapitel hat unsere Punktkonfiguration $\A$ keine Indexmenge. Anstatt der Indexmenge verwenden wir die Menge $\A$ selbst. Da alle Elemente von $\A$ paarweise verschieden sind, ist das in Ordnung.
}

\smallskip
Zur Einstimmung betrachten wir zunächst den ($2\times n$) Fall. Dort können wir mit der gleichen Technik wie im allgemeinen Fall, aber mit weniger Aufwand,  eine stärkere Aussage zeigen. 

\begin{Satz}
Torische Ideale von $(2\times n)$-Zellen sind im Grad zwei erzeugt. 
\end{Satz}

\begin{proof}
Sei wie üblich $\A$ die Menge der Gitterpunkte in der Zelle.
Angenommen, die Aussage ist falsch. Dann existiert ein Binom $\x^\u-\x^\v\in I_\A$ vom Grad $k\ge 3$, das sich nicht durch Binome kleineren Grades aus $I_\A$ darstellen lässt. 
Seien $\u$ und $\v$ so gewählt, dass $l:=\bar d(\u,\v)$ minimal ist über allen Binomen mit dieser Eigenschaft.

 Da unsere Zelle volldimensional ist, 
muss in jeder in $\A$ enthaltenen Matrix in jeder Spalte genau eine Eins und eine Null auftreten.
Seien $M^1\in \supp(\u)$ und $N^1\in\supp(\v)$ Matrizen mit $d(M^1,N^1)=l$. 
Nach Voraussetzung ist $l\ge 1$. Da $M^1$ und $N^1$ gleiche Zeilen- und Spaltensummen haben folgt automatisch $l\ge 4$.

Damit  können wir \oBdA voraussetzen, dass 
\begin{equation*}
M^1=
\begin{bmatrix}
1 & 0 &\multirow{2}{*}{\Large $M_*^1$}\\
0 & 1 &  
\end{bmatrix}
\quad\text{ und }\quad
N^1=
\begin{bmatrix}
0 & 1 &\multirow{2}{*}{\Large $N_*^1$}\\
1 & 0 &  
\end{bmatrix}
\end{equation*}
für geeignete ($m\times (n-2)$) Matrizen $M_*^1$ und $N_*^1$ gilt.

Sei $S:=\sum_{M\in\A} u_M M=\sum_{M\in\A} v_M M$. \OBdA sei $s_{11}\le s_{12}$ (andernfalls vertausche $\u$ und $\v$ sowie die ersten beiden Spalten der Matrizen). 
Dann muss  es eine Matrix $M^2\in \supp(\u)$ geben mit $m^2_{11}=0$ und $m^2_{12}=1$, also:
\begin{equation}
M^2=
\left[
\begin{matrix}
\hspace*{\Tabellenwegruecken}&\gc 0 & \gc 1 &\multirow{2}{*}{\Large $M_*^2$}\\
\hspace*{\Tabellenwegruecken}&1 & 0 &  
\end{matrix}\right]
\label{equation:2nMatrix}
\end{equation}

Es gilt also:
\begin{equation}
\label{equation:HilfsrelationzwoNZellen} 
M^1+M^2=
\begin{bmatrix}
1 & 0 &\multirow{2}{*}{\Large $M_*^1$}\\
0 & 1 &  
\end{bmatrix}
+
\begin{bmatrix}
0& 1 &\multirow{2}{*}{\Large $M_*^2$}\\
1 & 0 &  
\end{bmatrix}
=
\underbrace{
\begin{bmatrix}
0 & 1 &\multirow{2}{*}{\Large $M_*^1$}\\
1 & 0 &  
\end{bmatrix}
}_{:=\tilde{M}^1}
+
\underbrace{
\begin{bmatrix}
1 & 0 &\multirow{2}{*}{\Large $M_*^2$}\\
0 & 1 &  
\end{bmatrix}
}_{:=\tilde{M}^2}
\end{equation}

Sei $\vek{u'}:=\u+\vek{e_{\tilde{M}^1}} + \vek{e_{\tilde{M}^2}} - \vek{e_{{M}^1}} -\vek{e_{{M}^2}}$,
 d.\,h. $\x^{\vek{u'}}= 
\frac{ x_{\tilde{M}^1} \cdot x_{\tilde{M}^2}
}{
 x_{M^1}\cdot x_{M^2} 
}\x^\u$.

$\x^{\vek{u'}}-\x^\v$ ist ebenfalls ein Binom in $I_\A$ vom Grad $k$, das sich nicht durch Binome von kleinerem  Grad darstellen lässt 
(andernfalls ließe sich $\x^\u-\x^\v$ mit Hilfe dieser Binome und (\ref{equation:HilfsrelationzwoNZellen}) auch durch Binome kleineren Grades darstellen).

Es gilt aber $\bar d(\u',\v)=d(\tilde{M}^1,N^1)=\bar d(\u,\v)-4 < \bar d(\u,\v) $.\:  \wid 

$\,$
\end{proof}

Tatsächlich gilt sogar die stärkere Aussage, dass jedes $(2\times n)$-Transportpolytop  eine quadratfreie Gröbnerbasis im Grad zwei hat.
Ein $(2\times n)$-Transportpolytop $\transv{(r_1,r_2)}{(c_1,\ldots,c_n)}$ ist nämlich isomorph zum $r_1$-ten \emph{Hypersimplex} 
\begin{equation*}
H_{r_1}^n := \conv\left\{
\vek{e_{i_1}}+ \ldots + \vek{e_{i_{r_1}}} \,|\, i_1 < \ldots < i_{r_1} \right\}\subseteq \R^n\,.
\end{equation*}
Einen Isomorphismus liefert die Abbildung, die bei allen Matrizen aus dem Transportpolytop die zweite Zeile vergisst.
Nach \cite[Satz 14.2]{sturmfelsGBCP} haben die torischen Ideale von Hypersimplexen eine
quadratfreie Gröbnerbasis im Grad zwei.

\smallskip
Nun betrachten wir den ($m\times n$) Fall. Der Beweis ist von der Grundstruktur her ähnlich wie im ($2\times n$) Fall. Es gibt i.\,A. aber keine Matrix $M^2$, wie in Gleichung (\ref{equation:2nMatrix}), die die beiden Einsen von $\tilde{M}^1$ überdeckt, die  $M^1$ nicht überdeckt. Wir benötigen dafür zwei Matrizen $M^2$ und $M^3$ und können daher nur drei als Gradschranke zeigen. Darüberhinaus benötigen wir noch einige zusätzliche kombinatorische Tricks.

\begin{Satz}
\label{Satz:TransportGrad3}
Torische Ideale von ($m\times n$)-Zellen sind im Grad drei erzeugt. 
\end{Satz}

\begin{proof}
Wir beginnen ganz genau so wie im ($2\times n$) Fall: Sei $Z$ die betrachtete Zelle und $\A$ die Menge der Gitterpunkte der Zelle.  Wir nehmen an, die Aussage sei falsch. Dann gibt es  ein Binom $\x^\u-\x^\v\in I_\A$  vom Grad $k\ge 4$, das sich nicht durch Binome kleineren Grades aus $I_\A$ darstellen lässt. \OBdA seien $\u$ und $\v$ so gewählt, dass $l:=\bar d(\u,\v)$ minimal ist über allen Binomen mit dieser Eigenschaft.
$M^1\in \supp(\u)$ und $N^1\in\supp(\v)$ seien Matrizen mit $d(M^1,N^1)=l$. Wieder folgt $l\ge 4$.

Für den Beweis benötigen wir das folgende Lemma, das uns etwas über die Struktur von $M^1$ und $N^1$ verrät:

\begin{Lemma}  
\label{Lemma:ZentralZwo}
 
 Es kann keine Indizes $i_1,i_2, j_1,j_2$ geben mit 
 $m^1_{i_1j_1}=m^1_{i_2j_2}=n^1_{i_1j_2}=n^1_{i_2j_1}=1$ und
 $m^1_{i_1j_2}=m^1_{i_2j_1}=n^1_{i_1j_1}=n^1_{i_2j_2}=0$, d.\,h.
 wenn wir die zu diesen Indizes gehörenden Untermatrizen betrachten, so erhalten wir, 
dass nicht gleichzeitig
 \[
\bordermatrix[{[]}]
{
 & j_1 & j_2 \cr
i_1 & 1 & 0 \cr
i_2 & 0 & 1 \cr
}\:
 \text{ als Untermatrix von $M^1$ und }
\;\;\bordermatrix[{[]}]
{
 & j_1 & j_2 \cr
i_1 & 0 & 1 \cr
i_2 & 1 & 0 \cr
}\:
 \text{ als Untermatrix von $N^1$}
\]
  auftreten kann.

Die Aussage gilt immer noch, wenn wir in einer der beiden Untermatrizen eine Null durch eine Eins oder eine Eins durch eine Null ersetzen.
\end{Lemma}
\begin{proof}
Angenommen doch. Wir zeigen, dass man den Abstand von $M^1$ und $N^1$ dann weiter verkürzen kann.

Definiere eine Matrix $\tilde{M}$ gemäß:
\begin{equation}
\tilde{m}_{ij}:=
\begin{cases}
1 & \text{für $(i,j)=(i_1,j_2)$ oder $(i,j)=(i_2,j_1)$} \\
0 & \text{für $(i,j)=(i_1,j_1)$ oder $(i,j)=(i_2,j_2)$} \\
m_{ij} & \text{sonst} 
\end{cases}
\end{equation}

Es existieren  Matrizen $M^2$ und $M^3$ in $\supp(\u)$, sodass für $A:=M^2+M^3$ gilt:
$a_{i_1j_2}\ge \tilde{m}_{i_1j_2}=1$ und $a_{i_2j_1}\ge \tilde{m}_{i_2j_1}=1$.

Folglich ist $ M^1 + M^2 + M^3 - \tilde{M}^1 \ge \vek{0}$ und damit in $2\cdot Z$ enthalten. Nach dem verallgemeinerten Satz von Birkhoff und von Neumann (s. S. \pageref{Satz:VerallgemeinterterBvN}) 
 existieren also Matrizen $A^2$ und $A^3$ mit $M^1+M^2+M^3 = \tilde{M} + A^2 + A^3$.

Wähle $\vek{u'}$ so, dass $\x^{\vek{u'}}=
\frac{ x_{\tilde{M}} \cdot x_{A^2} \cdot x_{A^3}
}{
 x_{M^1}\cdot x_{M^2} \cdot x_{M^3} 
}\x^\u$
gilt.
Nach Voraussetzung muss eine Matrix $M_4\in\supp(\u)$ existieren mit $M_4\not=N$ für alle $N\in \supp(\v)$. Daraus folgt  $\x^{\u'}-\x^\v\not=0$.\,
$\x^{\vek{u'}}-\x^\v$ ist also ein Binom in $I_\A$ vom Grad $k$, das sich nicht durch Binome kleineren  Grades darstellen lässt, da sich sonst auch $\x^{\vek{u}}-\x^\v$ durch Binome
kleineren Grades darstellen ließe.

Es gilt $d(\tilde M,N^1) = d(M^1,N^1)-4$. 
Ersetzt man eine Null durch eine Eins oder andersherum wie im Lemma angegeben, so ist
$d(\tilde M,N^1) = d(M^1,N^1)-2$. 
Es gilt also stets
$\bar d(\u',\v)=d(\tilde M,N^1) < d(M^1,N^1)=\bar d(\u,\v) $.\:  \wid 

$\,$
\end{proof}

Nun können wir den Beweis von Satz \ref{Satz:TransportGrad3} beenden. $M^1$ und $N^1$ müssen sich unterscheiden. Es gibt also einen Eintrag, wo in $M^1$ eine Eins steht und in $N^1$ eine Null. Da beide die gleiche Spaltensumme haben, muss es in der gleichen Spalte einen Eintrag geben, wo dies umgekehrt ist. Wir können also  \oBdA schreiben:

\begin{equation}
M^1=
\begin{bmatrix}
\hspace*{\Tabellenwegruecken} &1 & \gc \ldots & \gc \ldots \\
\hspace*{\Tabellenwegruecken} &0 & \gc \ldots & \gc \ldots \\
\hspace*{\Tabellenwegruecken} &\gc & \gc \ldots & \gc \ldots  \\
\end{bmatrix}
\qquad
N^1=
\begin{bmatrix}
\hspace*{\Tabellenwegruecken} &0 & \gc \ldots & \gc \ldots \\
\hspace*{\Tabellenwegruecken} &1 & \gc \ldots & \gc \ldots \\
\hspace*{\Tabellenwegruecken} &\gc & \gc \ldots & \gc \ldots \\
\end{bmatrix}
\end{equation}
Da beide Matrizen die gleichen Zeilensummen haben, muss es in der ersten Zeile noch einen Eintrag geben, wo $N^1$ eine Eins hat und $M^1$ eine Null. Analog muss es in der zweiten Zeile einen Eintrag geben, wo $N^1$ eine Null hat und $M^1$ eine Eins. Wegen
Lemma \ref{Lemma:ZentralZwo} müssen diese Einträge in verschiedenen Spalten liegen.

Es gilt also \oBdA:

\begin{equation}
M^1=
\begin{bmatrix}
\hspace*{\Tabellenwegruecken} & 1 & \gc & \dc 0 & \gc \ldots \\
\hspace*{\Tabellenwegruecken} & 0 & \dc 1 &\gc & \gc \ldots \\
\hspace*{\Tabellenwegruecken} &\multicolumn{3}{c}{\gc \ldots}   \gc & \gc \ldots \\
\end{bmatrix}
\qquad
N^1=
\begin{bmatrix}
\hspace*{\Tabellenwegruecken} &0 & \gc & \dc 1 & \gc \ldots \\
\hspace*{\Tabellenwegruecken} &1 & \dc 0 & \gc & \gc \ldots \\
\hspace*{\Tabellenwegruecken} &\multicolumn{3}{c}{\gc \ldots}   \gc & \gc \ldots \\
\end{bmatrix}
\end{equation}

Wäre nun $n_{23}=0$, dann hätten wir eine nach Lemma \ref{Lemma:ZentralZwo} verbotene Untermatrix:
\begin{equation}
M^1=
\begin{bmatrix}
\hspace*{\Tabellenwegruecken} &\dc 1 & \gc & \dc 0 & \gc \ldots \\
\hspace*{\Tabellenwegruecken} &\dc 0 &  1 &\dc & \gc \ldots \\
\hspace*{\Tabellenwegruecken} &\multicolumn{3}{c}{\gc \ldots}   \gc & \gc \ldots \\
\end{bmatrix}
\qquad
N^1=
\begin{bmatrix}
\hspace*{\Tabellenwegruecken} &\dc 0 & \gc & \dc 1 & \gc \ldots \\
\hspace*{\Tabellenwegruecken} &\dc 1 & 0 & \dc 0 & \gc \ldots \\
\hspace*{\Tabellenwegruecken} &\multicolumn{3}{c}{\gc \ldots}   \gc & \gc \ldots \\
\end{bmatrix}
\end{equation}
Es muss also $n_{23}=1$ gelten.
Analog folgen die Werte für die anderen dunkel hinterlegten Einträge:
\begin{equation}
M^1=
\begin{bmatrix}
 \hspace*{\Tabellenwegruecken} & 1 & \dc 1 &  0 & \gc \ldots \\
 \hspace*{\Tabellenwegruecken} & 0 &  1 &\dc 0 & \gc \ldots \\
 \hspace*{\Tabellenwegruecken} & \multicolumn{3}{c}{\gc \ldots}   \gc & \gc \ldots \\
\end{bmatrix}
\qquad
N^1=
\begin{bmatrix}
\hspace*{\Tabellenwegruecken} & 0 & \dc  0 &  1 & \gc \ldots \\
\hspace*{\Tabellenwegruecken} &1 &  0 &  \dc 1 & \gc \ldots \\
\hspace*{\Tabellenwegruecken} &\multicolumn{3}{c}{\gc \ldots}   \gc & \gc \ldots \\
\end{bmatrix}
\end{equation}

Dieses Argument können wir nun erneut anwenden: Wegen der gleichen Zeilensummen  
muss es in der ersten Zeile einen weiteren Eintrag geben, wo $N^1$ eine Eins hat und $M^1$ eine Null, sowie einen weiteren Eintrag in einer anderen Spalte, wo $N^1$ eine Null hat und $M^1$ eine Eins. Wegen
Lemma \ref{Lemma:ZentralZwo} ist der darüber/darunter liegende Eintrag auch festgelegt und wir erhalten:

\begin{equation}
M^1=
\begin{bmatrix}
\hspace*{\Tabellenwegruecken} & 1 & 1 & 0 & \dc 1 & \dc 0 & \gc \ldots \\
\hspace*{\Tabellenwegruecken} & 0 & 1 & 0 & \dc 1 & \dc 0 & \gc \ldots \\
\hspace*{\Tabellenwegruecken} & \multicolumn{5}{c}{\gc \ldots}   \gc & \gc \ldots \\
\end{bmatrix}
\qquad
N^1=
\begin{bmatrix}
\hspace*{\Tabellenwegruecken} & 0 & 0 & 1 &\dc  0 &  \dc 1 & \gc \ldots \\
\hspace*{\Tabellenwegruecken} &1 &  0 & 1 & \dc 0 &  \dc 1 & \gc \ldots \\
\hspace*{\Tabellenwegruecken} &\multicolumn{5}{c}{\gc \ldots}   \gc & \gc \ldots \\
\end{bmatrix}
\end{equation}

Dieses Argument wenden wir nun immer weiter an. Irgendwann sind aber alle Spalten \glqq verbraucht\grqq. Es lässt sich also nicht vermeiden, dass eine nach Lemma 
\ref{Lemma:ZentralZwo} verbotene Untermatrix auftritt.\;\wid
$\,$
\end{proof}
Damit ist Satz \ref{Satz:TPsVonGradDreiErzeugt} bewiesen.

\subsection*{Verallgemeinerung auf Flusspolytope}
Wir werden nun zeigen, dass sich das Ergebnis aus dem vorigen Abschnitt auf Flusspolytope verallgemeinern lässt. Dies geschieht in zwei Schritten. Zunächst zeigen wir, dass die Aussage für bipartite Graphen gilt, bei denen alle Kanten von einer Farbklasse in die andere zeigen. Im zweiten Schritt reduzieren wir das Problem für beliebige Graphen auf den diesen Fall.

\begin{Satz}
Sei $\vec{G}$ ein bipartiter gerichteter Graph, bei dem alle Kanten von einer Farbklasse in die andere zeigen, und 
sei  $F$ ein Flusspolytop zum Graphen $\vec{G}$.

Dann ist das torische Ideal von $F$ im Grad drei erzeugt.
\end{Satz}
\begin{proof}
Wegen Korollar \ref{Korollar:ZellenReichen} genügt es zu zeigen, dass das torische Ideal von einer Zelle $Z$ des Flusspolytops im Grad drei erzeugt ist.

$\vec{G}$ ist ein Untergraph vom $\vec{K}_{m,n}$ für geeignete $m,n\in \N$ und $Z$ ist damit Seite einer Zelle $Z'$ des ($m\times n$)-Transportpolytops mit dem gleichen Bedarfsvektor. 

Wie wir gerade gezeigt haben, ist das torische Ideal $I_{Z'}$ im Grad drei  erzeugt und damit ist nach Satz \ref{Satz:SeiteTorischesIdeal} auch $I_Z$ im Grad drei  erzeugt. 
$\,$
\end{proof}

Diesen Satz verallgemeinern wir nun mittels einer Reduktion aus
\cite[21.6a]{schrijverCO} auf Zellen von beliebigen Flusspolytopen.

\begin{Satz}
Sei $F_{\vec{G}}=F_{\vec{G},\vek{d},\vek{u},\vek{l}}$ ein Flusspolytop.

Dann ist das torische Ideal $I_{F_{\vec{G}}}$ im Grad drei  erzeugt. 
\end{Satz}

\begin{proof}
Wir transformieren den Graphen $\vec{G}=(V,\vec{E})$ in einen bipartiten Graphen $\vec{G}'=(V',\vec{E}')$, bei dem alle Kanten von einer Farbklasse in die andere zeigen. 

Außerdem werden wir eine bijektive Abbildung $\phi : F_{\vec{G}} \cap \Z^{\vec{E}}
\to F_{\vec{G}'} \cap \Z^{\vec E'}$  zwischen den Gitterpunkten der beiden Flusspolytope 
angeben, die sich zu einem Isomorphismus zwischen den von den Gitterpunkten erzeugten Halbgruppen fortsetzen lässt.
Nach dem vorhergehenden Satz und Satz \ref{Satz:AdditivFolgtTorischeIdealeGleich}
auf Seite \pageref{Satz:AdditivFolgtTorischeIdealeGleich} sind wir dann fertig.

Wir transformieren unseren Graphen, indem wir jeden Knoten in zwei miteinander verbundene Knoten aufteilen. An den ersten der beiden Knoten werden alle eingehenden Kanten des ursprünglichen Knotens gehängt und an den zweiten alle ausgehenden.

Konkret geht das folgendermaßen:
Wir  teilen jeden Knoten $v\in V$ in zwei Knoten $v'$ und $v''$ auf.
Für jede Kante $(u,v)$ fügen wir eine Kante $(u',v'')$ ein, die die oberen und unteren Schranken der Kante $(u,v)$ erbt. Zusätzlich fügen wir für jedes $v\in V$ Kanten der Form $(v',v'')$ hinzu. 

Wir definieren $N:=\sum_{v\in V}|d_v|$ und setzen $d'_{v'}:=-N$, $d'_{v''}:=N+d_v$
sowie $u'_{(v',v'')}:=\infty$ und $l'_{(v',v'')}:=0$.
Für ein Beispiel siehe Abbildung \ref{figure:GraphtransformationBipartit}.

\begin{figure}[htbp]
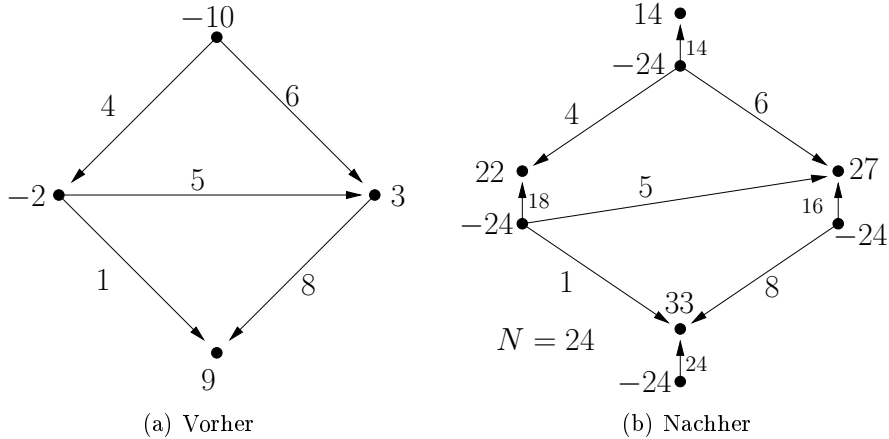

\centering
\subfloat[Vorher]
{

\scalebox{0.7}{
	\input{GraphMit4Ecken.pstex_t}
} 

}
\hspace*{2mm}
\subfloat[Nachher]
{

\scalebox{0.7}{
		\input{GraphMit4EckenBipartit.pstex_t}
}

}
\caption{Ein Graph $\vec{G}$ wird in einen bipartiten Graphen transformiert. Ein Bedarfsvektor und ein Fluss auf $\vec{G}$ werden mittransformiert.}
\label{figure:GraphtransformationBipartit}

\end{figure}

Zu einem ganzzahligen Fluss $f\in F_{\vec{G},\vek{d},\vek{u},\vek{l}}$ erhalten wir einen ganzzahligen Fluss $\phi(f)=f'\in  
F_{\vec{G}',\vek{d'},\vek{u'},\vek{l'}}$ gemäß: 

\begin{equation*}
f'(e')=
\begin{cases}
  f((u,v)) & \text{für $u,v\in V$, $u\not=v$ und $e'=(u',v'')$} \\
  N - \sum\limits_{\dout(e)=v}f(e)	& \text{für $v\in V$ und $e'=(v',v'')$} 
\end{cases}
\end{equation*}

Wie man leicht überprüfen kann, ist die so definierte Abbildung $\phi$ tatsächlich bijektiv und man kann sie zu einem Halbgruppenisomorphismus fortsetzen.

$\,$
\end{proof}

\section{Eine obere Schranke für Gröbnerbasen}
\label{section:SchrankeGB}
Wir zeigen in diesem Abschnitt, dass ($m\times n$)-Transportpolytope
reduzierte Gröbnerbasen im Grad $\left\lfloor\frac{m\cdot n}{2}\right\rfloor$ haben. 

Damit verbessern wir die allgemeine Gradschranke für Gröbnerbasen von torischen Idealen aus Korollar \ref{Korollar:GradschrankeTI} ungefähr um den Faktor zwei. %
Dazu verallgemeinern wir Theorem 14.8 aus
\cite{sturmfelsGBCP}.

\begin{Satz}[Gradschranke für Zellen]
Sei $Z$ eine Zelle eines ($m\times n$)-Transport\-poly\-tops $\trans{r}{c}$ mit $s:=\sum_i r_i=\sum_j c_j$. Sei $\A=\{\vek{a_i} \,|\, i\in I\}=Z\cap \Z^{m\times n}$ die Menge der Gitterpunkte in $Z$ und $\TO$ eine gradiert umgekehrt lexikographische Termordnung auf $k[\x]=k[x_{i}]_{i \in I}$.

Dann hat die reduzierte Gröbnerbasis $\cal G$ von $I_Z$ höchstens Grad $s$ und das Initialideal von $I_Z$ ist quadratfrei. 
\end{Satz}

\begin{proof}
\smallskip

Sei $\G$ die reduzierte Gröbnerbasis von $I_\A$ bezüglich $\TO$ und
sei $\vek{x}^\u-\vek{x}^\v \in \cal G$.
Es gilt %
$\supp(\u)\cap\supp(\v)=\emptyset$.
Sei $x_\rho$ die kleinste Variable, die $\vek{x}^\v$ teilt. Da unsere Termordnung $\TO$ umgekehrt lexikographisch ist, ist $x_\rho$ kleiner als jede Variable, die  $\vek{x}^\u$ teilt. 

Nach Voraussetzung gilt:
\begin{equation}
\hat\pi(\x^\u) = \vek{t}^{\sum_{i\in I} u_i \vek{a_i}} = \vek{t}^{\sum_{i\in I} v_i \vek{a_i}} = \hat\pi(\x^\v)
\label{torisch}
\end{equation}

$\hat\pi(x_\rho)=\vek{t}^\vek{a_\rho}$ teilt (\ref{torisch}) und damit auch $\vek{t}^{\sum_{i\in I} u_i \vek{a_i}}$. Also gibt es für jede Eins, die in der Matrix $\vek{a_\rho}$ auftritt, ein $i\in\supp(\u)$, sodass $\vek{a_i}$ an dieser Stelle auch eine Eins hat. Sei $J\subseteq \supp(\u)$ eine Menge minimaler Kardinalität, sodass $\sum_{j\in J} \vek{a_j} \ge \vek{a_\rho}$ gilt.

Sei $\vek{u'}\in \{0,1\}^I$ der Inzidenzvektor der Menge $J$. Dann ist $\vek{x}^\vek{u'}$ ein quadratfreies Monom, das höchstens Grad $s$ hat und es gilt $\left.\x^\vek{u'}|\, \x^\u\right.$.

Sei $M := \left(\sum_{i\in I} u'_i\vek{a_i}\right)-\vek{a_\rho}$. 
Die Matrix $M$ liegt im Polytop $k\cdot Z$ für ein $k\in \N$. Aus dem verallgemeinerten Satz von Birkhoff und von Neumann (s. S. \pageref{Satz:VerallgemeinterterBvN}) folgt nun,
dass sich $M$ als Summe von $k$  Matrizen aus $\A$ schreiben lässt, d.\,h. es existiert ein Vektor $\vek{v'}\in \N^I$  mit  $M=\sum_{i\in I}v_i\vek{a_i}$.

Also ist $\vek{x}^\vek{u'} - x_{\rho} \cdot \vek{x}^{\vek{v'}} \in I_{Z}$ und $\vek{x}^{\vek{u'}}$ ist der Leitterm, da $x_{\rho}$ kleiner ist als jede in $\x^\vek{u'}$ vorkommende Variable.

Da $\vek{x}^{\vek{u}}$ minimaler Erzeuger des Initialideales ist und ein Vielfaches von $\vek{x}^{\vek{u'}}$,
folgt $\vek{x}^{\vek{u}} = \vek{x}^{\vek{u'}}$. 
Also ist $\G$ eine quadratfreie Gröbnerbasis von $I_\A$, deren Grad höchstens $s$ ist.

$\,$
\end{proof}

\begin{Korollar}
Für $m,n\in\N$ beliebig gilt: Alle
($m\times n$)-Transportpolytope haben eine quadratfreie reduzierte Gröbnerbasis im Grad $\left\lfloor\frac{m\cdot n}{2}\right\rfloor$.
\end{Korollar}

\begin{proof}
Sei $\Delta_H$ die Hyperebenenunterteilung von $\trans{r}{c}\cap \Z^{m\times n}$. Diese verfeinern wir zu einer Triangulierung $\Delta$, indem wir an den Ecken aller Zellen ziehen.
Nach Satz \ref{Satz:PullingIstTriangulierung} erhalten wir so tatsächlich eine Triangulierung. Diese ist aufgrund von Pacos Lemma (Seite \pageref{Satz:PacosLemma})  unimodular und nach Satz \ref{Satz:PullingverfeinerungRegulaer} regulär. %

Wir zeigen nun, dass die Gröbnerbasis $\G_\Delta$ aus Satz \ref{Satz:GroebnerbasisNS} die gewünschte Eigenschaft hat.
Nach Korollar \ref{Korollar:ZellenReichen} genügt es dafür zu zeigen, dass 
für alle volldimensionalen Zellen $Z$ von $\trans{r}{c}$ die Gröbnerbasis $\G_{\Delta_Z}$
höchstens Grad $\left\lfloor\frac{m\cdot n}{2}\right\rfloor$ hat, wobei $\Delta_Z$ die Einschränkung von $\Delta$ auf $Z\cap \Z^{m\times n}$ bezeichnet.

Dies folgt aus dem eben bewiesenen Satz.
Jede Zelle $Z$ ist nämlich isomorph zu einer Zelle $Z'$, die die Eigenschaft hat, dass jeder Gitterpunkt $M\in Z'$ höchstens $\left\lfloor\frac{m\cdot n}{2}\right\rfloor$ viele Einträge hat, die verschieden von Null sind.
Es erfüllt  entweder $Z$ oder die komplementäre Zelle ($m_{ij} \mapsto 1-m_{ij}$) die Bedingung.
Aus den Korollaren \ref{Korollar:HomogenePkTranslationsinvariant} und \ref{Korollar:MinusAIdealGleich} folgt, dass die torischen Ideale der Zelle und ihres Komplements gleich sind.

$\,$
\end{proof}

\section{Konstruktion von Gröbnerbasen in hohem Grad}
\label{section:KonstruktionSchlechteGBs}
Wir zeigen in diesem Abschnitt, dass es glatte Transportpolytope und \glqq schlechte\grqq\ Termordnungen gibt, so dass die reduzierten Gröbnerbasen von den torischen Idealen dieser Transportpolytope bezüglich dieser Termordnungen einen hohen Grad haben.

Der Grad dieser Gröbnerbasen wird fast so groß sein wie die obere Schranke, die wir im vorigen Abschnitt bewiesen haben.

\begin{Satz}[Existenz von Gröbnerbasen in hohem Grad]

Seien $m$ und $n$ gerade. Dann existiert ein glattes $(m\times n)$-Transportpolytop $\trans{r}{c}$, so dass die reduzierte Gröbnerbasis $\G$ von $I_{\trans{r}{c}}$ mindestens Grad
$ \frac{m\cdot n}{2}-m=\frac{m(n-2)}{2}$
hat.
\end{Satz}

\begin{proof}

Zunächst definieren wir abkürzende Bezeichnungen für einige $\left(\frac m 2 \times \frac n 2\right)$-Ma\-tri\-zen mit Einträgen aus $\{0,1\}$. 
Mit $\Eins$ bzw. $\Null$ bezeichnen wir die Ma\-tri\-zen, bei der alle Einträge gleich Eins bzw.
Null sind.
 $\Einsij$ sei die Matrix, bei der der Eintrag $(i,j)$ gleich Eins ist und alle übrigen Null.
Vertauscht man in dieser Matrix die Rolle von Nullen und Einsen, so erhält man  $\Nullij:=\Eins-\Einsij$. 

Mit Hilfe dieser Matrizen definieren wir nun einige $\left(m  \times n \right)$-Matrizen:

\begin{align*}
A_{ij}&:=
\begin{bmatrix}
\hspace*{\Tabellenwegruecken} &\Einsi{1} & \gc \Nullij \\
\hspace*{\Tabellenwegruecken} &\gc\Nulli{1} & \Einsij \\
\end{bmatrix}
&
B_{ij}&:=
\begin{bmatrix}
\hspace*{\Tabellenwegruecken} &\Einsij & \gc \Nulli{1} \\
\hspace*{\Tabellenwegruecken} &\gc \Nullij & \Einsi{1} \\
\end{bmatrix}
\\[5mm]
E&:=
\begin{bmatrix}
\hspace*{\Tabellenwegruecken} &\Eins & \gc \Null \\
\hspace*{\Tabellenwegruecken} &\gc\Null & \Eins \\
\end{bmatrix}
&
D&:=
\begin{bmatrix}
\hspace*{\Tabellenwegruecken} &\Null & \gc \Eins \\
\hspace*{\Tabellenwegruecken} &\gc \Eins & \Null 
\end{bmatrix}
\end{align*}
\begin{align*}
C&:=
\begin{bmatrix}
\hspace*{\Tabellenwegruecken} &1 & 0 &\ldots & 0 & \gc 0 & \gc 1 &\gc \ldots & \gc 1 \\ 
\hspace*{\Tabellenwegruecken} &\vdots & & \ddots & & \gc \vdots & \gc & \ddots \gc & \gc\\
\hspace*{\Tabellenwegruecken} &1 & 0 & \ldots & 0 & \gc 0 & \gc  1 & \gc  \ldots & \gc 1\\
\hspace*{\Tabellenwegruecken} &\gc 0 & \gc 1 & \gc\ldots & \gc 1 &  1 &  0 &  \ldots & 0\\
\hspace*{\Tabellenwegruecken} &\gc \vdots & \gc  & \gc\ddots & \gc  &  \vdots &   &  \ddots & \\
\hspace*{\Tabellenwegruecken} &\gc 0 & \gc 1 & \gc\ldots & \gc 1 &  1 &  0 &  \ldots & 0\\
\end{bmatrix}
\end{align*}
\begin{align*}
\A^1&:=\left\{A_{ij} \,\middle|\, 1\le i \le \frac m 2, 2\le j\le \frac n 2 -1 \right\}
\cup \left\{B_{ij} \,\middle|\, 1\le i \le \frac m 2, 2\le j\le \frac n 2\right\}\\
\A^2&:=\left\{A_{ij} \,\middle|\, 1\le i \le \frac m 2, j= \frac n 2  \right\}
\end{align*}
In $\A^2$ sind also die Matrizen, die im unteren rechten Block in der letzten Spalte eine Eins haben, in $\A^1$ alle übrigen $A_{ij}$ und $B_{ij}$. 

Ein konkretes Beispiel für unsere Konstruktion befindet sich am Ende dieses Abschnitts.

\begin{Bemerkung}
Für jede Matrix  $A\in\A^1 \cup \A^2$ existieren eindeutige Indizes $(s,t)$ und $(k,l)$, sodass $A$ die einzige Matrix in $\A^1 \cup \A^2$ ist, 
mit $a_{st}=1$ und $a_{kl}=0$.

Genauer gesagt gilt $t=l=j+\frac n 2$ für $A_{ij}$ und für $B_{ij}$ gilt $t=l=j$. Die Matrizen aus $\A^2$ haben ihre eindeutige Null und Eins also in der letzten Spalte.
Die eindeutigen Einsen befinden sich alle im oberen linken oder im unteren rechten Block.
\label{Bemerkung:EindeutigNullEins}
\end{Bemerkung}

Es gilt:
{
\samepage
\ifthispageodd{}{\begin{addmargin}[-15mm]{0mm}}
\begin{align}
&\displaystyle
\qquad\sum\limits_{i=1}^{\frac m 2}\sum_{j=2}^{\frac n 2}(A_{ij} + B_{ij}) 
=  \left(\frac n 2-2\right)\cdot C + \left(\frac {m (n-2)-n} 2 +1\right)\cdot D + E= \label{equation:RelationVonGrossemGrad}
\\
\displaystyle
&
\begin{bmatrix}
\hspace*{\Tabellenwegruecken} & \frac n 2 -1 & 1 &\ldots & 1 & 
\gc \frac{m(n-2)-n} 2 + 1 & \gc \frac {m(n-2)} 2 -1 &\gc \ldots & \gc\frac{m(n-2)} 2 - 1 
\\ 
\hspace*{\Tabellenwegruecken} & \vdots & \vdots & \ddots & \vdots & \gc \vdots & \gc \vdots & \ddots \gc & \vdots\gc
\\
\hspace*{\Tabellenwegruecken} & \frac n 2 -1 & 1 & \ldots & 1 &
\gc \frac{m(n-2)-n} 2 + 1 & \gc  \frac {m(n-2)} 2 -1 & \gc  \ldots & \gc \frac {m(n-2)} 2 -1
\\
\hspace*{\Tabellenwegruecken} & \gc \frac{m(n-2)-n} 2 + 1 & \gc \frac {m(n-2)} 2 -1 & \gc\ldots & \gc \frac {m(n-2)} 2 -1 
&  \frac n 2 -1 &  1 &  \ldots & 1
\\
\hspace*{\Tabellenwegruecken} & \gc \vdots & \gc \vdots & \gc\ddots & \gc \vdots &  \vdots & \vdots  &  \ddots & \vdots
\\
\hspace*{\Tabellenwegruecken} & \gc \frac{m(n-2)-n} 2 + 1 & \gc \frac {m(n-2)} 2 -1 & \gc\ldots & \gc \frac {m(n-2)} 2 -1 
&  \frac n 2 -1 & 1 & \ldots & 1
\\
\end{bmatrix}
\nonumber
\end{align}
\ifthispageodd{}{\end{addmargin}}
}

Sei $\vek{r}=\left(\frac n 2+mn,\ldots,\frac n 2+mn \right)$ und $\vek{c}=\left(\frac m 2,\ldots,\frac m 2,m^2n + \frac m 2\right)$. Betrachte das zu diesen Vektoren gehörige Transportpolytop $\trans{r}{c}$. 
 Nach 
Satz \ref{Satz:KombinatorischGlatt} ist dieses Polytop glatt. Sei   $\A':=\trans{r}{c} \cap \Z^{m\times n}$.

(\ref{equation:RelationVonGrossemGrad}) ist eine Relation vom Grad $\frac {m( n-2)} 2 $ von Punkten der Punktkonfiguration $\A$, die man erhält, indem man $\A'$ um die Matrix
\[
N=
\begin{bmatrix}
0 & \ldots & 0 & -mn\\
\vdots & \ddots & \vdots  & \vdots \\
0 & \ldots & 0 & -mn
\end{bmatrix}
\]
verschiebt. Für $M\in \A$, $i\in[m]$ und $j\in[n-1]$ gilt $m_{ij}\ge 0$.

Wähle nun eine der 
\newcommand{\GRFAK}{\text{\Large!}}%
\begin{equation*}
\underbrace{\left(\frac {m(n-3)}{2}\right)}_{\abs{\A^1}}\GRFAK%
\underbrace{\bigg(\frac m 2\bigg)}_{\abs{\A^2}}\GRFAK
\underbrace{\left(\abs{\A}-\left(\frac{m(n-2)}{2}+1\right)\right)}_{\abs{\A \setminus (\A^1\cup \A^2 \cup \{E\})}
}\GRFAK
\end{equation*} %
vielen gradiert umgekehrt lexikographischen Termordnungen $\TO$ auf $\A$ mit der Eigenschaft, dass
$E \TO A^1 \TO A^2 \TO M$ für beliebige Matrizen $A^1 \in \A^1$,
$A^2 \in \A^2$ und 
$M\in \A\setminus (\A^1 \cup \A^2 \cup \{E\})$ gilt.

Sei $\x^\u-\x^\v\in I_\A$ das zur Relation (\ref{equation:RelationVonGrossemGrad}) gehörige Binom, d.\,h. $\u,\v\in \R^\A$ mit $u_M=1$ für $M\in\A^1\cup \A^2$, $v_C=\frac n 2 -2$,
$v_D=\frac{m (n-2)-n} 2 +1$, $v_E=1$ und alle anderen Werte sind Null.
$\x^\u$ ist der Leitterm bezüglich $\TO$, da $x_E \,|\, \x^\v$ und $E$  minimal ist bezüglich $\TO$.

Wir zeigen nun, dass $\x^{\u}-\x^{\v}$ in der reduzierten Gröbnerbasis $\G$ von $I_\A$ bezüglich $\TO$ liegt, 
indem wir zeigen, dass $\x^\u$ minimaler Erzeuger von $\lt{I_\A}{\TO}$ ist.

Angenommen, es existiert ein Binom $\x^{\u'}-\x^{\v'}\in I_\A$ mit Leitterm $\x^{\u'}, \u\not=\u'$ und $\x^{\u'}|\,\x^\u$, sowie $\supp(\vek{u'})\cap \supp(\vek{v'})=\emptyset$. 
Für $V:=\min(\supp(\vek{v'}))$ gilt dann $V\TO \min(\supp({\vek{u'}}))$.

Wir wollen nun einen Widerspruch herbeiführen. Dafür untersuchen wir, welche Werte $V$ annehmen kann.
\begin{Fallunterscheidung}
\Fall{$V=E$. Für jede Eins, die in $E$ auftritt, muss es eine Matrix in $\supp(\vek{u'})$ geben, die an dieser Stelle auch eine Eins hat. Mit 
Bemerkung \ref{Bemerkung:EindeutigNullEins} folgt dann aber $\u=\vek{u'}$.\: \wid
}
\Fall{$V\in \A^1$.
Nach Bemerkung \ref{Bemerkung:EindeutigNullEins} existiert ein eindeutiger
Index $(s,t)$ mit $t\le n-1$, sodass $V$ die einzige Matrix in $\A^1 \cup \A^2$ mit $v_{st}=1$ ist.
Für alle Matrizen $M$ aus $\supp(\vek{u'})$ gilt also $m_{st}=0$.
Für alle Matrizen $N\in\supp(\vek{v'})$ gilt $n_{st}\ge 0$. Damit kann aber nicht
 $\sum_{M\in \A} u_M' M = \sum_{M\in \A} v_M' M$ gelten, da die linke Summe an der Stelle $(s,t)$ kleiner ist als die rechte.

Also ist $\x^{\vek{u'}}-\x^{\vek{v'}}\not\in I_\A$.\:\wid 

}
\Fall{$V\in\A^2$, d.\,h. $V=A_{i\frac{n}{2}}$ für ein $i$. 
Aufgrund der Wahl unserer Termordnung folgt $\supp(\vek{u'})\subseteq \A^2$.
$V$ ist die einzige Matrix in $\A^2$ mit $v_{i1}=1$. 
Für alle Matrizen $M$ aus $\supp(\vek{u'})$ gilt also $m_{i1}=0$.
Für alle Matrizen $N\in\supp(\vek{v'})$ gilt $n_{i1}\ge 0$. Damit kann aber nicht
 $\sum_{M\in \A} u_M' M = \sum_{M\in \A} v_M' M$ gelten, da die linke Summe an der Stelle $(i,1)$ kleiner ist als die rechte.

Also ist $\x^{\vek{u'}}-\x^{\vek{v'}}\not\in I_\A$.\;\wid 
}
\end{Fallunterscheidung}

$\,$
\end{proof}

\begin{Beispiel}
Sei $n=m=6$. Wir betrachten also das Transportpolytop $\trans{r}{c}$ mit $\vek{r}=(39,39,39,39,39,39)$ und $\vek{c}=(3,3,3,3,3,219)$. Nach Translation sieht die Relation 
(\ref{equation:RelationVonGrossemGrad}) dann folgendermaßen aus: 
\begin{align*}
\underbrace{
\begin{sechssechsmatrix}
1 & \nix & \nix & \gc 1 & \gc \nix & \gc 1 \\
\nix & \nix & \nix & \gc 1 & \gc 1 & \gc 1 \\
\nix & \nix & \nix & \gc 1 & \gc 1 & \gc 1 \\
\gc \nix & \gc 1 & \gc 1 & \nix & 1 & \nix \\
\gc 1 & \gc 1 & \gc 1 & \nix & \nix & \nix \\
\gc 1 & \gc 1 & \gc 1 & \nix & \nix & \nix \\
\end{sechssechsmatrix}
}_{A_{12}\in\A^1}
+
\underbrace{
\begin{sechssechsmatrix}
1 & \nix & \nix & \gc 1 & \gc 1 & \gc \nix \\
\nix & \nix & \nix & \gc 1 & \gc 1 & \gc 1 \\
\nix & \nix & \nix & \gc 1 & \gc 1 & \gc 1 \\
\gc \nix & \gc 1 & \gc 1 & \nix & \nix & 1 \\
\gc 1 & \gc 1 & \gc 1 & \nix & \nix & \nix \\
\gc 1 & \gc 1 & \gc 1 & \nix & \nix & \nix \\
\end{sechssechsmatrix}
}_{A_{13}\in\A^2}
+
\underbrace{
\begin{sechssechsmatrix}
\nix & \nix & \nix & \gc 1 & \gc 1 & \gc 1 \\
1 & \nix & \nix & \gc 1 & \gc \nix & \gc 1 \\
\nix & \nix & \nix & \gc 1 & \gc 1 & \gc 1 \\
\gc 1 & \gc 1 & \gc 1 & \nix & \nix & \nix \\
\gc \nix & \gc 1 & \gc 1 & \nix & 1 & \nix \\
\gc 1 & \gc 1 & \gc 1 & \nix & \nix & \nix \\
\end{sechssechsmatrix}
}_{A_{22}\in\A^1}
+
\underbrace{
\begin{sechssechsmatrix}
\nix & \nix & \nix & \gc 1 & \gc 1 & \gc 1 \\
1 & \nix & \nix & \gc 1 & \gc 1 & \gc \nix \\
\nix & \nix & \nix & \gc 1 & \gc 1 & \gc 1 \\
\gc 1 & \gc 1 & \gc 1 & \nix & \nix & \nix \\
\gc \nix & \gc 1 & \gc 1 & \nix & \nix & 1 \\
\gc 1 & \gc 1 & \gc 1 & \nix & \nix & \nix \\
\end{sechssechsmatrix}
}_{A_{23}\in\A^2}
+
\underbrace{
\begin{sechssechsmatrix}
\nix & \nix & \nix & \gc 1 & \gc 1 & \gc 1 \\
\nix & \nix & \nix & \gc 1 & \gc 1 & \gc 1 \\
1 & \nix & \nix & \gc 1 & \gc \nix & \gc 1 \\
\gc 1 & \gc 1 & \gc 1 & \nix & \nix & \nix \\
\gc 1 & \gc 1 & \gc 1 & \nix & \nix & \nix \\
\gc \nix & \gc 1 & \gc 1 & \nix & 1 & \nix \\
\end{sechssechsmatrix}
}_{A_{32}\in\A^1}
+
\underbrace{
\begin{sechssechsmatrix}
\nix & \nix & \nix & \gc 1 & \gc 1 & \gc 1 \\
\nix & \nix & \nix & \gc 1 & \gc 1 & \gc 1 \\
1 & \nix & \nix & \gc 1 & \gc 1 & \gc \nix \\
\gc 1 & \gc 1 & \gc 1 & \nix & \nix & \nix \\
\gc 1 & \gc 1 & \gc 1 & \nix & \nix & \nix \\
\gc \nix & \gc 1 & \gc 1 & \nix & \nix & 1 \\
\end{sechssechsmatrix}
}_{A_{33}\in\A^2}
&
\displaybreak[1]
\\
+
\underbrace{
\begin{sechssechsmatrix}
\nix & 1 & \nix & \gc \nix & \gc 1 & \gc 1 \\
\nix & \nix & \nix & \gc 1 & \gc 1 & \gc 1 \\
\nix & \nix & \nix & \gc 1 & \gc 1 & \gc 1 \\
\gc 1 & \gc \nix & \gc 1 & 1 & \nix & \nix \\
\gc 1 & \gc 1 & \gc 1 & \nix & \nix & \nix \\
\gc 1 & \gc 1 & \gc 1 & \nix & \nix & \nix \\
\end{sechssechsmatrix}
}_{B_{12}\in\A^1}
+
\underbrace{
\begin{sechssechsmatrix}
\nix & \nix & 1 & \gc \nix & \gc 1 & \gc 1 \\
\nix & \nix & \nix & \gc 1 & \gc 1 & \gc 1 \\
\nix & \nix & \nix & \gc 1 & \gc 1 & \gc 1 \\
\gc 1 & \gc 1 & \gc \nix & 1 & \nix & \nix \\
\gc 1 & \gc 1 & \gc 1 & \nix & \nix & \nix \\
\gc 1 & \gc 1 & \gc 1 & \nix & \nix & \nix \\
\end{sechssechsmatrix}
}_{B_{13}\in\A^1}
+
\underbrace{
\begin{sechssechsmatrix}
\nix & \nix & \nix & \gc 1 & \gc 1 & \gc 1 \\
\nix & 1 & \nix & \gc \nix & \gc 1 & \gc 1 \\
\nix & \nix & \nix & \gc 1 & \gc 1 & \gc 1 \\
\gc 1 & \gc 1 & \gc 1 & \nix & \nix & \nix \\
\gc 1 & \gc \nix & \gc 1 & 1 & \nix & \nix \\
\gc 1 & \gc 1 & \gc 1 & \nix & \nix & \nix \\
\end{sechssechsmatrix}
}_{B_{22}\in\A^1}
+
\underbrace{
\begin{sechssechsmatrix}
\nix & \nix & \nix & \gc 1 & \gc 1 & \gc 1 \\
\nix & \nix & 1 & \gc \nix & \gc 1 & \gc 1 \\
\nix & \nix & \nix & \gc 1 & \gc 1 & \gc 1 \\
\gc 1 & \gc 1 & \gc 1 & \nix & \nix & \nix \\
\gc 1 & \gc 1 & \gc \nix & 1 & \nix & \nix \\
\gc 1 & \gc 1 & \gc 1 & \nix & \nix & \nix \\
\end{sechssechsmatrix}
}_{B_{23}\in\A^1}
+
\underbrace{
\begin{sechssechsmatrix}
\nix & \nix & \nix & \gc 1 & \gc 1 & \gc 1 \\
\nix & \nix & \nix & \gc 1 & \gc 1 & \gc 1 \\
\nix & 1 & \nix & \gc \nix & \gc 1 & \gc 1 \\
\gc 1 & \gc 1 & \gc 1 & \nix & \nix & \nix \\
\gc 1 & \gc 1 & \gc 1 & \nix & \nix & \nix \\
\gc 1 & \gc \nix & \gc 1 & 1 & \nix & \nix \\
\end{sechssechsmatrix}
}_{B_{32}\in\A^1}
+
\underbrace{
\begin{sechssechsmatrix}
\nix & \nix & \nix & \gc 1 & \gc 1 & \gc 1 \\
\nix & \nix & \nix & \gc 1 & \gc 1 & \gc 1 \\
\nix & \nix & 1 & \gc \nix & \gc 1 & \gc 1 \\
\gc 1 & \gc 1 & \gc 1 & \nix & \nix & \nix \\
\gc 1 & \gc 1 & \gc 1 & \nix & \nix & \nix \\
\gc 1 & \gc 1 & \gc \nix & 1 & \nix & \nix \\
\end{sechssechsmatrix}
}_{B_{33}\in\A^1}
&
\displaybreak[1]
\\
\quad=
\begin{sechssechsmatrix}
2 & 1 & 1 & \gc 10 & \gc 11 & \gc 11 \\
2 & 1 & 1 & \gc 10 & \gc 11 & \gc 11 \\
2 & 1 & 1 & \gc 10 & \gc 11 & \gc 11 \\
\gc 10 & \gc 11 & \gc 11 & 2 & 1 & 1 \\
\gc 10 & \gc 11 & \gc 11 & 2 & 1 & 1 \\
\gc 10 & \gc 11 & \gc 11 & 2 & 1 & 1 \\
\end{sechssechsmatrix}
=
1\cdot
\underbrace{
\begin{sechssechsmatrix}
1 & \nix & \nix & \gc  & \gc 1 & \gc 1 \\
1 & \nix & \nix & \gc  & \gc 1 & \gc 1 \\
1 & \nix & \nix & \gc  & \gc 1 & \gc 1 \\
\gc \nix & \gc 1 & \gc 1 & 1 & \nix & \nix \\
\gc \nix & \gc 1 & \gc 1 & 1 & \nix & \nix \\
\gc \nix & \gc 1 & \gc 1 & 1 & \nix & \nix \\
\end{sechssechsmatrix}
}_{C}
+
10\cdot
\underbrace{
\begin{sechssechsmatrix}
\nix & \nix & \nix & \gc 1 & \gc 1 & \gc 1 \\
\nix & \nix & \nix & \gc 1 & \gc 1 & \gc 1 \\
\nix & \nix & \nix & \gc 1 & \gc 1 & \gc 1 \\
\gc 1 & \gc 1 & \gc 1 & \nix & \nix & \nix \\
\gc 1 & \gc 1 & \gc 1 & \nix & \nix & \nix \\
\gc 1 & \gc 1 & \gc 1 & \nix & \nix & \nix \\
\end{sechssechsmatrix}
}_{D}
+
\underbrace{
\begin{sechssechsmatrix}
1 & 1 & 1 & \gc \nix & \gc \nix & \gc \nix \\
1 & 1 & 1 & \gc \nix & \gc \nix & \gc \nix \\
1 & 1 & 1 & \gc \nix & \gc \nix & \gc \nix \\
\gc \nix & \gc \nix & \gc \nix & 1 & 1 & 1 \\
\gc \nix & \gc \nix & \gc \nix & 1 & 1 & 1 \\
\gc \nix & \gc \nix & \gc \nix & 1 & 1 & 1 \\
\end{sechssechsmatrix}
}_{E%
}
\hspace*{1.3cm}  &
\end{align*}

\end{Beispiel}

\section{Glatte $(3\times 4)$-Transportpolytope}
\label{section:3Kreuz4Glatt}

In diesem Abschnitt zeigen wir, dass die torischen Ideale von glatten 
($3\times 4$)-Trans\-port\-poly\-to\-pen im Grad zwei erzeugt sind. 
In \cite{christian-andreas-GBTP} wurde auf ähnliche Weise die gleiche Aussage für glatte ($3\times 3$)-Transportpolytope gezeigt.
\medskip

 \OBdA seien alle  auftretenden Transportpolytope maximaldimensional. Nach Bemerkung 
\ref{Bemerkung:DimensionLieblingspolytope} sind 
 also alle auftretenden Transportpolytope  \mbox{$(m-1)(n-1)=6$} dimensional.
Wir verwenden die Zellunterteilungsmethode aus Abschnitt \ref{section:Zellunterteilungsmethode}.
Wir beschäftigen uns also mit Zellen der Form 
 $\cell{\vek{r}}{\vek{c}}$ mit $\vek{r}=(r_1,r_2,r_3)$ und $\vek{c}=(c_1,c_2,c_3,c_4)$. 
 
 Aus Satz \ref{Satz:WannIstZelleVolldimensional}
folgt, dass für volldimensionale Zellen $ 1 \le r_i \le 3$ und $ 1 \le c_j \le 2$ gelten muss. 
\OBdA gelte sogar $1 \le r_1 \le r_2 \le r_3 \le 3$ und  $1 \le c_1 \le c_2 \le c_3 \le c_4 \le 2$.
Damit können nur die in Tabelle \ref{table:AuftretendeZelltypen} aufgelisteten acht verschiedenen Zelltypen auftreten.

\begin{table}[htbp]
\begin{center}
\renewcommand{\arraystretch}{1.5} %
\begin{tabular}[h]{|r|l|}\hline
Summe aller Einträge & Zelltypen\\\hline
4 & $\cell{ 1,1,2 }{ 1,1,1,1 }$ \\
5 & $\cell{ 1,2,2 }{ 1,1,1,2 }$, $\cell{ 1,1,3 }{ 1,1,1,2 }$\\
6 & $\cell{ 2,2,2 }{ 1,1,2,2 }$, $\cell{ 1,2,3 }{ 1,1,2,2 }$\\
7 & $\cell{ 2,2,3 }{ 1,2,2,2 }$, $\cell{ 1,3,3 }{ 1,2,2,2 }$\\
8 & $\cell{ 2,3,3 }{ 2,2,2,2 }$\\\hline
\end{tabular}
\caption{In ($3\times 4$)-Transportpolytopen auftretende Zelltypen}
\label{table:AuftretendeZelltypen}
\end{center}
\end{table}

Unter diesen Zelltypen gibt es drei Paare von Zellen, die  
mittels der Abbildung $[m_{ij} \mapsto (1-m_{ij})]$ isomorph sind. 
Diese Paare sind
\[\cell{ 1,1,2 }{ 1,1,1,1 } \iso \cell{ 3,3,2 }{ 2,2,2,2 },
\quad
\cell{ 1,2,2 }{ 1,1,1,2 } \iso  \cell{ 3,2,2 }{ 2,2,2,1 }
\;\mbox{ und }\;
\cell{ 1,1,3 }{ 1,1,1,2 } \iso \cell{ 3,3,1 }{ 2,2,2,1 },\]
wobei rechts $\vek{r}$ und $\vek{c}$ jeweils absteigend sortiert sind, damit man besser sieht, dass die Zellen isomorph sind.

Aus den Korollaren \ref{Korollar:HomogenePkTranslationsinvariant} und \ref{Korollar:MinusAIdealGleich} folgt, dass die von den isomorphen Zellen erzeugten torischen Ideale jeweils gleich sind.

Es müssen also nur noch fünf Zelltypen untersucht werden.
Für diese haben wir mittels 4ti2 (\cite{4ti2}) minimale Erzeugendensysteme ausgerechnet (s. Tabelle \ref{table:Erzeuger3x4Zellen}).
 Zum Erzeugen der Eingabe und zur Weiterverarbeitung der Ausgabe wurden in der Skriptsprache
Perl (\cite{perl})
 geschriebene Skripte verwendet.

\renewcommand{\arraystretch}{1} %

\newlength{\erzeugerabsatzbreite}
\setlength{\erzeugerabsatzbreite}{10.5cm}
\newlength{\binombreite}
\setlength{\binombreite}{3.4cm}

\begin{table}[htbp]
\begin{tabularx}{14.17cm}{|r|c|lll|}\hline
Zelltyp & \#GP
& Erzeuger &  & \\\hline &&&&\\[-6pt]
$\cell{ 1,1,2 } { 1,1,1,1 }$ &
12
&
$ x_{5}x_{9}x_{11} - x_{6} x_{8} x_{12} $, &
$ x_{2}x_{9}x_{10} - x_{3} x_{7} x_{12} $, &
$x_{1}x_{6}x_{10} - x_{3} x_{4} x_{11} $, \hspace*{7pt}\\
&& $ x_{1}x_{5}x_{7} - x_{2} x_{4} x_{8} $, 
& $ x_{7}x_{11} - x_{8} x_{10} $, 
& $ x_{4}x_{12} - x_{5} x_{10} $, \\
&& $ x_{4}x_{9} - x_{6} x_{7} $, 
& $ x_{2}x_{6} - x_{3} x_{5} $, 
& $ x_{1}x_{12} - x_{2} x_{11} $, \\
&& $ x_{1}x_{9} - x_{3} x_{8} $ && 
\\[6pt]
$\cell{ 1,2,2 }{ 1,1,1,2 }$ & 12 
&$x_{1}x_{4}x_{5} - x_{2} x_{3} x_{6}$, 
&$x_{7}x_{12} - x_{8} x_{11}$, 
& $x_{7}x_{12} - x_{9} x_{10}$, \\ 
&& $x_{5}x_{10} - x_{6} x_{8}$,
& $x_{5}x_{11} - x_{6} x_{9}$,
& $x_{3}x_{10} - x_{4} x_{7}$, \\
&&$x_{3}x_{12} - x_{4} x_{9}$, 
& $x_{1}x_{8} - x_{2} x_{7}$, 
& $x_{1}x_{12} - x_{2} x_{11}$ 
\\[6pt]
$\cell{ 1,1,3 }{ 1,1,1,2 }$ &  
7
&
\multicolumn{2}{l}{Unimodularer Simplex, also Nullideal} &
\\[6pt]
$\cell{ 2,2,2 }{ 1,1,2,2 }$ &
15
&
$x_{11}x_{14} - x_{12} x_{13}$, 
& $x_{10}x_{15} - x_{11} x_{14}$, 
& $x_{8}x_{14} - x_{9} x_{10}$, \\
&& $x_{8}x_{15} - x_{9} x_{11}$, 
& $x_{6}x_{13} - x_{7} x_{10}$, 
&$x_{6}x_{15} - x_{7} x_{12}$, \\
&& $x_{4}x_{12} - x_{5} x_{10}$, 
& $x_{4}x_{15} - x_{5} x_{13}$, 
& $x_{2}x_{11} - x_{3} x_{10}$, \\
&& $x_{2}x_{15} - x_{3} x_{14}$, 
& $x_{1}x_{10} - x_{2} x_{8}$,
& $x_{1}x_{14} - x_{2} x_{9}$, \\
&& $x_{1}x_{11} - x_{3} x_{8}$, 
& $x_{1}x_{15} - x_{3} x_{9}$, 
& $x_{1}x_{10} - x_{4} x_{6}$, \\
&& $x_{1}x_{13} - x_{4} x_{7}$, 
& $x_{1}x_{12} - x_{5} x_{6}$, 
& $x_{1}x_{15} - x_{5} x_{7}$ 
\\
$\cell{ 1,2,3 }{ 1,1,2,2 }$ & 
8
&
$x_{3}x_{7} - x_{4} x_{6}$ && 
\\[3pt]\hline
\end{tabularx}
\caption{Erzeuger der Zellen von ($3\times 4$)-Transportpolytopen}
\label{table:Erzeuger3x4Zellen}
\end{table}

Wie man der Tabelle entnehmen kann, sind die torischen Ideale von drei der fünf  Zellen  im Grad zwei erzeugt. Die torischen Ideale der Zellen $\cell{ 1,1,2 } { 1,1,1,1 }$ und $\cell{ 1,2,2 }{ 1,1,1,2 }$ hingegen haben minimale Erzeuger vom Grad drei. 

Die beiden problematischen Zellen sind zwar Gitterpolytope, aber keine glatten Transportpolytope. Zu jeder Zelle muss es also immer noch 
mindestens eine weitere benachbarte Zelle im Polytop geben. Wir werden zeigen, dass man stets in einer 
solchen benachbarten Zelle Gitterpunkte finden kann, mit deren Hilfe man die  Relationen
vom Grad drei durch Relationen vom Grad zwei darstellen kann. Diese nennen wir \emph{Retter}.

\begin{Definition}
Sei $\trans{r}{c}$ ein ($3\times 4$)-Transportpolytop und sei $\A$ die Menge der Gitterpunkte von $\trans{r}{c}$. Seien $A,B,C,D,E,F \in \A$  und sei $A + B + C = D + E + F$ eine Relation vom Grad drei von $\A$. 
\begin{itemize}
\item $R \in \A$ heißt \emph{1-Retter}, wenn gilt: 
\begin{equation}
R + A = E + F
\label{1-RetterA}
\end{equation}
\begin{equation}
R + D = B + C
\label{1-RetterB}
\end{equation}
\item  Ein Tupel $(R_1,R_2,R_3)$ mit $R_1,R_2,R_3 \in \A$ heißt \emph{3-Retter}, wenn gilt:
\begin{equation}
B+C = R_1 + R_2
\label{3-RetterA}
\end{equation}
\begin{equation}
A + R_1 = D + R_3
\label{3-RetterB}
\end{equation}
\begin{equation}
R_2 + R_3 = E + F
\label{3-RetterC}
\end{equation}
\end{itemize}
\end{Definition}

\begin{Lemma}
Existiert ein 1-Retter oder ein \mbox{3-Retter} für eine Relation vom Grad drei, so lässt sich diese durch
Relationen vom Grad zwei ausdrücken.
\label{Lemma:Retterlemma}
\end{Lemma}
\begin{proof}
Die Aussage folgt trivial durch Einsetzen:
\begin{equation*}
A + B + C \gleich{(\ref{1-RetterA})} A + R + D \gleich{(\ref{1-RetterB})} D + E + F
\end{equation*}
\begin{equation*}
A + B + C \gleich{ (\ref{3-RetterA})} A + R_1 + R_2 \gleich{(\ref{3-RetterB}) } D + R_3 + R_2 \gleich{(\ref{3-RetterC}) } D + E + F
\end{equation*}
\vspace{-25pt}

$\,$
\end{proof}

Wir werden nun Retter für die Relationen vom Grad drei in den beiden problematischen Zellen suchen.

\paragraph{Die Zelle $\cell{ 1,1,2 } { 1,1,1,1 }$:}
Diese Zelle enthält die folgenden zwölf Gitterpunkte:

\newcommand{\ZelleEinsEins}{
\begin{dreiviermatrix}
1 & 0 & 0 & 0  \\
0 & 1 & 0 & 0  \\
0 & 0 & 1 & 1  \\
\end{dreiviermatrix}
}

\newcommand{\ZelleEinsZwo}{
\begin{dreiviermatrix}
1 & 0 & 0 & 0  \\
0 & 0 & 1 & 0  \\
0 & 1 & 0 & 1  \\
\end{dreiviermatrix}
}

\newcommand{\ZelleEinsDrei}{
\begin{dreiviermatrix}
1 & 0 & 0 & 0  \\
0 & 0 & 0 & 1  \\
0 & 1 & 1 & 0  \\
\end{dreiviermatrix}
}

\newcommand{\ZelleEinsVier}{
\begin{dreiviermatrix}
0 & 1 & 0 & 0  \\
1 & 0 & 0 & 0  \\
0 & 0 & 1 & 1  \\
\end{dreiviermatrix}
}

\newcommand{\ZelleEinsFuenf}{
\begin{dreiviermatrix}
0 & 1 & 0 & 0  \\
0 & 0 & 1 & 0  \\
1 & 0 & 0 & 1  \\
\end{dreiviermatrix}
}

\newcommand{\ZelleEinsSechs}{
\begin{dreiviermatrix}
0 & 1 & 0 & 0  \\
0 & 0 & 0 & 1  \\
1 & 0 & 1 & 0  \\
\end{dreiviermatrix}
}

\newcommand{\ZelleEinsSieben}{
\begin{dreiviermatrix}
0 & 0 & 1 & 0  \\
1 & 0 & 0 & 0  \\
0 & 1 & 0 & 1  \\
\end{dreiviermatrix}
}

\newcommand{\ZelleEinsAcht}{
\begin{dreiviermatrix}
0 & 0 & 1 & 0  \\
0 & 1 & 0 & 0  \\
1 & 0 & 0 & 1  \\
\end{dreiviermatrix}
}

\newcommand{\ZelleEinsNeun}{
\begin{dreiviermatrix}
0 & 0 & 1 & 0  \\
0 & 0 & 0 & 1  \\
1 & 1 & 0 & 0  \\
\end{dreiviermatrix}
}

\newcommand{\ZelleEinsZehn}{
\begin{dreiviermatrix}
0 & 0 & 0 & 1  \\
1 & 0 & 0 & 0  \\
0 & 1 & 1 & 0  \\
\end{dreiviermatrix}
}

\newcommand{\ZelleEinsElf}{
\begin{dreiviermatrix}
0 & 0 & 0 & 1  \\
0 & 1 & 0 & 0  \\
1 & 0 & 1 & 0  \\
\end{dreiviermatrix}
}

\newcommand{\ZelleEinsZwoelf}{
\begin{dreiviermatrix}
0 & 0 & 0 & 1  \\
0 & 0 & 1 & 0  \\
1 & 1 & 0 & 0  \\
\end{dreiviermatrix}
}

\begin{align*}
M_{1}&=\ZelleEinsEins
&
M_{2}&=\ZelleEinsZwo
&
M_{3}&=\ZelleEinsDrei
&
M_{4}&=\ZelleEinsVier
\\
M_{5}&=\ZelleEinsFuenf
&
M_{6}&=\ZelleEinsSechs
&
M_{7}&=\ZelleEinsSieben
&
M_{8}&=\ZelleEinsAcht
\\
M_{9}&=\ZelleEinsNeun
&
M_{10}&=\ZelleEinsZehn
&
M_{11}&=\ZelleEinsElf
&
M_{12}&=\ZelleEinsZwoelf
\end{align*}

Die von 4ti2 (\cite{4ti2}) berechneten Relationen vom Grad drei sind:
\begin{eqnarray}
M_{5} + M_{9} + M_{11} &=& M_{6} + M_{8} + M_{12}\\
M_{2} + M_{9} + M_{10} &=& M_{3} + M_{7} + M_{12}\\
M_{1} + M_{6} + M_{10} &=& M_{3} + M_{4} + M_{11}\\
M_{1} + M_{5} + M_{7} &=& M_{2} + M_{4} + M_{8}
\end{eqnarray}

\addtocounter{equation}{-4} %

\newcommand{\ZelleEinsA}{
\begin{dreiviermatrix}
\gc 0 & 1 & 0 & 0 \\ %
\gc 0 & 0 & 1 & 0 \\
\gc 1 & 0 & 0 & 1 
\end{dreiviermatrix}
}

\newcommand{\ZelleEinsB}{
\begin{dreiviermatrix}
\gc 0 & 0 & 1 & 0 \\ %
\gc 0 & 0 & 0 & 1 \\
\gc 1 & 1 & 0 & 0 \\
\end{dreiviermatrix}
}

\newcommand{\ZelleEinsC}{
\begin{dreiviermatrix}
\gc 0 & 0 & 0 & 1 \\ %
\gc 0 & 1 & 0 & 0 \\
\gc 1 & 0 & 1 & 0 \\
\end{dreiviermatrix}
}

\newcommand{\ZelleEinsD}{
\begin{dreiviermatrix}
\gc 0 & 1 & 0 & 0 \\ %
\gc 0 & 0 & 0 & 1 \\
\gc 1 & 0 & 1 & 0 \\
\end{dreiviermatrix}
}

\newcommand{\ZelleEinsE}{
\begin{dreiviermatrix}
\gc 0 & 0 & 1 & 0 \\ %
\gc 0 & 1 & 0 & 0 \\
\gc 1 & 0 & 0 & 1 \\
\end{dreiviermatrix}
}

\newcommand{\ZelleEinsF}{
\begin{dreiviermatrix}
\gc 0 & 0 & 0 & 1 \\ %
\gc 0 & 0 & 1 & 0 \\
\gc 1 & 1 & 0 & 0 \\
\end{dreiviermatrix}
}
\begin{equation}
\label{relation112_1111}
\ZelleEinsA
+
\ZelleEinsB
+
\ZelleEinsC
 =
\ZelleEinsD
+
\ZelleEinsE
+
\ZelleEinsF
\end{equation}
\begin{equation}
\begin{dreiviermatrix}
1 & \gc 0 & 0 & 0 \\
0 & \gc 0 & 1 & 0 \\
0 & \gc 1 & 0 & 1 \\
\end{dreiviermatrix}
+
\begin{dreiviermatrix}
0 & \gc 0 & 1 & 0 \\
0 & \gc 0 & 0 & 1 \\
1 & \gc 1 & 0 & 0 \\
\end{dreiviermatrix}
+
\begin{dreiviermatrix}
0 & \gc 0 & 0 & 1 \\
1 & \gc 0 & 0 & 0 \\
0 & \gc 1 & 1 & 0 \\
\end{dreiviermatrix}
 =
\begin{dreiviermatrix}
1 & \gc 0 & 0 & 0 \\
0 & \gc 0 & 0 & 1 \\
0 & \gc 1 & 1 & 0 \\
\end{dreiviermatrix}
+
\begin{dreiviermatrix}
0 & \gc 0 & 1 & 0 \\
1 & \gc 0 & 0 & 0 \\
0 & \gc 1 & 0 & 1 \\
\end{dreiviermatrix}
+
\begin{dreiviermatrix}
0 & \gc 0 & 0 & 1 \\
0 & \gc 0 & 1 & 0 \\
1 & \gc 1 & 0 & 0 \\
\end{dreiviermatrix}
\end{equation}
\begin{equation}
\begin{dreiviermatrix}
1 & 0 & \gc 0 & 0 \\
0 & 1 & \gc 0 & 0 \\
0 & 0 & \gc 1 & 1 \\
\end{dreiviermatrix}
+
\begin{dreiviermatrix}
0 & 1 & \gc 0 & 0 \\
0 & 0 & \gc 0 & 1 \\
1 & 0 & \gc 1 & 0 \\
\end{dreiviermatrix}
+
\begin{dreiviermatrix}
0 & 0 & \gc 0 & 1 \\
1 & 0 & \gc 0 & 0 \\
0 & 1 & \gc 1 & 0 \\
\end{dreiviermatrix}
 =
\begin{dreiviermatrix}
1 & 0 & \gc 0 & 0 \\
0 & 0 & \gc 0 & 1 \\
0 & 1 & \gc 1 & 0 \\
\end{dreiviermatrix}
+
\begin{dreiviermatrix}
0 & 1 & \gc 0 & 0 \\
1 & 0 & \gc 0 & 0 \\
0 & 0 & \gc 1 & 1 \\
\end{dreiviermatrix}
+
\begin{dreiviermatrix}
0 & 0 & \gc 0 & 1 \\
0 & 1 & \gc 0 & 0 \\
1 & 0 & \gc 1 & 0 \\
\end{dreiviermatrix}
\end{equation}
\begin{equation}
\begin{dreiviermatrix}
1 & 0 & 0 & \gc 0 \\
0 & 1 & 0 & \gc 0 \\
0 & 0 & 1 & \gc 1 \\
\end{dreiviermatrix}
+
\begin{dreiviermatrix}
0 & 1 & 0 & \gc 0 \\
0 & 0 & 1 & \gc 0 \\
1 & 0 & 0 & \gc 1 \\
\end{dreiviermatrix}
+
\begin{dreiviermatrix}
0 & 0 & 1 & \gc 0 \\
1 & 0 & 0 & \gc 0 \\
0 & 1 & 0 & \gc 1 \\
\end{dreiviermatrix}
 =
\begin{dreiviermatrix}
1 & 0 & 0 & \gc 0 \\
0 & 0 & 1 & \gc 0 \\
0 & 1 & 0 & \gc 1 \\
\end{dreiviermatrix}
+
\begin{dreiviermatrix}
0 & 1 & 0 & \gc 0 \\
1 & 0 & 0 & \gc 0 \\
0 & 0 & 1 & \gc 1 \\
\end{dreiviermatrix}
+
\begin{dreiviermatrix}
0 & 0 & 1 & \gc 0 \\
0 & 1 & 0 & \gc 0 \\
1 & 0 & 0 & \gc 1 \\
\end{dreiviermatrix}
\label{relation112_1111letzte}
\end{equation}

Wie man sieht, ist bei allen Relationen eine Spalte fest und von der Form $(0,0,1)^T$
(grau markiert). Projiziert man diese
Spalte weg, so handelt es sich bei allen vier Relationen um die Relation zwischen den Gitterpunkten im Birkhoffpolytop $B_3$.\footnote{Im 
Birkhoffpolytop $B_3$ gibt es genau eine Relation zwischen den Gitterpunkten: Die Summe der Ecken mit Determinante Eins ist gleich der Summe der Ecken mit Determinante minus Eins.}

Es genügt also, wenn wir nur noch die Relation (\ref{relation112_1111}) betrachten. 
Wir werden zeigen, dass stets ein 1-Retter oder ein 3-Retter für diese Relation existiert, wenn sie in einem glatten Transportpolytop auftritt.

\smallskip
Jede Zelle $Z$ grenzt an höchstens 24 andere Zellen.
Damit die Zelle an eine andere grenzt, muss nämlich \mbox{$[a_{ij} \ge 0]$} oder $[a_{ij} \le 1]$ eine Facette der Zelle sein, die gleichzeitig keine Facette des Polytops ist.
Dies wird in der folgenden Abbildung verdeutlicht:
\medskip

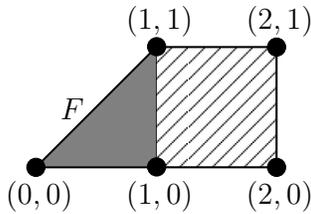
\begin{minipage}{4.5cm}
\vspace{1mm}
\scalebox{0.8}{
\input{ZelleMitNachbarn.pstex_t}
} 

\end{minipage}
\hfill
\begin{minipage}{9cm}
Betrachte das Polytop $P\subseteq \R^2$ in der Abbildung links. Es hat zwei volldimensionale Zellen: Die grau gefärbte Zelle $Z$ und die schraffierte Zelle $Z'$. 
$[x_1\le 1]$ definiert eine Facette von $Z$, aber keine Facette von $P$. \glqq Hinter\grqq\ dieser Facette von $Z$ liegt die Nachbarzelle $Z'$. 

Hinter der Facette $F=\{(x_1,x_2)\in P \,|\,$ $ x_1=x_2\}$ kann keine andere Zelle liegen.
 
 \end{minipage}
\smallskip

Wir überlegen uns nun, warum unsere Zelle keine Facette der Form \mbox{$[a_{ij} \le 1]$} hat.
Die volldimensionalen Zellen sind sechsdimensionale Gitterpolytope. In einer Facette müssen also mindestens sechs Gitterpunkte liegen, die einen fünfdimensionalen affinen Raum aufspannen.
$a_{ij}=1$ für $i\in\{1,2\}$ gilt jeweils nur für drei Gitterpunkte. $a_{3j}=1$ gilt jeweils für sechs Gitterpunkte. Diese sind aber mit den Relationen (\ref{relation112_1111})-(\ref{relation112_1111letzte}) affin abhängig, spannen also einen höchstens vierdimensionalen affinen Raum auf.
\smallskip

Demnach können Nachbarzellen nur hinter einer Facette der Form $[a_{ij}\ge 0]$ liegen. Diese zwölf Ungleichungen sind tatsächlich genau die Ungleichungen, die die zwölf Facetten der Zelle $\cell{ 1,1,2 } { 1,1,1,1 }$ definieren, wie eine Berechnung mit dem Programm Polymake (\cite{math.CO/0507273}) gezeigt hat.

\begin{Lemma}
\label{Lemma:NachbarnInGlattenPolytopen}
Sei $Z=\cell{1,1,2}{1,1,1,1}$ eine (verschobene) Zelle eines glatten ($3\times4$)-Trans\-port\-poly\-tops $\trans{r}{c}$. 

Dann gilt:
Mindestens drei der zwölf Facetten der Zelle sind nicht Facetten des Polytops, d.h. jede Zelle von diesem Typ ist zu mindestens drei anderen Zellen benachbart. 
\end{Lemma}

\begin{proof}

Angenommen mindestens zehn Facetten der Zelle sind auch Facetten des Polytops. Die Facetten haben alle die Form 
$[a_{ij}\ge 0]$. Die Matrizen in der Zelle enthalten jeweils acht Nulleinträge. 

Folglich kann man nun in jedem Fall eine Matrix $M\in Z$ finden, die in mindestens sieben Facetten von $Z$ enthalten ist, die auch gleichzeitig Facetten von $\trans{r}{c}$ sind.

Unser Polytop ist aber glatt und damit nach Satz \ref{Satz:KombinatorischGlatt} einfach. Also darf kein Punkt in mehr als $\dim(\trans{r}{c})=6$ Facetten enthalten sein.
\:\wid
$\,$
\end{proof}

\noindent Wir wollen nun mit Hilfe dieses Lemmas Retter finden. Dazu unterscheiden wir zwei Fälle:

\begin{Fallunterscheidung}
\Fall{
Es gibt eine Facette $[a_{ij} \ge 0]$  mit $i \in \{1,2,3\}, j\in \{2,3,4\}$, hinter der eine weitere Zelle $Z'$ liegt. \OBdA sei $i=1,\, j=2$. 
Dann ist 
\begin{equation*}
R=
\begin{dreiviermatrix}
\gc 0 & -1 & 1 & 1 \\  
\gc 0 &  1 & 0 & 0 \\
\gc 1 &  1 & 0 & 0 
\end{dreiviermatrix}\in Z'\,\footnote{Genau genommen liegt dieser Retter in einer verschobenen Version von $Z'$.}
\end{equation*}
ein 1-Retter für (\ref{relation112_1111}), denn:
\[
\underbrace{
\begin{dreiviermatrix}
\gc 0 & -1 & 1 & 1 \\  
\gc 0 &  1 & 0 & 0 \\
\gc 1 &  1 & 0 & 0 
\end{dreiviermatrix}
}_{R}
+
\underbrace{
\begin{dreiviermatrix}
\gc 0 &  1 & 0 & 0\\  
\gc 0 &  0 & 0 & 1 \\
\gc 1 &  0 & 1 & 0 
\end{dreiviermatrix}
}_{M_6}
=
\underbrace{
\begin{dreiviermatrix}
\gc 0 &  0 & 0 & 1\\  
\gc 0 &  1 & 0 & 0 \\
\gc 1 &  0 & 1 & 0 
\end{dreiviermatrix}
}_{M_{11}}
+
\underbrace{
\begin{dreiviermatrix}
\gc 0 &  0 & 1 & 0\\  
\gc 0 &  0 & 0 & 1 \\
\gc 1 &  1 & 0 & 0 
\end{dreiviermatrix}
}_{M_9}
\]
\[
\underbrace{
\begin{dreiviermatrix}
\gc 0 & -1 & 1 & 1 \\  
\gc 0 &  1 & 0 & 0 \\
\gc 1 &  1 & 0 & 0 
\end{dreiviermatrix}
}_{R}
+
\underbrace{
\begin{dreiviermatrix}
\gc 0 &  1 & 0 & 0\\  
\gc 0 &  0 & 1 & 0 \\
\gc 1 &  0 & 0 & 1 
\end{dreiviermatrix}
}_{M_5}
=
\underbrace{
\begin{dreiviermatrix}
\gc 0 &  0 & 1 & 0\\  
\gc 0 &  1 & 0 & 0 \\
\gc 1 &  0 & 0 & 1 
\end{dreiviermatrix}
}_{M_8}
+
\underbrace{
\begin{dreiviermatrix}
\gc 0 &  0 & 0 & 1\\  
\gc 0 &  0 & 1 & 0 \\
\gc 1 &  1 & 0 & 0 
\end{dreiviermatrix}
}_{M_{12}}
\]
}
\Fall{
Die drei Facetten, hinter denen eine benachbarte Zelle liegt, sind:
$[a_{11}\ge 0]$, $[a_{21}\ge 0]$ und $[a_{31}\ge 0]$.
Dann ist

\newcommand{\ZelleEinsDreiRetterA}{
\begin{dreiviermatrix} %
\gc -1 & 0 & 1 & 1 \\  
\gc 1 &  0 & 0 & 0 \\
\gc 1 &  1 & 0 & 0 
\end{dreiviermatrix}
}

\newcommand{\ZelleEinsDreiRetterB}{
\begin{dreiviermatrix} %
\gc 1 & 0 & 0 & 0 \\  
\gc -1 &  1 & 0 & 1 \\
\gc 1 &  0 & 1 & 0 
\end{dreiviermatrix}
}
\newcommand{\ZelleEinsDreiRetterC}{
\begin{dreiviermatrix} %
\gc -1 & 1 & 1 & 0 \\  
\gc 1 &  0 & 0 & 0 \\
\gc 1 &  0 & 0 & 1 
\end{dreiviermatrix}
}

\[
(R_1,R_2,R_3)=
\left(
\ZelleEinsDreiRetterA,
\ZelleEinsDreiRetterB,
\ZelleEinsDreiRetterC
\right)
\]
ein 3-Retter, denn:
\begin{align*}
\underbrace{\ZelleEinsB}_{M_9} + \underbrace{\ZelleEinsC}_{M_{11}} 
&= \underbrace{\ZelleEinsDreiRetterA}_{R_1}  + \underbrace{\ZelleEinsDreiRetterB}_{R_2} 
\displaybreak[2]
\\
\underbrace{\ZelleEinsA}_{M_5} + \underbrace{\ZelleEinsDreiRetterA}_{R_1}
&=  \underbrace{\ZelleEinsF}_{M_{12}} + \underbrace{\ZelleEinsDreiRetterC}_{R_3}
\displaybreak[2]
\\
\underbrace{\ZelleEinsDreiRetterB}_{R_2} + \underbrace{\ZelleEinsDreiRetterC}_{R_3} 
&= \underbrace{\ZelleEinsE}_{M_8}  + \underbrace{\ZelleEinsD}_{M_6}
\end{align*}

}
\end{Fallunterscheidung}

Wegen Lemma \ref{Lemma:Retterlemma}  
können wir also alle Relationen vom Grad drei, die zwischen Gitterpunkten von Zellen dieses Typs auftreten, durch Relationen vom Grad zwei von Gitterpunkten im ganzen Polytop darstellen.

\paragraph{Die Zelle $\cell{ 1,2,2 }{ 1,1,1,2 }$:}
Diese Zelle enthält die folgenden zwölf Gitterpunkte:

\begin{align*}
M_{1}&=\begin{dreiviermatrix}
1 & 0 & 0 & 0  \\
0 & 1 & 0 & 1  \\
0 & 0 & 1 & 1  \\
\end{dreiviermatrix}
&
M_{2}&=\begin{dreiviermatrix}
1 & 0 & 0 & 0  \\
0 & 0 & 1 & 1  \\
0 & 1 & 0 & 1  \\
\end{dreiviermatrix}
&
M_{3}&=\begin{dreiviermatrix}
0 & 1 & 0 & 0  \\
1 & 0 & 0 & 1  \\
0 & 0 & 1 & 1  \\
\end{dreiviermatrix}
&
M_{4}&=\begin{dreiviermatrix}
0 & 1 & 0 & 0  \\
0 & 0 & 1 & 1  \\
1 & 0 & 0 & 1  \\
\end{dreiviermatrix}
\\
M_{5}&=\begin{dreiviermatrix}
0 & 0 & 1 & 0  \\
1 & 0 & 0 & 1  \\
0 & 1 & 0 & 1  \\
\end{dreiviermatrix}
&
M_{6}&=\begin{dreiviermatrix}
0 & 0 & 1 & 0  \\
0 & 1 & 0 & 1  \\
1 & 0 & 0 & 1  \\
\end{dreiviermatrix}
&
M_{7}&=\begin{dreiviermatrix}
0 & 0 & 0 & 1  \\
1 & 1 & 0 & 0  \\
0 & 0 & 1 & 1  \\
\end{dreiviermatrix}
&
M_{8}&=\begin{dreiviermatrix}
0 & 0 & 0 & 1  \\
1 & 0 & 1 & 0  \\
0 & 1 & 0 & 1  \\
\end{dreiviermatrix}
\\
M_{9}&=\begin{dreiviermatrix}
0 & 0 & 0 & 1  \\
1 & 0 & 0 & 1  \\
0 & 1 & 1 & 0  \\
\end{dreiviermatrix}
&
M_{10}&=\begin{dreiviermatrix}
0 & 0 & 0 & 1  \\
0 & 1 & 1 & 0  \\
1 & 0 & 0 & 1  \\
\end{dreiviermatrix}
&
M_{11}&=\begin{dreiviermatrix}
0 & 0 & 0 & 1  \\
0 & 1 & 0 & 1  \\
1 & 0 & 1 & 0  \\
\end{dreiviermatrix}
&
M_{12}&=\begin{dreiviermatrix}
0 & 0 & 0 & 1  \\
0 & 0 & 1 & 1  \\
1 & 1 & 0 & 0  \\
\end{dreiviermatrix}
\end{align*}

Die von 4ti2 (\cite{4ti2}) berechnete kubische Relation ist:

\newcommand{\ZelleZwoA}{
\begin{dreiviermatrix}
1 & 0 & 0 & \gc 0 \\ %
0 & 1 & 0 & \gc 1 \\
0 & 0 & 1 & \gc 1 \\
\end{dreiviermatrix}
}
\newcommand{\ZelleZwoB}{
\begin{dreiviermatrix}
0 & 1 & 0 & \gc 0 \\ %
0 & 0 & 1 & \gc 1 \\
1 & 0 & 0 & \gc 1 \\
\end{dreiviermatrix}
}
\newcommand{\ZelleZwoC}{
\begin{dreiviermatrix}
0 & 0 & 1 & \gc 0 \\ %
1 & 0 & 0 & \gc 1 \\
0 & 1 & 0 & \gc 1 \\
\end{dreiviermatrix}
}
\newcommand{\ZelleZwoD}{
\begin{dreiviermatrix}
1 & 0 & 0 & \gc 0 \\ %
0 & 0 & 1 & \gc 1 \\
0 & 1 & 0 & \gc 1 \\
\end{dreiviermatrix}
}
\newcommand{\ZelleZwoE}{
\begin{dreiviermatrix}
0 & 1 & 0 & \gc 0 \\ %
1 & 0 & 0 & \gc 1 \\
0 & 0 & 1 & \gc 1 \\
\end{dreiviermatrix}
}
 \newcommand{\ZelleZwoF}{
\begin{dreiviermatrix}
0 & 0 & 1 & \gc 0 \\ %
0 & 1 & 0 & \gc 1 \\
1 & 0 & 0 & \gc 1 \\
\end{dreiviermatrix}
}

\addtocounter{equation}{-1}

\begin{align}\nonumber
M_{1} + M_{4} + M_{5} &= M_{2} + M_{3} + M_{6} 
\displaybreak[1]\\
\ZelleZwoA
+ 
\ZelleZwoB
+
\ZelleZwoC
&=
\ZelleZwoD
+ 
\ZelleZwoE
+
\ZelleZwoF
\end{align}

Wie man sieht, ist die vierte Spalte konstant $(0,1,1)^T$ und die Relation ist nach Projektion wieder die Relation aus dem Birkhoffpolytop $B_3$.

Eine Berechnung mit Polymake (\cite{math.CO/0507273}) liefert, dass die Zelle die folgenden elf Facetten hat:
$[a_{ij} \ge 0]$ für $i \in \{1,2,3\}$ und $j\in \{1,2,3\}$, sowie $[a_{24} \le 1]$ und $[a_{34} \le 1]$.

\smallskip
\begin{minipage}{2cm}
$\begin{dreiviermatrix} 
0 & 0 & 0 & x \\  
0  & 0 & 0 & 1 \\
0  & 0 & 0 & 1 
\end{dreiviermatrix}$ 
\end{minipage}
\hfill
\begin{minipage}{11.3cm}
Die links stehende Matrix liegt also für beliebiges $x$ im Schnitt aller elf oben angegebenen facettendefinierenden Hyperebenen.
\end{minipage}
\smallskip

Wir wissen, dass hinter einer der Facetten der Zelle noch eine weitere Zelle liegen muss, da die  Zelle in einem glatten Polytop liegt.
Wir unterscheiden wieder zwei Fälle:

\begin{Fallunterscheidung}
\Fall{Die Zelle hat eine Facette
$[a_{ij} \ge 0]$, hinter der eine Nachbarzelle $Z'$ liegt mit $i \in \{1,2,3\}$ und $j\in \{1,2,3\}$.   \OBdA sei $i=j=1$.

Dann ist 

\newcommand{\ZelleZwoEinsRetter}{ %
\begin{dreiviermatrix} %
-1 & 1 & 1 & \gc 0 \\  
1  & 0 & 0 & \gc 1 \\
1  & 0 & 0 & \gc 1 
\end{dreiviermatrix}
}
\[
R=\ZelleZwoEinsRetter\in Z'
\] ein 1-Retter, denn:

\begin{align*}
\underbrace{\ZelleZwoEinsRetter}_{R}  + \underbrace{\ZelleZwoD}_{M_2} &=  \underbrace{\ZelleZwoB}_{M_4} + \underbrace{\ZelleZwoC}_{M_5} 
\\
\underbrace{\ZelleZwoEinsRetter}_{R}  + \underbrace{\ZelleZwoA}_{M_1} &=  
\underbrace{\ZelleZwoF}_{M_6} + \underbrace{\ZelleZwoE}_{M_3}
\end{align*}
}

\Fall{
Liegt hinter keiner der neun Facetten der Zelle, die wir im ersten Fall betrachtet haben, eine Nachbarzelle, so müssen hinter den durch $[a_{24} \le 1]$ und 
 $\mbox{$[a_{34} \le 1]$}$ definierten Facetten Nachbarzellen liegen.

Wäre dies nicht so, so läge beispielsweise die Matrix $M_2$ in sieben Facetten des Polytops, und wir erhielten wie im Beweis von Lemma \ref{Lemma:NachbarnInGlattenPolytopen}
einen Widerspruch.

Es ist dann

\newcommand{\ZelleZwoDreiRetterA}{
\begin{dreiviermatrix}
0 & 0 & 1 & \gc 0 \\  
0 & 0 & 0 & \gc 2 \\
1 & 1 & 0 & \gc 0 
\end{dreiviermatrix}
}
\newcommand{\ZelleZwoDreiRetterB}{
\begin{dreiviermatrix}
0 & 1 & 0 & \gc 0 \\  
1 & 0 & 1 & \gc 0 \\
0 & 0 & 0 & \gc 2 
\end{dreiviermatrix}
}
\newcommand{\ZelleZwoDreiRetterC}{
\begin{dreiviermatrix}
1 & 0 & 0 & \gc 0 \\  
0 & 0 & 0 & \gc 2 \\
0 & 1 & 1 & \gc 0 
\end{dreiviermatrix}
}

\[
(R_1,R_2,R_3)=
\left(
\ZelleZwoDreiRetterA,
\ZelleZwoDreiRetterB,
\ZelleZwoDreiRetterC
\right)
\]
ein 3-Retter, denn:

\begin{align*}
\underbrace{\ZelleZwoB}_{M_4} + \underbrace{\ZelleZwoC}_{M_5} &= \underbrace{\ZelleZwoDreiRetterA}_{R_1}  + \underbrace{\ZelleZwoDreiRetterB}_{R_2}
\\
\underbrace{\ZelleZwoA}_{M_1} + \underbrace{\ZelleZwoDreiRetterA}_{R_1} &= 
\underbrace{\ZelleZwoF}_{M_6} + \underbrace{\ZelleZwoDreiRetterC}_{R_3} 
\\
\underbrace{\ZelleZwoDreiRetterB}_{R_2} + \underbrace{\ZelleZwoDreiRetterC}_{R_3} &= \underbrace{\ZelleZwoE}_{M_3}  + \underbrace{\ZelleZwoD}_{M_2}
\end{align*}

}

\end{Fallunterscheidung}

Wir können also auch hier wieder aufgrund
von Lemma \ref{Lemma:Retterlemma} alle Relationen  vom Grad drei zwischen Gitterpunkten von Zellen dieses Typs durch Relationen vom Grad zwei von Gitterpunkten im ganzen Polytop darstellen.
\smallskip

Damit ist  der folgende Satz bewiesen:
\begin{Satz}
Sei $\trans{r}{c}$ ein glattes ($3\times 4 $)-Transportpolytop.
Dann ist das torische Ideal  $I_{\trans{r}{c}}$  im Grad zwei erzeugt.
\end{Satz}

\subsection*{Beispiele}
\begin{Beispiel}[Ein glattes Transportpolytop]
Betrachte nochmal das Transportpolytop $\transv{(1,1,1,10)}{(3,3,3,3)}\cong \transv{(1,1,1,6)}{(2,2,2,2)}$
aus den Beispielen \ref{Beispiel:IsoTPglattKomb} (S.~\pageref{Beispiel:IsoTPglattKomb}) und \ref{Beispiel1110_3333} (S.~\pageref{Beispiel1110_3333}).
Das zugehörige torische Ideal ist tatsächlich im Grad zwei erzeugt, nämlich von den folgenden Binomen, wie eine Berechnung mit 4ti2 (\cite{4ti2}) ergeben hat:
\begin{tabbing}
$x_{11}x_{16} - x_{12} x_{15} $,
\= $x_{10}x_{15} - x_{11} x_{14} $,
\= $x_{10}x_{16} - x_{12} x_{14} $,
\= $x_{9}x_{14} - x_{10} x_{13} $,
\= $x_{9}x_{15} - x_{11} x_{13} $, \\
$x_{9}x_{16} - x_{12} x_{13} $,
\> $x_{7}x_{12} - x_{8} x_{11} $,
\> $x_{7}x_{16} - x_{8} x_{15} $,
\> $x_{6}x_{11} - x_{7} x_{10} $,
\> $x_{6}x_{15} - x_{7} x_{14} $, \\
$x_{6}x_{12} - x_{8} x_{10} $,
\> $x_{6}x_{16} - x_{8} x_{14} $,
\> $x_{5}x_{10} - x_{6} x_{9} $,
\> $x_{5}x_{14} - x_{6} x_{13} $,
\> $x_{5}x_{11} - x_{7} x_{9} $, \\
$x_{5}x_{15} - x_{7} x_{13} $,
\> $x_{5}x_{12} - x_{8} x_{9} $,
\> $x_{5}x_{16} - x_{8} x_{13} $,
\> $x_{3}x_{8} - x_{4} x_{7} $,
\> $x_{3}x_{12} - x_{4} x_{11} $, \\
$x_{3}x_{16} - x_{4} x_{15} $,
\> $x_{2}x_{7} - x_{3} x_{6} $,
\> $x_{2}x_{11} - x_{3} x_{10} $,
\> $x_{2}x_{15} - x_{3} x_{14} $,
\> $x_{2}x_{8} - x_{4} x_{6} $, \\
$x_{2}x_{12} - x_{4} x_{10} $,
\> $x_{2}x_{16} - x_{4} x_{14} $,
\> $x_{1}x_{6} - x_{2} x_{5} $,
\> $x_{1}x_{10} - x_{2} x_{9} $,
\> $x_{1}x_{14} - x_{2} x_{13} $, \\
$x_{1}x_{7} - x_{3} x_{5} $,
\> $x_{1}x_{11} - x_{3} x_{9} $,
\> $x_{1}x_{15} - x_{3} x_{13} $,
\> $x_{1}x_{8} - x_{4} x_{5} $,
\> $x_{1}x_{12} - x_{4} x_{9} $, \\
$x_{1}x_{16} - x_{4} x_{13} $
\end{tabbing}
\end{Beispiel}

\begin{Beispiel}[Ein Transportpolytop, das nicht glatt ist]
$\,$\\
Torische Ideale von Transportpolytopen, die nicht glatt sind, sind i.\,A. nicht im Grad zwei erzeugt. Betrachte z.\,B. das Polytop $\transv{(113)}{(1112)}$, das nach Satz
\ref{Satz:KombinatorischGlatt} nicht glatt ist.

Es gilt:
\newcommand{\einrueckenXXX}{\hspace{10pt}} %
\newcommand{\einrueckenXXXX}{\hspace{9pt}} %

\begin{align*}
\A:=\transv{(113)}{(1112)} \cap \Z^{3\times 4}  =  &\left\{
\begin{dreiviermatrix}
1 & 0 & 0 & 0  \\
0 & 1 & 0 & 0  \\
0 & 0 & 1 & 2  \\
\end{dreiviermatrix}
,\;
\begin{dreiviermatrix}
1 & 0 & 0 & 0  \\
0 & 0 & 1 & 0  \\
0 & 1 & 0 & 2  \\
\end{dreiviermatrix}
,\;
\begin{dreiviermatrix}
1 & 0 & 0 & 0  \\
0 & 0 & 0 & 1  \\
0 & 1 & 1 & 1  \\
\end{dreiviermatrix}
,\;
\begin{dreiviermatrix}
0 & 1 & 0 & 0  \\
1 & 0 & 0 & 0  \\
0 & 0 & 1 & 2  \\
\end{dreiviermatrix}
,\;
\right.
\\
&
\einrueckenXXX
\begin{dreiviermatrix}
0 & 1 & 0 & 0  \\
0 & 0 & 1 & 0  \\
1 & 0 & 0 & 2  \\
\end{dreiviermatrix}
,\;
\begin{dreiviermatrix}
0 & 1 & 0 & 0  \\
0 & 0 & 0 & 1  \\
1 & 0 & 1 & 1  \\
\end{dreiviermatrix}
,\;
\begin{dreiviermatrix}
0 & 0 & 1 & 0  \\
1 & 0 & 0 & 0  \\
0 & 1 & 0 & 2  \\
\end{dreiviermatrix}
,\;
\begin{dreiviermatrix}
0 & 0 & 1 & 0  \\
0 & 1 & 0 & 0  \\
1 & 0 & 0 & 2  \\
\end{dreiviermatrix}
,\;
\\
&
\einrueckenXXX
\begin{dreiviermatrix}
0 & 0 & 1 & 0  \\
0 & 0 & 0 & 1  \\
1 & 1 & 0 & 1  \\
\end{dreiviermatrix}
,\;
\begin{dreiviermatrix}
0 & 0 & 0 & 1  \\
1 & 0 & 0 & 0  \\
0 & 1 & 1 & 1  \\
\end{dreiviermatrix}
,\;
\begin{dreiviermatrix}
0 & 0 & 0 & 1  \\
0 & 1 & 0 & 0  \\
1 & 0 & 1 & 1  \\
\end{dreiviermatrix}
,\;
\begin{dreiviermatrix}
0 & 0 & 0 & 1  \\
0 & 0 & 1 & 0  \\
1 & 1 & 0 & 1  \\
\end{dreiviermatrix}
,\;
\\
&
\einrueckenXXXX
\left.
\begin{dreiviermatrix}
0 & 0 & 0 & 1  \\
0 & 0 & 0 & 1  \\
1 & 1 & 1 & 0  \\
\end{dreiviermatrix}
\right
\}
\end{align*}

$\transv{(113)}{(1112)}$
hat lediglich zwei volldimensionale Zellen, nämlich  $\cell{1,1,3}{1,1,1,2}$ (ohne Verschiebung) und
$\cell{1,1,2}{1,1,1,1}$ (nach Verschiebung durch $a'_{34}:=a_{34}-1$).

Da das Polytop nicht glatt ist, kann man Lemma \ref{Lemma:NachbarnInGlattenPolytopen} nicht anwenden. Die Zelle $\cell{1,1,2}{1,1,1,1}$ hat nur eine Nachbarzelle. Diese liegt hinter der Facette $[a_{34}\ge 0]$.

Die Relation 
\begin{equation}
\label{equation:Grad3RelationUnrettbar}
\begin{dreiviermatrix}
1 & 0 & 0 & 0 \\
0 & 1 & 0 & 0 \\
0 & 0 & 1 & 2 \\
\end{dreiviermatrix}
+
\begin{dreiviermatrix}
0 & 1 & 0 & 0 \\
0 & 0 & 1 & 0 \\
1 & 0 & 0 & 2 \\
\end{dreiviermatrix}
+
\begin{dreiviermatrix}
0 & 0 & 1 & 0 \\
1 & 0 & 0 & 0 \\
0 & 1 & 0 & 2 \\
\end{dreiviermatrix}
 =
\begin{dreiviermatrix}
1 & 0 & 0 & 0 \\
0 & 0 & 1 & 0 \\
0 & 1 & 0 & 2 \\
\end{dreiviermatrix}
+
\begin{dreiviermatrix}
0 & 1 & 0 & 0 \\
1 & 0 & 0 & 0 \\
0 & 0 & 1 & 2 \\
\end{dreiviermatrix}
+
\begin{dreiviermatrix}
0 & 0 & 1 & 0 \\
0 & 1 & 0 & 0 \\
1 & 0 & 0 & 2 \\
\end{dreiviermatrix}
\end{equation}
lässt sich nicht durch Relationen vom Grad zwei ausdrücken.
Alle in der Relation vorkommenden Punkte liegen nämlich in der Hyperebene 
$H=\{A\in\R^{3\times 4}\,|\, a_{34}=2 \}$. Da diese Hyperebene eine Seite des Polytops definiert, müsste jeder Retter für diese Relation auch in dieser Hyperebene liegen
(Argumentation genauso wie im Beweis von Satz \ref{Satz:SeiteTorischesIdeal}).
Dies ist aber nicht möglich, da im Polytop keine weiteren Gitterpunkte enthalten sind, die in dieser Hyperebene liegen.

Also hat jedes minimale Erzeugendensystem des torischen Ideals zu diesem Transportpolytop Grad drei. 
Ein mittels 4ti2 (\cite{4ti2}) berechnetes minimales Erzeugendensystem ist:
\newlength{\TabBreite}
\setlength{\TabBreite}{2.6cm}
\begin{tabbing}
 \hspace*{3cm}
 \=\hspace*{\TabBreite}
 \=\hspace*{\TabBreite}
 \=\hspace*{\TabBreite}
 \=\hspace*{\TabBreite}
\kill
$x_{1}x_{5}x_{7} - x_{2} x_{4} x_{8}$,
\>$x_{7}x_{11} - x_{8} x_{10}$,
\>$x_{7}x_{13} - x_{9} x_{10}$,
\>$x_{5}x_{13} - x_{6} x_{12}$,
\>$x_{4}x_{12} - x_{5} x_{10}$, \\
$x_{4}x_{9} - x_{6} x_{7}$,
\>$x_{4}x_{13} - x_{6} x_{10}$,
\>$x_{2}x_{6} - x_{3} x_{5}$,
\>$x_{2}x_{13} - x_{3} x_{12}$,
\>$x_{8}x_{13} - x_{9} x_{11}$, \\
$x_{1}x_{12} - x_{2} x_{11}$,
\>$x_{1}x_{9} - x_{3} x_{8}$,
\>$x_{1}x_{13} - x_{3} x_{11}$
\end{tabbing}
Das Binom $x_{1}x_{5}x_{7} - x_{2} x_{4} x_{8}$ entspricht der Relation (\ref{equation:Grad3RelationUnrettbar}).
\end{Beispiel}

\section{Ausblick}
Wir haben in diesem Kapitel einige Antworten auf die in der Einführung aufgeworfenen Fragen gefunden. Es bleiben aber weiterhin  offene Fragen.

Die Frage nach dem Grad von minimalen Erzeugendensystemen haben wir vollständig beantwortet: Wir haben gezeigt, dass die torischen Ideale von \emph{allen} Flusspolytopen im Grad drei erzeugt sind und wir kennen Beispiele, für die diese Schranke scharf ist. 

Bei der Frage nach Gradschranken für Gröbnerbasen haben wir uns auf die torischen Ideale von Transportpolytopen beschränkt. Wir haben gezeigt, dass die reduzierten Gröbnerbasen von torischen Idealen von $(m\times n)$-Transportpolytopen bezüglich einer beliebigen umgekehrt lexikographischen Termordnung höchstens Grad $\left\lfloor\frac{m\cdot n}{2}\right\rfloor$ haben.

Eine interessante Frage ist, ob sich ein ähnliches Ergebnis auch für Flusspolytope erzielen lässt.
Unser Beweis lässt sich leider nicht direkt übertragen, da wir die Tatsache benutzen, dass alle ganzzahligen Matrizen in einer Zelle eines Transportpolytops eine konstante Anzahl Einsen enthalten. Eine analoge Aussage für Flusspolytope gilt leider nicht.

Eine weitere  Frage ist, ob diese Schranke scharf ist.
Wir konnten dafür keine Beispiele finden. Alle reduzierten Gröbnerbasen von ($m\times n$)-Zellen, die wir am Computer ausgerechnet haben, hatten höchstens Grad  $\min(m,n)$.
Wir haben allerdings gezeigt, dass wir mehr Mühe auf die Wahl der Termordnung verwenden müssen, wenn wir eine bessere Schranke beweisen wollen. Es gibt nämlich ungünstige Termordnungen, für die unsere obere Schranke annähernd scharf ist, sogar dann, wenn man sich auf glatte Transportpolytope beschränkt.

Wir haben gezeigt, dass torische Ideale von glatten $(3\times 4)$-Transportpolytopen  im Grad zwei erzeugt sind. Als nächstes könnte man darüber nachdenken, ob sie auch Gröbnerbasen im Grad zwei haben. Gegenbeispiele sind uns nicht bekannt. 

Es wäre auch interessant sich zu überlegen, ob man dieses Ergebnis auf Transportpolytope größerer Dimension oder sogar Flusspolytope verallgemeinern kann.
Dazu ist es aber notwendig, sich eine neue Beweistechnik einfallen zu lassen. 
Denn einerseits sind sicher sehr viele Fallunterscheidungen notwendig, wenn man für alle 
problematischen Relationen von Zellen größerer Transportpolytope Retter finden will.
Andererseits stößt man auch recht schnell auf große Hürden, wenn man mit Hilfe eines Computers Erzeugendensysteme und Gröbnerbasen der Zellen berechnen will.
Ein aktueller PC mit aktueller Software (4ti2) ist beispielsweise nicht mehr in der Lage minimale Erzeugendensysteme für die torischen Ideale von allen Zellen von ($6\times 6$)-Transportpolytopen auszurechnen.

\cleardoublepage

\bibliographystyle{alpha}
\bibliography{../Allgemein/literatur}

\clearpage %
\setlength{\parskip}{0pt} %

\noindent\thispagestyle{empty}{\bf Eidesstattliche Erkl\"arung:}\\
Hiermit erkl\"are ich an Eides statt, die vorliegende Arbeit eigenh\"andig
angefertigt zu haben. Alle verwendeten Hilfsmittel sind aufgef\"uhrt. Des
Weiteren versichere ich, dass ich diese Arbeit nicht in dieser oder \"ahnlicher
Form an einer anderen Universit\"at im Rahmen eines Pr\"ufungsverfahrens
eingereicht habe. 

\vspace{1.5cm}
\noindent Berlin, den \DatumErklaerung

\vspace{1.5cm}
\noindent Matthias Lenz

\end{document}